\documentclass[12pt]{article}

\usepackage[utf8]{inputenc}
\usepackage[T1]{fontenc}

\usepackage[a4paper, hmargin=2cm, vmargin=2cm]{geometry}
\usepackage{ragged2e}

\usepackage[pagebackref, hidelinks]{hyperref}
\hypersetup{
colorlinks = true, 
citecolor = cyan, 
linkcolor = blue, 
urlcolor = magenta, 
anchorcolor = purple
}

\usepackage{amsmath, amsthm, amsfonts}
\usepackage{amscd, amssymb, latexsym}

\usepackage[all]{xy}
\usepackage{mathtools} 
\usepackage{complexity}

\usepackage{comment, xcolor, soul}
\usepackage{enumitem}
\usepackage{verbatim}

\usepackage{subfig}
\usepackage{float}
\usepackage{multicol}
\usepackage{import}

\usepackage{tikz}
\usepackage{scalefnt, anyfontsize}

\usepackage{adjustbox} 

\newtheorem{definition}{Definition}[section]

\newtheorem{theorem}[definition]{Theorem}
\newtheorem{algorithm}[definition]{Algorithm}
\newtheorem{lemma}[definition]{Lemma}
\newtheorem{proposition}[definition]{Proposition}
\newtheorem{corollary}[definition]{Corollary}
\newtheorem{question}[definition]{Question}

\newtheorem{remark}[definition]{Remark}
\newtheorem{example}[definition]{Example}

\newcommand{\Id}{1}

\newcommand{\N}{\mathbb{N}}
\newcommand{\Z}{\mathbb{Z}}
\renewcommand{\R}{\mathbb{R}}

\newcommand{\Free}{\mathcal{F}}
\newcommand{\Tree}{\mathcal{T}}
\newcommand{\Surf}{\mathbb{F}}
\newcommand{\Sph}{\mathbb{S}}
\newcommand{\Regions}{\ensuremath{R}}

\DeclareMathOperator{\Disc}{\mathbb{D}}
\DeclareMathOperator{\Ball}{\mathbb{B}}

\DeclareMathOperator{\im}{im}

\newcommand{\true}{1} 
\newcommand{\false}{0} 
\DeclareMathOperator{\MoB}{MoB}

\DeclareMathOperator{\Parts}{\mathcal{P}}
\DeclareMathOperator{\Card}{Card}

\DeclareMathOperator{\len}{len}
\DeclareMathOperator{\si}{si}
\DeclareMathOperator{\ti}{i}
\DeclareMathOperator{\lk}{lk}

\DeclareMathOperator{\cord}{cord}
\DeclareMathOperator{\cross}{cross}
\DeclareMathOperator{\val}{val}

\title{The pinning ideal of a multiloop}

\author{Christopher-Lloyd Simon and Ben Stucky}
\date{\today}

\begin{document}

\maketitle

\begin{abstract}
A multiloop $\gamma\colon \sqcup_1^s \mathbb{S}^1 \looparrowright \mathbb{F}$ is a generic immersion of a finite union of circles into an oriented surface, considered up to homeomorphisms. 
A pinning set is a set of points $P\subset \mathbb{F}\setminus \im(\gamma)$, such that in the punctured surface $\mathbb{F} \setminus P$, the immersion $\gamma$ has the minimal number of double points in its homotopy class.
The collection of pinning sets of $\gamma$ forms a poset under inclusion called the pinning ideal $\mathcal{PI}(\gamma)$ which is endowed with the cardinal function whose minimum defines the pinning number $\varpi(\gamma)$.

We show that the decision problem associated to computing the pinning number of a multiloop is \textsf{NP}-complete, even for loops in the sphere. 
We give two proofs that it is \textsf{NP}: First, we implement a polynomial algorithm to check if a point-set is pinning, adapting methods of Birman--Series and Cohen--Lustig for computing intersection numbers of curves in surfaces.
Second, for loops in the sphere we reduce the problem in polynomial time to a variant of boolean satisfiability by applying a theorem of Hass--Scott characterizing taut loops, and adapting algorithms of Blank and Shor--Van Wyk which decide when a curve in the plane bounds an immersed disc. To show that it is \textsf{NP}-hard we reduce the vertex cover problem for graphs to the pinning problem for plane loops.

We use our algorithms to compute the pinning ideals for $\approx 1000$ of the smallest multiloops in the sphere, available in the online \href{https://christopherlloyd.github.io/LooPindex/}{LooPindex} catalog.
\end{abstract}

\begin{figure}[H]
\centering
\vspace{-1.5cm}
{\scalefont{1.5}\rotatebox{-45}{\scalebox{0.4}{\definecolor{c990000}{RGB}{153,0,0}

\def \globalscale {1.000000}
\begin{tikzpicture}[y=1cm, x=1cm, yscale=\globalscale,xscale=\globalscale, every node/.append style={scale=\globalscale}, inner sep=0pt, outer sep=0pt]
  \path[draw=black,line cap=round,line width=0.1588cm] (3.1168, 5.1911) -- (3.1168, 0.5345);

  \path[draw=black,line cap=round,line width=0.1588cm] (3.1168, 0.5345) -- (5.4451, 0.5345);

  \path[draw=black,line cap=round,line width=0.1588cm] (5.4451, 0.5345) -- (5.4451, 12.1761);

  \path[draw=black,line cap=round,line width=0.1588cm] (5.4451, 12.1761) -- (12.4301, 12.1761);

  \path[draw=black,line cap=round,line width=0.1588cm] (12.4301, 12.1761) -- (12.4301, 7.5195);

  \path[draw=black,line cap=round,line width=0.1588cm] (12.4301, 7.5195) -- (0.7885, 7.5195);

  \path[draw=black,line cap=round,line width=0.1588cm] (0.7885, 7.5195) -- (0.7885, 2.8628);

  \path[draw=black,line cap=round,line width=0.1588cm] (0.7885, 2.8628) -- (7.7735, 2.8628);

  \path[draw=black,line cap=round,line width=0.1588cm] (7.7735, 2.8628) -- (7.7735, 9.8478);

  \path[draw=black,line cap=round,line width=0.1588cm] (7.7735, 9.8478) -- (10.1018, 9.8478);

  \path[draw=black,line cap=round,line width=0.1588cm] (10.1018, 9.8478) -- (10.1018, 5.1911);

  \path[draw=black,line cap=round,line width=0.1588cm] (10.1018, 5.1911) -- (3.1168, 5.1911);

  \path[draw=red,fill=red,line width=0.0cm] (1.3123, 6.313) ellipse (0.291cm and 0.291cm);

  \tikzstyle{every node}=[font=\fontsize{18}{18}\selectfont]
  \node[text=white,anchor=south,line width=0.0cm,rotate=45.0] (text24) at (1.4291, 6.1292){$\mathbf{A}$};

  \path[draw=c990000,fill=c990000,line width=0.0cm] (3.6407, 3.9846) ellipse (0.291cm and 0.291cm);

  \node[text=white,anchor=south,line width=0.0cm,rotate=45.0] (text26) at (3.7593, 3.7961){$\mathbf{B}$};

  \path[draw=red,fill=red,line width=0.0cm] (3.6407, 1.6563) ellipse (0.291cm and 0.291cm);

  \node[text=white,anchor=south,line width=0.0cm,rotate=45.0] (text28) at (3.7574, 1.4726){$\mathbf{A}$};

  \path[draw=c990000,fill=c990000,line width=0.0cm] (4.281, 1.6563) circle (0.291cm);

  \node[text=white,anchor=south,line width=0.0cm,rotate=45.0] (text29) at (4.3996, 1.4678){$\mathbf{B}$};

  \path[draw=c990000,fill=c990000,line width=0.0cm] (5.969, 6.313) ellipse (0.291cm and 0.291cm);

  \node[text=white,anchor=south,line width=0.0cm,rotate=45.0] (text33) at (6.0877, 6.1245){$\mathbf{B}$};

  \path[draw=red,fill=red,line width=0.0cm] (5.969, 3.9846) circle (0.291cm);

  \node[text=white,anchor=south,line width=0.0cm,rotate=45.0] (text35) at (6.0857, 3.8009){$\mathbf{A}$};

  \path[draw=red,fill=red,line width=0.0cm] (8.2973, 8.6413) circle (0.291cm);

  \node[text=white,anchor=south,line width=0.0cm,rotate=45.0] (text37) at (8.4141, 8.4576){$\mathbf{A}$};

  \path[draw=c990000,fill=c990000,line width=0.0cm] (8.9376, 8.6413) circle (0.291cm);

  \node[text=white,anchor=south,line width=0.0cm,rotate=45.0] (text38) at (9.0563, 8.4528){$\mathbf{B}$};

\end{tikzpicture}}}}
\qquad
{\scalefont{1.7}\rotatebox{-45}{\scalebox{0.4}{\definecolor{c007700}{RGB}{0,119,0}
\definecolor{lime}{RGB}{0,255,0}
\definecolor{c00bb00}{RGB}{0,187,0}

\def \globalscale {1.000000}
\begin{tikzpicture}[y=1cm, x=1cm, yscale=\globalscale,xscale=\globalscale, every node/.append style={scale=\globalscale}, inner sep=0pt, outer sep=0pt]
  \path[draw=black,line cap=round,line width=0.1588cm] (3.1168, 5.1911) -- (3.1168, 0.5345);

  \path[draw=black,line cap=round,line width=0.1588cm] (3.1168, 0.5345) -- (5.4451, 0.5345);

  \path[draw=black,line cap=round,line width=0.1588cm] (5.4451, 0.5345) -- (5.4451, 12.1761);

  \path[draw=black,line cap=round,line width=0.1588cm] (5.4451, 12.1761) -- (12.4301, 12.1761);

  \path[draw=black,line cap=round,line width=0.1588cm] (12.4301, 12.1761) -- (12.4301, 7.5195);

  \path[draw=black,line cap=round,line width=0.1588cm] (12.4301, 7.5195) -- (0.7885, 7.5195);

  \path[draw=black,line cap=round,line width=0.1588cm] (0.7885, 7.5195) -- (0.7885, 2.8628);

  \path[draw=black,line cap=round,line width=0.1588cm] (0.7885, 2.8628) -- (7.7735, 2.8628);

  \path[draw=black,line cap=round,line width=0.1588cm] (7.7735, 2.8628) -- (7.7735, 9.8478);

  \path[draw=black,line cap=round,line width=0.1588cm] (7.7735, 9.8478) -- (10.1018, 9.8478);

  \path[draw=black,line cap=round,line width=0.1588cm] (10.1018, 9.8478) -- (10.1018, 5.1911);

  \path[draw=black,line cap=round,line width=0.1588cm] (10.1018, 5.1911) -- (3.1168, 5.1911);

  \path[draw=c007700,fill=c007700,line width=0.0cm] (1.3123, 6.313) ellipse (0.291cm and 0.291cm);

  \tikzstyle{every node}=[font=\fontsize{18}{18}\selectfont]
  \node[text=white,anchor=south,line width=0.0cm,rotate=45.0] (text24) at (1.4381, 6.1557){$\mathbf{c}$};

  \path[draw=lime,fill=lime,line width=0.0cm] (3.6407, 3.9846) ellipse (0.291cm and 0.291cm);

  \node[text=white,anchor=south,line width=0.0cm,rotate=45.0] (text26) at (3.7629, 3.8359){$\mathbf{a}$};

  \path[draw=c00bb00,fill=c00bb00,line width=0.0cm] (4.281, 3.9846) circle (0.291cm);

  \node[text=white,anchor=south,line width=0.0cm,rotate=45.0] (text27) at (4.4502, 3.8469){$\mathbf{b}$};

  \path[draw=c007700,fill=c007700,line width=0.0cm] (4.9213, 3.9846) circle (0.291cm);

  \node[text=white,anchor=south,line width=0.0cm,rotate=45.0] (text28) at (5.0471, 3.8274){$\mathbf{c}$};

  \path[draw=lime,fill=lime,line width=0.0cm] (3.6407, 1.6563) ellipse (0.291cm and 0.291cm);

  \node[text=white,anchor=south,line width=0.0cm,rotate=45.0] (text30) at (3.7629, 1.5076){$\mathbf{a}$};

  \path[draw=c00bb00,fill=c00bb00,line width=0.0cm] (4.281, 1.6563) circle (0.291cm);

  \node[text=white,anchor=south,line width=0.0cm,rotate=45.0] (text31) at (4.4502, 1.5185){$\mathbf{b}$};

  \path[draw=c007700,fill=c007700,line width=0.0cm] (4.9213, 1.6563) circle (0.291cm);

  \node[text=white,anchor=south,line width=0.0cm,rotate=45.0] (text32) at (5.0471, 1.4991){$\mathbf{c}$};

  \path[draw=lime,fill=lime,line width=0.0cm] (5.969, 10.9696) circle (0.291cm);

  \node[text=white,anchor=south,line width=0.0cm,rotate=45.0] (text35) at (6.0913, 10.8209){$\mathbf{a}$};

  \path[draw=c007700,fill=c007700,line width=0.0cm] (6.6093, 10.9696) circle (0.291cm);

  \node[text=white,anchor=south,line width=0.0cm,rotate=45.0] (text36) at (6.7351, 10.8124){$\mathbf{c}$};

  \path[draw=c00bb00,fill=c00bb00,line width=0.0cm] (5.969, 3.9846) circle (0.291cm);

  \node[text=white,anchor=south,line width=0.0cm,rotate=45.0] (text39) at (6.1383, 3.8469){$\mathbf{b}$};

  \path[draw=lime,fill=lime,line width=0.0cm] (8.2973, 8.6413) circle (0.291cm);

  \node[text=white,anchor=south,line width=0.0cm,rotate=45.0] (text41) at (8.4196, 8.4926){$\mathbf{a}$};

  \path[draw=c00bb00,fill=c00bb00,line width=0.0cm] (8.9376, 8.6413) circle (0.291cm);

  \node[text=white,anchor=south,line width=0.0cm,rotate=45.0] (text42) at (9.1069, 8.5035){$\mathbf{b}$};

  \path[draw=c007700,fill=c007700,line width=0.0cm] (9.5779, 8.6413) circle (0.291cm);

  \node[text=white,anchor=south,line width=0.0cm,rotate=45.0] (text43) at (9.7037, 8.4841){$\mathbf{c}$};

  \path[draw=lime,fill=lime,line width=0.0cm] (8.2973, 6.313) ellipse (0.291cm and 0.291cm);

  \node[text=white,anchor=south,line width=0.0cm,rotate=45.0] (text45) at (8.4196, 6.1643){$\mathbf{a}$};

  \path[draw=c00bb00,fill=c00bb00,line width=0.0cm] (8.9376, 6.313) ellipse (0.291cm and 0.291cm);

  \node[text=white,anchor=south,line width=0.0cm,rotate=45.0] (text46) at (9.1069, 6.1752){$\mathbf{b}$};

\end{tikzpicture}}}}

\vspace{-3cm}
    \hspace{2cm}
	{\scalefont{1}{\adjustbox{}{\rotatebox{0}{\scalebox{0.4}{\input{images/tikz/front_page_lattice.tex}}}}}}
	\hspace{2cm}
	{\scalefont{1.5}{\adjustbox{raise=1.75cm}{\rotatebox{-45}{\scalebox{0.4}{\definecolor{lime}{RGB}{0,255,0}
\definecolor{cff0008}{RGB}{255,0,8}

\def \globalscale {1.000000}
\begin{tikzpicture}[y=1cm, x=1cm, yscale=\globalscale,xscale=\globalscale, every node/.append style={scale=\globalscale}, inner sep=0pt, outer sep=0pt]
  \path[fill=black,fill opacity=0.2,line cap=butt,line join=miter,line width=0.0265cm] (3.0113, 2.7573) -- (3.0113, 5.1559) -- (7.8086, 5.1559) -- (7.8086, 2.7573) -- (3.0113, 2.7573);

  \path[fill=blue,fill opacity=0.2,line cap=butt,line join=miter,line width=0.0265cm] (5.4099, 2.7573) -- (5.4099, 7.5546) -- (0.6126, 7.5546) -- (0.6126, 2.7573) -- (5.4099, 2.7573);

  \path[fill=red,fill opacity=0.2,line cap=butt,line join=miter,line width=0.0265cm] (10.2073, 7.5546) -- (0.6126, 7.5546) -- (0.6126, 2.7573) -- (3.0113, 2.7573) -- (3.0113, 5.1559) -- (10.2073, 5.1559) -- cycle;

  \path[draw=white,fill=lime,fill opacity=0.2,line cap=butt,line join=miter,line width=0.0265cm] (10.2073, 7.5546) -- (12.606, 7.5546) -- (12.606, 12.352) -- (5.4099, 12.352) -- (5.4099, 2.7573) -- (7.8086, 2.7573) -- (7.8086, 9.9533) -- (10.2073, 9.9533) -- cycle;

  \begin{scope}[draw=black,line cap=round,line width=0.1588cm]
    \path[draw=black,line cap=round,line width=0.1588cm] (3.0113, 5.1559) -- (3.0113, 0.3586);

    \path[draw=black,line cap=round,line width=0.1588cm] (3.0113, 0.3586) -- (5.4099, 0.3586);

    \path[draw=black,line cap=round,line width=0.1588cm] (5.4099, 0.3586) -- (5.4099, 12.352);

    \path[draw=black,line cap=round,line width=0.1588cm] (5.4099, 12.352) -- (12.606, 12.352);

    \path[draw=black,line cap=round,line width=0.1588cm] (12.606, 12.352) -- (12.606, 7.5546);

    \path[draw=black,line cap=round,line width=0.1588cm] (12.606, 7.5546) -- (0.6126, 7.5546);

    \path[draw=black,line cap=round,line width=0.1588cm] (0.6126, 7.5546) -- (0.6126, 2.7573);

    \path[draw=black,line cap=round,line width=0.1588cm] (0.6126, 2.7573) -- (7.8086, 2.7573);

    \path[draw=black,line cap=round,line width=0.1588cm] (7.8086, 2.7573) -- (7.8086, 9.9533);

    \path[draw=black,line cap=round,line width=0.1588cm] (7.8086, 9.9533) -- (10.2073, 9.9533);

    \path[draw=black,line cap=round,line width=0.1588cm] (10.2073, 9.9533) -- (10.2073, 5.1559);

    \path[draw=black,line cap=round,line width=0.1588cm] (10.2073, 5.1559) -- (3.0113, 5.1559);

  \end{scope}
  \tikzstyle{every node}=[font=\fontsize{20}{20}\selectfont]
  \node[text=black,anchor=south,line cap=butt,line join=miter,line width=0.0cm,miter limit=4.0,rotate=45.0] (text2) at (8.4449, 9.2281){$4$};

  \node[text=black,anchor=south,line cap=butt,line join=miter,line width=0.0cm,miter limit=4.0,rotate=45.0] (text3) at (8.5325, 6.8271){$5$};

  \node[text=black,anchor=south,line cap=butt,line join=miter,line width=0.0cm,miter limit=4.0,rotate=45.0] (text4) at (6.1416, 6.8339){$6$};

  \node[text=black,anchor=south,line cap=butt,line join=miter,line width=0.0cm,miter limit=4.0,rotate=45.0] (text5) at (6.1451, 4.4263){$2$};

  \node[text=black,anchor=south,line cap=butt,line join=miter,line width=0.0cm,miter limit=4.0,rotate=45.0] (text6) at (3.7409, 4.4279){$3$};

  \node[text=black,anchor=south,line cap=butt,line join=miter,line width=0.0cm,miter limit=4.0,rotate=45.0] (text7) at (1.3429, 6.8259){$8$};

  \node[text=black,anchor=south,line cap=butt,line join=miter,line width=0.0cm,miter limit=4.0,rotate=45.0] (text8) at (3.7092, 2.0579){$1$};

  \node[text=black,anchor=south,line cap=butt,line join=miter,line width=0.0cm,miter limit=4.0,rotate=45.0] (text9) at (6.1398, 2.0286){$9$};

  \node[text=black,anchor=south,line cap=butt,line join=miter,line width=0.0cm,miter limit=4.0,rotate=45.0] (text13) at (6.1312, 11.6298){$7$};

  \path[fill=lime,even odd rule,line cap=round,line width=0.0865cm] (5.4099, 2.7573) circle (0.3198cm);

  \path[fill=black,even odd rule,line cap=round,line width=0.0865cm] (7.8086, 5.1559) circle (0.2041cm);

  \path[fill=yellow,even odd rule,line cap=round,line width=0.0865cm] (10.2073, 7.5546) circle (0.3198cm);

  \path[fill=magenta,even odd rule,line cap=round,line width=0.0865cm] (10.2073, 7.5546) circle (0.1701cm);

  \path[fill=lime,even odd rule,line cap=round,line width=0.0865cm] (10.2073, 7.5546) circle (0.3198cm);

  \path[fill=blue,even odd rule,line cap=butt,line join=miter,line width=0.0529cm,miter limit=4.0] (5.4099, 2.7573) circle (0.2041cm);

  \path[fill=cff0008,even odd rule,line cap=round,line width=0.0865cm] (3.0113, 2.7573) circle (0.3198cm);

  \path[fill=black,even odd rule,line cap=round,line width=0.0865cm] (3.0113, 2.7573) circle (0.2041cm);

  \path[fill=cff0008,even odd rule,line cap=round,line width=0.0865cm] (10.2073, 7.5546) circle (0.2041cm);

  \path[fill=blue,even odd rule,line cap=butt,line join=miter,line width=0.0529cm,miter limit=4.0] (5.4099, 7.5546) circle (0.2041cm);

  \node[text=black,anchor=south,line width=0.0cm,rotate=45.0] (text23) at (9.9877, 3.3796){$1\wedge4\wedge(2\vee3)\wedge(3\vee8)\wedge(2\vee6\vee7)\wedge(5\vee6\vee8)$};

  \path[draw=blue,fill=blue,fill opacity=0.2,line cap=butt,line join=miter,line width=0.0423cm] (8.9706, 1.6256).. controls (8.9276, 2.4425) and (9.8831, 2.5388) .. (9.8831, 2.5388).. controls (9.8831, 2.5388) and (9.7875, 1.5826) .. (8.9706, 1.6256) -- cycle;

  \path[draw=blue,fill=blue,even odd rule,line cap=butt,line join=miter,line width=0.0529cm,miter limit=4.0,rotate around={45.0:(0.0, 12.9646)}] (-1.6747, -1.3965) circle (0.1701cm);

  \path[draw=blue,fill=blue,even odd rule,line cap=butt,line join=miter,line width=0.0529cm,miter limit=4.0,rotate around={45.0:(0.0, 12.9646)}] (-0.3628, -1.3963) circle (0.1701cm);

  \path[draw=blue,line cap=butt,line join=miter,line width=0.0529cm] (8.9706, 1.6256).. controls (8.9276, 2.4425) and (9.8831, 2.5388) .. (9.8831, 2.5388).. controls (9.8831, 2.5388) and (9.7875, 1.5826) .. (8.9706, 1.6256) -- cycle;

  \path[draw=red,line cap=butt,line join=miter,line width=0.0529cm] (13.9156, 6.6235).. controls (13.8726, 7.4404) and (14.8281, 7.5367) .. (14.8281, 7.5367).. controls (14.8281, 7.5367) and (14.7325, 6.5805) .. (13.9156, 6.6235) -- cycle;

  \path[draw=red,fill=red,even odd rule,line cap=round,line width=0.0865cm,rotate around={45.0:(0.0, 12.9646)}] (5.356, -1.3591) circle (0.1701cm);

  \path[draw=red,fill=red,even odd rule,line cap=round,line width=0.0865cm,rotate around={45.0:(0.0, 12.9646)}] (6.6678, -1.3589) circle (0.1701cm);

  \path[fill=red,fill opacity=0.2,line cap=butt,line join=miter,line width=0.0423cm] (13.9156, 6.6235).. controls (13.8726, 7.4404) and (14.8281, 7.5367) .. (14.8281, 7.5367).. controls (14.8281, 7.5367) and (14.7325, 6.5805) .. (13.9156, 6.6235) -- cycle;

  \path[draw=lime,fill=lime,fill opacity=0.2,line cap=butt,line join=miter,line width=0.0265cm] (11.2579, 3.9129).. controls (11.2149, 4.7298) and (12.1704, 4.8261) .. (12.1704, 4.8261).. controls (12.1704, 4.8261) and (12.0748, 3.8699) .. (11.2579, 3.9129) -- cycle;

  \path[draw=lime,fill=lime,even odd rule,line cap=round,line width=0.0865cm,rotate around={45.0:(0.0, 12.9646)}] (1.56, -1.3965) circle (0.1701cm);

  \path[draw=lime,fill=lime,even odd rule,line cap=round,line width=0.0865cm,rotate around={45.0:(0.0, 12.9646)}] (2.8719, -1.3963) circle (0.1701cm);

  \path[draw=lime,line cap=butt,line join=miter,line width=0.0423cm,miter limit=4.0] (11.2579, 3.9129).. controls (11.2149, 4.7298) and (12.1704, 4.8261) .. (12.1704, 4.8261).. controls (12.1704, 4.8261) and (12.0748, 3.8699) .. (11.2579, 3.9129) -- cycle;

  \begin{scope}[cm={ 0.7071,0.7071,-0.7071,0.7071,(2.6993, -5.8923)}]
    \path[fill=black,fill opacity=0.2,line cap=butt,line join=miter,line width=0.0265cm] (7.045, 0.8816).. controls (7.5922, 1.4896) and (8.336, 0.8821) .. (8.336, 0.8821).. controls (8.336, 0.8821) and (7.5922, 0.2735) .. (7.045, 0.8816) -- cycle;

    \path[fill=black,even odd rule,line cap=round,line width=0.0865cm] (7.045, 0.8816) circle (0.1701cm);

    \path[fill=black,even odd rule,line cap=round,line width=0.0865cm] (8.3569, 0.8817) circle (0.1701cm);

    \path[draw=black,fill opacity=0.2,line cap=butt,line join=miter,line width=0.0423cm,miter limit=4.0] (7.045, 0.8816).. controls (7.5922, 1.4896) and (8.336, 0.8821) .. (8.336, 0.8821).. controls (8.336, 0.8821) and (7.5922, 0.2735) .. (7.045, 0.8816) -- cycle;

  \end{scope}

\end{tikzpicture}}}}}}
    
\vspace{-2cm}
\caption{Top: The \href{https://christopherlloyd.github.io/LooPindex/multiloops/9\%5E1_5.html}{loop $9^1_{5}$} has $2$ optimal pinning sets, and $3$ other minimal pinning sets. \\
Bottom: Part of the pinning ideal of the \href{https://christopherlloyd.github.io/LooPindex/multiloops/9\%5E1_5.html}{loop $9^1_{5}$}, and its mobidisc formula.}
 \label{fig:intro-figure}
\end{figure}

\newpage

\renewcommand{\contentsname}{Plan of the paper}
\setcounter{tocdepth}{2}
\tableofcontents

\section{Introduction}

\subsection*{Multiloops and their pinning ideals}

Let $\Surf$ be a closed oriented smooth surface (we will often focus on the case of a sphere, which already contains most of the intricacies).
For $P\subset \Surf$ we write $\Surf_P=\Surf\setminus P$. 
Denote by $\sqcup_1^s \Sph^1$ a disjoint union of $s\in \N$ oriented circles.

A \emph{multiloop} with $s$ strands is a generic immersion $\gamma \colon \sqcup_1^s \Sph^1 \looparrowright \Surf_P$ considered up to orientation preserving diffeomorphisms of the source and the target.
In particular, it has a finite number of multiple points, and those are all transverse double points: we denote this number by $\# \gamma \in \N$.

The \emph{regions} of $\gamma$ are the connected components of $\Surf \setminus \im(\gamma)$, indexed by $\Regions =\pi_0(\Surf \setminus \im(\gamma))$.
The multiloop $\gamma$ is \emph{filling} $\Surf$ when $\Surf\setminus \im(\gamma)$ is homotopic to a disjoint union of discs. 

A \emph{multicurve} is a homotopy class of multiloops.
The set of multicurves in $\Surf_P$ is endowed with the \emph{self-intersection} function $\si_P$ counting the minimal number of double points among its generic representatives.
A multiloop $\gamma\colon \sqcup \Sph^1 \looparrowright \Surf_P$ is \emph{taut} when it has the minimal number of double points in its homotopy class, in formula $\si_P(\gamma)=\# \gamma$.

\begin{definition}[pinning set]
Consider a multiloop $\gamma \colon \sqcup_1^s \Sph^1 \looparrowright \Surf$ with regions $\Regions$.

For a subset $P\subset \Regions$, we denote by $\si_P(\gamma)$ the self-intersection number of the multicurve representing $\gamma$ in $\Surf_P$, and we say that $P$ is \emph{pinning} $\gamma$ when $\si_P(\gamma)=\# \gamma$, namely $\gamma$ becomes taut after puncturing the surface $\Surf$ at $P$.

The pinning sets form the \emph{pinning ideal} $\mathcal{PI}(\gamma)\subset \Parts(\Regions)$: a sub-poset which is absorbing under union and which contains $\Regions$.
A pinning set is called \emph{minimal} if it is minimal with respect to inclusion, and \emph{optimal} if it has the minimum cardinal among all pinning sets.
The \emph{pinning number} $\varpi(\gamma)\in \N$ is the cardinal of its optimal pinning sets.
\end{definition}

\begin{question}[goal]
Given a multiloop $\gamma \colon \sqcup_1^s \Sph^1 \looparrowright \Surf$, how (efficiently) can we:
\begin{itemize}[noitemsep]
    \item Construct minimal or optimal pinning sets? 
    \item Find which regions belong to every pinning set?
    \item Compute the pinning number $\varpi(\gamma)$?
\end{itemize}
More generally, we wonder what can be the shape of the pinning ideal of a (multi)loop, and what are its statistics for certain random models.
\end{question}

We will address all of these questions, with a particular emphasis on the following.
\begin{theorem}[\textsc{MultiLooPinNum}]
\label{thm:Intro:MultiLooPinNum-is-SNPC}
    The \textsc{MultiLooPinNum} problem defined by 
    \begin{itemize}[align=left, noitemsep]
        \item[Instance:] A filling multiloop $\gamma \colon \sqcup_1^s \Sph^1 \looparrowright \Surf$, and an integer $p\in \N$.
        \item[Question:] Does $\gamma$ have pinning number $\varpi(\gamma)\le p$?
    \end{itemize}
    is \textsf{NP}-complete, and so is its restriction to loops in the sphere $\gamma\colon \Sph^1 \looparrowright \Sph^2$.
\end{theorem}

\subsection{Combinatorial group theory: counting self-intersections}

A nontrivial filling multiloop $\gamma\colon \sqcup_1^s \Sph^1 \looparrowright \Surf$ corresponds, up to orientation of its strands, to a $4$-regular map $\Gamma \subset \Surf$, which can be encoded by a pair of permutations over its set of half-edges.
The size of such an encoding yields a measure for its complexity: it is proportional to $\Card(\Regions)$.

\begin{algorithm}[self-intersection]
\label{algo:Intro:self-intersection}
    We will describe an algorithm which has:
    \begin{itemize}[noitemsep, align=left]
        \item[Input:] A multiloop $\gamma\colon \sqcup_1^s \Sph^1 \looparrowright \Surf$ with a subset of regions $P\subset \Regions$.
        \item[Output:] The self-intersection number $\si_P(\gamma)$ of the multicurve $\gamma \colon \sqcup_1^s \Sph^1 \looparrowright \Surf_P$.
        \item[Complexity:] Running time $O(\Card(\Regions)^2)$.
    \end{itemize}
\end{algorithm}

We will present Algorithm \ref{algo:Intro:self-intersection} for multiloops in the sphere with $P\ne \emptyset$, in such a way that it holds without essential modification to multiloops in any surface with $P\ne \emptyset$.
It first computes a convenient presentation of the free group $\pi_1(\Surf_\Regions)$. 
Its generators correspond to the edges of a spanning tree of the dual map of $\Gamma\subset \Surf$.
The homotopy classes of the loops $\gamma_i\subset \Surf_\Regions$ correspond to conjugacy classes in that free group, whence to cyclically reduced words in the generators.
Filling the regions $Q = \Regions\setminus P$ yields explicit relations providing a presentation of $\pi_1(\Surf_P)$, and we deduce the conjugacy classes associated to $\gamma_i \subset \Surf_P$ by applying simple rewriting rules.
Finally we adapt algorithms of Birman--Series \cite{Birman-Series_algo-simple-closed-curves_1984} and Cohen--Lustig \cite{Cohen-Lustig_algo-geometric-intersection-number_1987} to compute the intersection number between homotopy classes.

\begin{remark}[closed surfaces]
The Algorithm \ref{algo:Intro:self-intersection} can be adapted to handle the case $P=\emptyset$.

The case of the closed sphere is trivial.
The case of the closed torus is easy: every loop is homotopic to the power of a simple loop, whose self-intersection and mutual intersection numbers can be computed from primitivity exponents and algebraic intersection numbers.
The case of closed surfaces of genus $g\ge 2$ is more subtle, but our description of Algorithm \ref{algo:Intro:self-intersection} can be extended by adapting \cite{Birman-Series_Dehn-algo-simple-closed-curves_1987, Lustig_Paths-geodesics-geometric-intersection-II_1987}.

Still, all algorithms should take care of these special cases individually. 
Our cautionary Examples \ref{ex:increasing-si-any-puncture} explain why one cannot reduce to the case $P\ne \emptyset$ so naively.
\end{remark}

\begin{corollary}[\textsc{Pinning} and \textsc{PinMin} are \textsf{P}]
\label{cor:intro:LooPinMin-is-P}
Given a plane multiloop $\gamma \colon \sqcup_1^s \Sph^1 \looparrowright \Sph^2$ and a subset of regions $P\subset \Regions$, we can successively:
    \begin{itemize}[noitemsep, align=left]
        \item[\textsc{Pinning}:] certify that $P$ is pinning in time $O(\Card(\Regions)^2)$, and if so
        \item[\textsc{PinMin}:] check if $P$ is a minimal pinning set in time $O(\Card(P)\Card(\Regions)^2)$, if not 
        \item[\textsc{PinMinC}:] construct a minimal pinning set $M\subset P$ in time $O(\Card(P)\Card(\Regions)^2)$.
    \end{itemize}
\end{corollary}

\begin{corollary}
\label{cor:Intro:MultiLooPinNum-is-NP}
    The \textsc{MultiLooPinNum} problem is in \textsf{NP}.
\end{corollary}

\subsection{Geometric topology: computing immersed discs}

We now focus on loops, namely multiloops with one strand.
Let us recall a characterization of taut loops $\gamma\colon \Sph^1\looparrowright \Surf_P$ obtained in the proof of \cite[Theorem 4.2]{Hass-Scott_Intersection-curves-surfaces_1985}.
First define a \emph{singular monorbigon} of $\gamma$ as a null-homotopic subloop delimited by one or two double points.
Next, consider the singular monorbigons that bound an immersed disc in the surface $\Surf_P$: such a disc consists of a collection of regions $D\subset \Regions$ that we call a \emph{mobidisc}.
Denote by $\MoB(\gamma)\subset \Parts(\Regions)$ the collection of mobidiscs, so $\MoB(\gamma)\setminus\{R\}$ refers to the collection of \emph{proper} mobidiscs.

A close reading of the proofs in \cite{Hass-Scott_Intersection-curves-surfaces_1985} yield an improvement of their theorem which we reformulate in terms of pinning sets as follows.

\begin{theorem}[pinning mobidiscs]
    \label{thm:Intro:pinning-mobidiscs}
    For a loop $\gamma \colon \Sph^1 \looparrowright \Surf$ with regions $\Regions$, a set $P\in \Parts(\Regions)$ is pinning if and only if it intersects every proper mobidisc $D\in \MoB(\gamma)\setminus\{\Regions\}$.
\end{theorem}

\begin{remark}[monotone CNF]
    Theorem \ref{thm:Intro:pinning-mobidiscs} implies that the pinning sets of the loop $\gamma$ correspond to the solutions of the monotone conjunctive normal form associated to $\MoB(\gamma)\subset \Parts(\Regions)$, or equivalently to the vertex covers of the hypergraph whose vertices are indexed by the regions $\Regions$ and hyperedges correspond to the mobidiscs $\MoB(\gamma)\subset \Parts(\Regions)$.
    Thus, it is possible that a particular type of \textsc{SAT}-solver may be used to find the minimum cardinal for these solutions efficiently.
\end{remark}

\begin{remark}[just for loops]
    Theorem \ref{thm:Intro:pinning-mobidiscs} only holds for loops.
    Indeed, the Figure \ref{fig:no_singular_monorbigon_non_taut_multiloop}, which is a slight modification of \cite[Figure 0.1]{Hass-Scott_Intersection-curves-surfaces_1985} exhibits a multiloop with $2$ strands in a thrice-punctured sphere, which is not taut but has no singular monorbigons at all. 
\end{remark}

For a loop in the plane, its subset of immersed monorbigons can be computed in polynomial time by an algorithm of Blank \cite{Blank_extending-immersions-circle_1967, Poenaru_extension-immersions-codim-1_1995} and Shor--Van Wyk \cite{Shor-VanWyk_Detecting-decomposing-self-overlapping-curves_1992}, and their associated mobidiscs may also be computed in polynomial time.
Hence for a spherical loop $\gamma\colon \Sph^1\looparrowright \Sph^2$, we may run over all possible choices of the region at infinity to compute its collection of proper mobidiscs.

\begin{algorithm}[pinning mobidiscs]
    \label{algo:Intro:mobidisc}
    For a loop $\gamma \colon \Sph^1 \looparrowright \Sph^2$ with regions $\Regions$, we may compute its collection of proper mobidiscs $\MoB(\gamma)\setminus\{\Regions\}$ in time $O(\Card(\Regions)^5\log(\Card(\Regions)))$.

\end{algorithm}

\begin{remark}[computing mobidiscs in higher genus]
In higher genus surfaces, the characterization of loops which bound immersed discs and their computation has been attempted in \cite{Frisch_classification-immersions-curves-in-surfaces_2010}, and several works address related problems 
\cite{Ezell-Marx_Branched-extensions-curves-orientable-surfaces_1980, Pappas_Extensions-codim-1-immersions_1996}, but it is not clear whether any of these provides a clear computable criterion, let alone of polynomial time complexity.
\end{remark}

\subsection{From vertex covers of planar graphs to pinning sets of loops}

We know from Corollary \ref{cor:Intro:MultiLooPinNum-is-NP} that the \textsc{MultiLooPinNum} problem is in \textsf{NP}. 

\begin{algorithm}[reduction from \textsc{planar vertex cover} to \textsc{LooPinNum}]
\label{algo:intro:planar-vertex-cover-to-LooPinNum}
    We describe a polynomial time algorithm which to a planar graph $G=(V,E)$ associates a plane loop $\gamma \colon \Sph^1 \looparrowright \R^2$ with $O(\Card(E)^2)$ double points, such that the $k$-vertex-covers of $G$ are in one to one correspondence with the $(k+6\Card(E)$)-pinning-sets of $\gamma$.
\end{algorithm}

Note that the \textsc{LooPinNum} problem contains as a sub-problem the problem of pinning loops in the sphere, so its complexity is at least as hard as the \textsc{planar vertex cover} problem, which is known to be \textsf{NP}-complete.

\begin{corollary}
    The \textsc{LooPinNum} problem is \textsf{NP}-hard, even in the sphere.
\end{corollary}

To justify the reduction, we will derive from \cite[Theorem 4.2]{Hass-Scott_Intersection-curves-surfaces_1985} a Lemma \ref{lem:link(singular-monorbigons)=0} ensuring a sufficient condition that a set of regions is pinning in terms of linking numbers of singular monorbigons with pins.

\subsection{The pinning ideal of a multiloop}

Consider a multiloop $\gamma\colon \sqcup_1^s \Sph^1 \looparrowright \Surf$ with regions $\Regions$.
The pinning ideal $\mathcal{PI}\subset \Parts(\Regions)$ carries all the information related to the pinning problem for the multiloop.
As an ideal, it is generated by the minimal pinning sets.
It contains the \emph{pinning semi-lattice}, which is the sub-poset obtained by all possible unions of minimal pinning sets. 

Using Algorithm \ref{algo:Intro:self-intersection} and some functionalities from \texttt{plantri} \cite{plantri}, we computed the pinning ideals and semi-lattices of all irreducible indecomposible spherical multiloops with at most $12$ regions (unoriented and up to reflection), and the statistics of certain numerical parameters.
The results are available in the \href{https://github.com/ChristopherLloyd/LooPindex}{LooPindex} \cite{Simon-Stucky_LooPindex_2024}.

We make a few observations in section \ref{sec:pinning-ideal}, either suggesting certain heuristics for approximation algorithms, or invalidating some naive conjectures about the behavior of pinning ideals.
Let us record the main questions and observations.

\begin{question}[structure of pinning ideals] 
Which ideals arise as pinning ideals of filling loops or multiloops? 

How do pinning ideals of filling multiloops behave under Reidemeister moves, flypes, and crossing resolutions? 

How do pinning ideals of loops behave under the operations of spheric-sums and toric-sums introduced in \cite{CLS_ChorDiaGraFiloop_2023}?
\end{question}

The Figures \ref{fig:pinnum-R3} and \ref{fig:pin-flype} show that the pinning number can change under $R3$ moves and flypes as well as mutations, even for indecomposible loops in the sphere.

\begin{figure}[h]
    \centering
    {\scalefont{1}\scalebox{0.7}{\definecolor{cd82626}{RGB}{216,38,38}
\definecolor{c2626d8}{RGB}{38,38,216}
\definecolor{lime}{RGB}{0,255,0}
\definecolor{c006600}{RGB}{0,102,0}
\definecolor{c00cc00}{RGB}{0,204,0}
\definecolor{c009900}{RGB}{0,153,0}

\def \globalscale {1.000000}
\begin{tikzpicture}[font=\bf,y=1cm, x=1cm, yscale=\globalscale,xscale=\globalscale, every node/.append style={scale=\globalscale}, inner sep=0pt, outer sep=0pt]
  \path[draw=cd82626,fill=cd82626] (6.6146, 8.296) circle (0.0529cm);

  \path[draw=cd82626,fill=cd82626] (10.4951, 8.296) circle (0.0529cm);

  \path[draw=cd82626,fill=cd82626] (10.4951, 0.5437) circle (0.0529cm);

  \path[draw=cd82626,fill=cd82626] (4.6743, 0.5437) circle (0.0529cm);

  \path[draw=cd82626,fill=cd82626] (4.6743, 8.296) circle (0.0529cm);

  \path[draw=cd82626,fill=cd82626] (0.7938, 8.296) circle (0.0529cm);

  \path[draw=cd82626,fill=cd82626] (0.7938, 2.4818) circle (0.0529cm);

  \path[draw=cd82626,fill=cd82626] (6.6146, 2.4818) circle (0.0529cm);

  \path[draw=c2626d8,fill=c2626d8] (2.734, 4.4199) circle (0.0529cm);

  \path[draw=c2626d8,fill=c2626d8] (2.734, 6.3579) circle (0.0529cm);

  \path[draw=c2626d8,fill=c2626d8] (8.5549, 6.3579) circle (0.0529cm);

  \path[draw=c2626d8,fill=c2626d8] (8.5549, 2.4818) circle (0.0529cm);

  \path[draw=c2626d8,fill=c2626d8] (12.4354, 2.4818) circle (0.0529cm);

  \path[draw=c2626d8,fill=c2626d8] (12.4354, 4.4199) circle (0.0529cm);

  \path[draw=cd82626,line width=0.1588cm] (6.6146, 8.296) -- (10.4951, 8.296);

  \path[draw=cd82626,line width=0.1588cm] (10.4951, 8.296) -- (10.4951, 0.5437);

  \path[draw=cd82626,line width=0.1588cm] (10.4951, 0.5437) -- (4.6743, 0.5437);

  \path[draw=cd82626,line width=0.1588cm] (4.6743, 0.5437) -- (4.6743, 8.296);

  \path[draw=cd82626,line width=0.1588cm] (4.6743, 8.296) -- (0.7938, 8.296);

  \path[draw=cd82626,line width=0.1588cm] (0.7938, 8.296) -- (0.7938, 2.4818);

  \path[draw=cd82626,line width=0.1588cm] (0.7938, 2.4818) -- (6.6146, 2.4818);

  \path[draw=cd82626,line width=0.1588cm] (6.6146, 2.4818) -- (6.6146, 8.296);

  \path[draw=c2626d8,line width=0.1588cm] (2.734, 4.4199) -- (2.734, 6.3579);

  \path[draw=c2626d8,line width=0.1588cm] (2.734, 6.3579) -- (8.5549, 6.3579);

  \path[draw=c2626d8,line width=0.1588cm] (8.5549, 6.3579) -- (8.5549, 2.4818);

  \path[draw=c2626d8,line width=0.1588cm] (8.5549, 2.4818) -- (12.4354, 2.4818);

  \path[draw=c2626d8,line width=0.1588cm] (12.4354, 2.4818) -- (12.4354, 4.4199);

  \path[draw=c2626d8,line width=0.1588cm] (12.4354, 4.4199) -- (2.734, 4.4199);

  \path[draw=red,fill=red,line width=0.0cm] (6.8822, 8.0284) ellipse (0.1384cm and 0.1384cm);

  \node[text=white,anchor=south,line width=0.0cm] (text1081) at (6.8822, 7.9375){A};

  \path[draw=lime,fill=lime,line width=0.0cm] (7.1868, 8.0284) ellipse (0.1384cm and 0.1384cm);

  \node[text=white,anchor=south,line width=0.0cm] (text8556) at (7.1868, 7.9375){a};

  \path[draw=c006600,fill=c006600,line width=0.0cm] (7.4913, 8.0284) ellipse (0.1384cm and 0.1384cm);

  \node[text=white,anchor=south,line width=0.0cm] (text1879) at (7.4913, 7.9375){d};

  \path[draw=lime,fill=lime,line width=0.0cm] (4.9419, 2.2142) circle (0.1384cm);

  \node[text=white,anchor=south,line width=0.0cm] (text8675) at (4.9419, 2.1167){a};

  \path[draw=red,fill=red,line width=0.0cm] (10.7628, 4.1522) circle (0.1384cm);

  \node[text=white,anchor=south,line width=0.0cm] (text1008) at (10.7628, 4.0481){A};

  \path[draw=lime,fill=lime,line width=0.0cm] (11.0673, 4.1522) circle (0.1384cm);

  \node[text=white,anchor=south,line width=0.0cm] (text9574) at (11.0673, 4.0481){a};

  \path[draw=c00cc00,fill=c00cc00,line width=0.0cm] (11.3719, 4.1522) circle (0.1384cm);

  \node[text=white,anchor=south,line width=0.0cm] (text8975) at (11.3719, 4.0481){b};

  \path[draw=c009900,fill=c009900,line width=0.0cm] (11.6764, 4.1522) circle (0.1384cm);

  \node[text=white,anchor=south,line width=0.0cm] (text2121) at (11.6764, 4.0481){c};

  \path[draw=c006600,fill=c006600,line width=0.0cm] (11.981, 4.1522) circle (0.1384cm);

  \node[text=white,anchor=south,line width=0.0cm] (text3289) at (11.981, 4.0481){d};

  \path[draw=c00cc00,fill=c00cc00,line width=0.0cm] (8.8225, 4.1522) ellipse (0.1384cm and 0.1384cm);

  \node[text=white,anchor=south,line width=0.0cm] (text742) at (8.8225, 4.0481){b};

  \path[draw=c009900,fill=c009900,line width=0.0cm] (9.1271, 4.1522) ellipse (0.1384cm and 0.1384cm);

  \node[text=white,anchor=south,line width=0.0cm] (text8618) at (9.1271, 4.0481){c};

  \path[draw=c006600,fill=c006600,line width=0.0cm] (9.4316, 4.1522) ellipse (0.1384cm and 0.1384cm);

  \node[text=white,anchor=south,line width=0.0cm] (text6318) at (9.4316, 4.0481){d};

  \path[draw=lime,fill=lime,line width=0.0cm] (1.0614, 8.0284) circle (0.1384cm);

  \node[text=white,anchor=south,line width=0.0cm] (text7390) at (1.0614, 7.9375){a};

  \path[draw=c00cc00,fill=c00cc00,line width=0.0cm] (1.3659, 8.0284) circle (0.1384cm);

  \node[text=white,anchor=south,line width=0.0cm] (text5880) at (1.3659, 7.9375){b};

  \path[draw=c006600,fill=c006600,line width=0.0cm] (1.6705, 8.0284) circle (0.1384cm);

  \node[text=white,anchor=south,line width=0.0cm] (text52) at (1.6705, 7.9375){d};

  \path[draw=red,fill=red,line width=0.0cm] (3.0017, 6.0903) ellipse (0.1384cm and 0.1384cm);

  \node[text=white,anchor=south,line width=0.0cm] (text8270) at (3.0017, 5.9796){A};

  \path[draw=lime,fill=lime,line width=0.0cm] (3.3062, 6.0903) circle (0.1384cm);

  \node[text=white,anchor=south,line width=0.0cm] (text3723) at (3.3062, 5.9796){a};

  \path[draw=c00cc00,fill=c00cc00,line width=0.0cm] (3.6108, 6.0903) circle (0.1384cm);

  \node[text=white,anchor=south,line width=0.0cm] (text864) at (3.6108, 5.9796){b};

  \path[draw=c009900,fill=c009900,line width=0.0cm] (3.9153, 6.0903) circle (0.1384cm);

  \node[text=white,anchor=south,line width=0.0cm] (text5626) at (3.9153, 5.9796){c};

  \path[draw=c006600,fill=c006600,line width=0.0cm] (4.2199, 6.0903) circle (0.1384cm);

  \node[text=white,anchor=south,line width=0.0cm] (text617) at (4.2199, 5.9796){d};

  \path[draw=c006600,fill=c006600,line width=0.0cm] (4.9419, 6.0903) circle (0.1384cm);

  \node[text=white,anchor=south,line width=0.0cm] (text2570) at (4.9419, 5.9796){d};

  \path[draw=red,fill=red,line width=0.0cm] (4.9419, 4.1522) ellipse (0.1384cm and 0.1384cm);

  \node[text=white,anchor=south,line width=0.0cm] (text4189) at (4.9419, 4.0481){A};

  \path[draw=c009900,fill=c009900,line width=0.0cm] (5.2465, 4.1522) ellipse (0.1384cm and 0.1384cm);

  \node[text=white,anchor=south,line width=0.0cm] (text7324) at (5.2465, 4.0481){c};

  \path[draw=c00cc00,fill=c00cc00,line width=0.0cm] (6.8822, 6.0903) circle (0.1384cm);

  \node[text=white,anchor=south,line width=0.0cm] (text6896) at (6.8822, 5.9796){b};

  \path[draw=c009900,fill=c009900,line width=0.0cm] (7.1868, 6.0903) circle (0.1384cm);

  \node[text=white,anchor=south,line width=0.0cm] (text6527) at (7.1868, 5.9796){c};

\end{tikzpicture}}}
    \qquad
    {\scalefont{1}\scalebox{0.6}{\input{images/tikz/10_2_16lattice.tex}}}
    \caption{The \href{https://christopherlloyd.github.io/LooPindex/multiloops/10\%5E2_16.html}{multiloop $10^2_{16}$} and its pinning semi-lattice obtained by unions of minimal pinning sets (in red or green), together with the set of all regions.}
    \label{fig:some-pinning-ideal}
\end{figure}

\subsection{Relations with previous works}

\paragraph{Computing self-intersection}

The algorithmic computation for the self-intersection number of (multi)loops in various surfaces has been approached from several perspectives (geometric group theory and combinatorial topology).

Our Algorithm \ref{algo:Intro:self-intersection} is not entirely new as it adapts \cite{Birman-Series_algo-simple-closed-curves_1984, Birman-Series_Dehn-algo-simple-closed-curves_1987, Cohen-Lustig_algo-geometric-intersection-number_1987, Lustig_Paths-geodesics-geometric-intersection-II_1987}, but our discussion does explain some details not covered in those references (such as dealing with non-primitive loops, as highlighted in \cite{Despres-Lazarus_Computing-intersection_2019}).
Moreover, its complexity is less tight than the one proposed and analyzed in \cite{Despres-Lazarus_Computing-intersection_2019} which relies on different methods, but our aim here is to present an algorithm that has been implemented to suit our pinning problems with enough detail so as to easily justify its correctness while controlling its complexity.
Indeed, we found it both robust and versatile as it combines with dynamical programming approaches to pinning problems which evolve in time (when removing and adding pins, or modifying the isotopy class of the multiloop).

\paragraph{Tautening multiloops}

\begin{figure}
    \centering
    \includegraphics[scale=0.15]{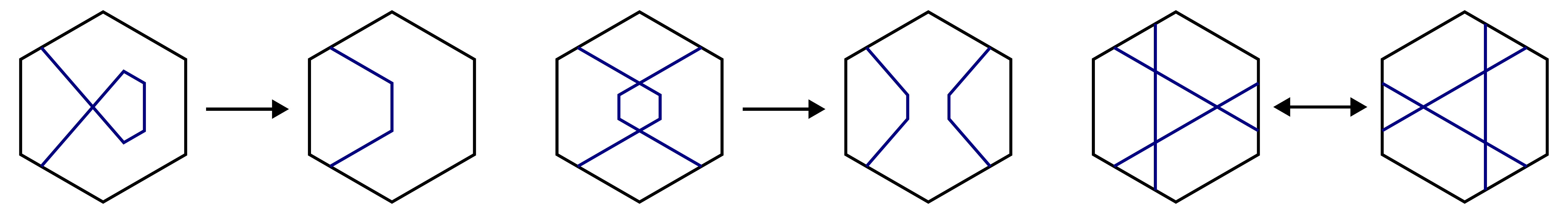}
    \caption{The shortening Reidemeister moves of type $\Vec{R1}$, $\Vec{R2}$, $R3$.}
    \label{fig:Reidemeister}
\end{figure}

In a surface $\Surf_P$, if two multiloops are homotopic then they are connected by a sequence of isotopies and Reidemeister moves $R1, R2, R3$ (depicted in Figure \ref{fig:Reidemeister}), and one may define the \emph{combinatorial length} as the total (weighted) number of Reidemeister moves (for a certain choice of weights).
Moreover, a multiloop can be connected to a taut representative by a sequence of isotopies and \emph{shortening}-Reidemeister moves $\Vec{R1}, \Vec{R2}, R3$ (which never increase the number of double points).
The existence of such shortening homotopies was shown for loops by \cite{Hass-Scott_Shortening-curves_1994} and made explicit for multiloops in terms of Reidemeister moves by \cite{Graaf-Schrijver_Reidemeister-shortening_1997}.

Given a multiloop in a surface, the work \cite{Chang-Mesmay_Reidemeister-shortening_2022} demonstrate an algorithm of polynomial complexity which constructs a curve shortening homotopy to some taut representative (in particular, the combinatorial length of a shortening homotopy is at most a polynomial).
On the other hand, the work \cite{Chang-Mesmay-al_Tightening-curves-local-moves_2018} construct families of loops  for which this combinatorial length is at least quadratic function.

\paragraph{Computing immersed monorbigons}

Our Theorem \ref{thm:Intro:pinning-mobidiscs} reveals that the problem of computing mobidiscs in closed surfaces may always be reduced to that of punctured surfaces in polynomial time by running over all possible punctures.

Given a loop in a surface, the description and computation of all its extensions to immersions of a disc was raised by Hopf and Thom \cite{Thom_inequadiff-globales_1959}.

When the surface is the plane, it was first solved by Blank \cite{Blank_extending-immersions-circle_1967, Poenaru_extension-immersions-codim-1_1995} who gave a combinatorial characterization by a polynomial-time algorithm.
Then Shor--Van Wyk \cite{Shor-VanWyk_Detecting-decomposing-self-overlapping-curves_1992} proposed another algorithmic solution of polynomial complexity.
This yields our polynomial Algorithm \ref{algo:Intro:mobidisc} for loops in the sphere.

In higher genus however, the characterization and computation of loops that bound immersed discs remains unclear.
Frisch worked on generalizing the methods of Blank to oriented surfaces of higher genus in his Thesis  \cite{Frisch_classification-immersions-curves-in-surfaces_2010} (and one may extract a polynomial algorithm from \cite[Theorem 3.1.5 and Theorem 3.2.4 and Theorem 5.2.7]{Frisch_classification-immersions-curves-in-surfaces_2010}), but restricted to the case of the sphere in his preprint \cite{Frisch_extending-immersions-sphere_2010}, and none of these works have been published.

More general problems have been approached by various methods, but it is not clear whether one may cunningly extract a polynomial algorithm from any of these.
For instance the work of \cite{Pappas_Extensions-codim-1-immersions_1996} uses Morse theory to describe codimension one extensions (in any dimension). 
The work \cite{Ezell-Marx_Branched-extensions-curves-orientable-surfaces_1980} applies methods to those of Blank and Shor--Van Wyk to describe analytic extensions to more general surfaces.

Let us finally mention that one may hope to adapt the polynomial algorithm in \cite{Chang-Mesmay_Reidemeister-shortening_2022} for tightening a multiloop in any surface to deduce an algorithm to compute mobidiscs of loops in any surface, but that remains to be done.

\subsection{Further directions of research}

The notion of pinning ideal and the problems related to its computation prompts improvement of past works about the combinatorial topology of multiloops and raises new challenges.

\paragraph{Multiloops versus loops.}

The characterizations of taut loops in terms of singular monorbigons and mobidiscs leading to Theorem \ref{thm:Intro:pinning-mobidiscs} do not have any analogs for multiloops.

\begin{question}[characterizing taut multiloops]
    Can one find a topological notion generalizing mobidiscs to multiloops so as to characterize which ones are not taut?
\end{question}

For the purposes of pinning problems, Theorem \ref{thm:Intro:pinning-mobidiscs} implies that the pinning ideal of a loop $\gamma$ with regions $\Regions$ is isomorphic to the solution ideal of the monotone conjunctive normal form associated to $\MoB(\gamma)\subset \Parts(\Regions)$.
The pinning ideal of multiloops remains mysterious.

\begin{question}[basis of the pinning ideal]
    Given a multiloop, can one compute a basis for its pinning ideal in polynomial time?
\end{question}

The following aims at qualitative differences between pinning ideals of loops and multiloops.

\begin{question}[counting strands]
    Can we compute or estimate lower bounds for the number of strands of a multiloop from its pinning ideal?
\end{question}

Let us focus on a ``simple'' family of multiloops. 
A loop has one strand, so all double points involve that single strand.
At the other extreme are the multiloops whose double points only involve distinct strands, which we call \emph{simple multiloops} as they consists of simple strands.

\begin{question}[computing the pinning number of a simple multiloops]
    We may restrict the \textsc{MultiLooPinNum} problem to simple multiloops.
    It can be approached using methods similar to those in subsection \ref{subsec:monorbigons}, relying on \cite[Proposition 1.7]{Farb-Margalit_MappingClassGroups_2012}.
    Indeed, one may show that the problem is in \textsf{P} if we restrict to $s\le 3$ strands, and \textsf{NP}-complete if we do not restrict $s\in \N$.
    Is there a critical value of $s\in \N$ for which it is \textsf{NP}-complete?
\end{question}

\paragraph{Approximate solutions and heuristics.}

When facing an \textsf{NP}-complete problem, one wonders whether it is possible to find approximate solutions in polynomial time.

\begin{question}[approximability]
    Do there exist $\epsilon_1 ,\epsilon_2 >0$ and a polynomial time algorithm which: given a multiloop $\gamma$ and $k\in \N$, either confirms that $\varpi(\gamma) > k$ or else constructs $P\subset \Regions$ with $\Card P \le (1+\epsilon_1) k$ and $\si_P(\gamma) > (1-\epsilon_2). \si_\Regions(\gamma)$?
\end{question}

A weaker expectation would be an algorithm based on certain heuristics which computes an approximate solution with high probability (that is on most entries).

\begin{question}[heuristics]
    Are there heuristics building on the idea of pinning regions of small degree and certain embedded monorbigons which lead to algorithms that find approximate solutions with high probability?
    We discuss some in section \ref{sec:pinning-ideal}.
\end{question}

\paragraph{Pinning statistics for random multiloops.}

The statics of pinning ideals of random multiloops would shed light on the distribution of hard instances to certain pinning problems, hence guide the search (or discourage the hope) for approximate solutions and heuristics.

\begin{question}
    Consider a random model for unoriented multiloops in closed orientable surfaces (in fixed genus or not), by choosing a distribution on a set of $4$-regular maps whose weights are inverse proportional to cardinal of the automorphism groups.
    What is the distribution of:
\begin{itemize}[noitemsep]
    \item[] the size of an optimal pinning set, that is the pinning number?
    \item[] the number of optimal pinning sets, and of minimal pinning sets?
\end{itemize}
\end{question}

\paragraph{Other variational problems.}

For a filling multiloop $\gamma \colon \sqcup_1^s \Sph^1 \looparrowright \Surf$, consider the increasing functionals on its lattice of regions $\Parts(\Regions)$, which to a set $P\in \Regions$ associate: 
\begin{itemize}[noitemsep]
    \item its cardinal $\Card(P) = \sum_{P_i\in P} 1$, its \emph{total-degree} $\operatorname{Deg}(P) = \sum_{P_i\in P} \deg(P_i)$;
    \item minus the self-intersection number $-\si_P(\gamma)$ of the multicurve $\gamma \colon \sqcup_1^s \Sph^1 \looparrowright \Surf_P$.
\end{itemize}
Note that $\Card(P)$ and $\deg(P)$ are valuations in the sense of lattice theory.

\begin{question}[multiloops with the same self-intersection functionals]
    We have many examples of loops in the sphere with isomorphic pinning ideals.

    More precisely, we may say that two multiloops $\alpha, \beta$ to be $\si$-equivalent when their self-intersection functionals $\Parts(\Regions(\alpha)) \ni P\mapsto \si_P(\alpha) \in \N$ and $\Parts(\Regions(\beta)) \ni P\mapsto \si_P(\beta) \in \N$ are conjugate by a bijection $\varphi \colon \Regions(\alpha) \to \Regions(\beta)$.
    Can we describe the $\si$-equivalence classes ?
\end{question}

The pinning problems addressed in this work concern the local and global minima in $\Parts(\Regions)$ of the valuation $\Card(P)$ in restriction to the level-set $\{P\in \Regions \mid \si_P(\gamma)=n\}$ where $n$ is the number of double points of the multiloop $\gamma$.

If we change $\Card(P)$ for $\operatorname{Deg}(P)$, then we may apply Algorithm \ref{algo:Intro:self-intersection} to deduce an analogous Corollary \ref{cor:intro:LooPinMin-is-P} showing that certain pinning problems are in \textsf{P}.

\begin{question}[Changing the valuation]
    We believe that the \textsc{LooPinDeg} problem, that is the $\operatorname{Deg}$-weighted analogue of the \textsc{LooPinNum}, problem is also \textsf{NP}-complete.
    
    Can one adapt the Algorithm \ref{algo:intro:planar-vertex-cover-to-LooPinNum} to reduce the problem of finding a vertex cover for a plane graph to that of finding the minimum total region degree of a corresponding plane loop?
\end{question}

\begin{question}[extremizing other functionals]
There are other interesting functionals to extremize, which from a game-theoretic perspective can be thought of as varying the cost or gain function, and that may yield other interesting topological invariants of filling multiloops.

Which monotonous functions depending on $(\Card(P), \deg(P))$ and $-\si_P(\gamma)$ have minimum whose computation is in \textsf{P}? 
What are the topological interpretations of these quantities?

In particular, what is the complexity of computing the minimal values of the monotonous functionals $\Card(P)-\si_P(\gamma)$ or $\deg(P)-\si_P(\gamma)$?
What about non-monotonous functionals of the form $4\Card(P)-\deg(P)- 3\si_P(\gamma)$?
\end{question}

\begin{remark}[taut multiloops are shortest geodesics]
It follows from \cite{Neumann-Coto_shortest-geodesics-surfaces_2001} that a multiloop is taut if and only if it is isotopic to a union of shortest geodesics for some Riemannian metric. We emphasize that the term ``shortest'' is important here.

However, the multiloops arising from geodesics in negatively curved or hyperbolic surfaces are much more subtle to describe (see for instance \cite{Hass-Scott_configurations-curves_1999}).
\end{remark}

\paragraph{Pinning invariants for homeomorphisms.}

Let $\psi \in \operatorname{Homeo}_0(\Surf)$ be a homeomorphism which is isotopic to the identity.
For a multiloop $\gamma \colon \sqcup_1^s \Sph^1 \looparrowright \Surf$, we define
\begin{equation*}
    \varpi(\psi,\gamma) = \min \{\varpi(\gamma_0 \cup \psi(\gamma_0)) \mid \text{generic immersions $\gamma_0$ in the isotopy class $\gamma$} \}
\end{equation*}
This quantity measures the minimum cardinal of a puncture-set $P\subset \Surf$ yielding a topological model for the action of a class $[\psi]\in \operatorname{Map}(\Surf_P)$ on $\gamma \subset \Surf_P$.
It may be compared with $\varpi(\gamma)+\varpi(\psi(\gamma))=2\varpi(\gamma)$
Now, for a certain class of multiloops $\mathcal{C}$ (say all multiloops, or those without double points), we define 
\begin{equation*}
    \varpi(\psi,\mathcal{C})=\min\{\varpi(\psi,\gamma)\mid \gamma \in \mathcal{C}\}.
\end{equation*}

The case where $\mathcal{C}$ are loops with no double points appears closely related to the work \cite{Bowden-Hensel-Webb_quasi-morphisms-Diff-0_2022} about quasi-morphisms on the groups $\operatorname{Homeo}_0(\Surf)$ and $\operatorname{Diff}_0(\Surf)$.
Indeed, their main tool is to relate the action of $\psi \in \operatorname{Diff}_0(\Surf)$ on the fine curve graph of $\Surf$ to the actions of $[\psi]\in \operatorname{Map}(\Surf_P)$ on the surviving curve graphs of $\Surf_P$.
We believe that the quantities $\varpi(\psi,\mathcal{C})$ may be of interest for such investigations.

\section{Group theory: counting self-intersections}
\label{sec:combinatorial-group-theory-self-intersection}

In Subsection \ref{subsec:maps-permutations} we encode nontrivial filling multiloops $\gamma\colon \sqcup_1^s \Sph^1 \looparrowright \Surf$ by $4$-regular maps. 
We used this encoding extensively in the implementation of our algorithms, but one may gloss over it without much hindrance to understanding the rest of our work.

In Subsection \ref{subsec:MuLoops-to-FreeGroups} we use this encoding to derive, for non-empty $P\subset \Regions$, a presentation of the fundamental group $\pi_1(\Surf_P)$ and the conjugacy classes associated to the strands of $\gamma$.

In Subsection \ref{subsec:si-algorithm} we explain an algorithm to compute, for non-empty $P\subset \Regions$, the self-intersection number of the multicurve associated to $\gamma$ in $\Surf_P$.

Altogether they imply Corollary \ref{cor:MultiLooPinMin-P-easy}, saying certain pinning problems are in \textsf{P}.

\subsection{Encoding filling multiloops as \texorpdfstring{$4$-regular}{4-valent} maps}
\label{subsec:maps-permutations}

Let us first explain the permutation encoding of filling multigraphs and filling multiloops, which we will use extensively in our algorithms.
A general reference for this is \cite{Lando-Zvonkin_graphs-surfaces_2004}.

A \emph{map} $\Gamma\subset \Surf$ is an embedded multigraph in $\Surf$ whose complement is homotopic to a disjoint union of discs, considered up to orientation-preserving diffeomorphisms.

\begin{lemma}
A map $\Gamma \subset \Surf$ with $n$ edges corresponds to a pair of permutations $\epsilon, \sigma \in \mathfrak{S}_{2n}$ considered up to conjugacy, such that:
\begin{itemize}[noitemsep]
    \item[-] $\epsilon$ has orbits of size $2$ corresponding to the edges
    \item[-] $\sigma$ has orbits corresponding to the vertices
    \item[-] $\varphi = (\sigma \epsilon)^{-1}$ has orbits corresponding to the regions
\end{itemize}
The \emph{degree} of a vertex or region is the cardinal of the corresponding orbit for $\sigma$ or $\varphi$.
The Euler characteristic of $\Surf$ is the alternating sum of the number of orbits of $\sigma,\epsilon, \varphi$.
The dual map of $\Gamma\subset \Surf$ is obtained by exchanging $\sigma$ and $\varphi$.

\end{lemma}

\begin{remark}[Computing regions]
    \label{rem:compute-union-regions}
    For a map $\epsilon,\sigma \in \mathfrak{S}_{2n}$ and a union of orbits of $\varphi=(\epsilon\sigma)^{-1}$, the corresponding union of closed regions and its interior have homotopy types which may be computed in polynomial time on $n$.
    For example, we often wish to determine their number of connected components and their genera. 
\end{remark}

A filling multiloop $\gamma\colon \sqcup_1^s \Sph^1 \looparrowright \Surf$ is determined, modulo orientation of its strands, by its associated $4$-regular map $\Gamma \subset \Surf$.
Note that the orientation of $\gamma$ corresponds to an orientation of the edges of $\Gamma$ which go straight on at each vertex.

\begin{lemma}
A $4$-regular map $\Gamma \subset \Surf$ with $n$ edges corresponds to a pair of permutations $\epsilon, \sigma \in \mathfrak{S}_{2n}$ considered up to conjugacy, such that:
\begin{itemize}[noitemsep]
    \item[-] $\epsilon$ has orbits of size $2$ corresponding to the edges
    \item[-] $\sigma$ has orbits of size $4$ corresponding to the vertices
    \item[-] $\varphi = (\sigma \epsilon)^{-1}$ has orbits corresponding to the regions
\end{itemize}
Moreover, the orbits of $\delta = \epsilon \sigma^2$ define the oriented strands of $\Gamma$, and these are paired by $\epsilon$ to form its unoriented strands. 
If a multiloop $\gamma$ yields the map $\Gamma$, the oriented strands of $\gamma$ select one orbit of $\delta$ in each $\epsilon$-pair.
\end{lemma}

\begin{example}[permutation representation]
    In Figure \ref{fig:4-valent-map-permutations}, the permutations $\epsilon, \sigma, \varphi, \delta$ are:
    \small{
	\begin{equation*}
    \label{eq:permrepexampleagain}
		\begin{split}
			\epsilon= &(-1,1)(-2,2)(-3,3)(-4,4)(-5,5)(-6,6)(-7,7)(-8,8)\\
			 &(-9,9)(-10,10)(-11,11)(-12,12)(-13,13)(-14,14)(-15,15)(-16,16), \\
			\sigma=&(9,8,-10,-1)(5,2,-6,-3)(3,12,-4,-13)(13,4,-14,-5)\\
			&(14,7,-15,-8)(10,15,-11,-16)(1,16,-2,-9)(6,11,-7,-12), \\
			\varphi=&(-1,-9)(-2,5,-14,-8,9)(-3,-13,-5)(-4,13)(-6,-12,3)\\
			&(-7,14,4,12)(-10,-16,1)(-11,6,2,16)(-15,10,8)(7,11,15), \\
			\delta=&(1,2,3,4,5,6,7,8)(-8,-7,-6,-5,-4,-3,-2,-1)\\
			&(9,10,11,12,13,14,15,16)(-16,-15,-14,-13,-12,-11,-10,-9).
		\end{split}
	\end{equation*}
    }
    \begin{figure}[H]
    	\centering
    	\scalebox{0.9}{\definecolor{cd82626}{RGB}{216,38,38}

\def \globalscale {1.000000}
\begin{tikzpicture}[y=1cm, x=1cm, yscale=\globalscale,xscale=\globalscale, every node/.append style={scale=\globalscale}, inner sep=0pt, outer sep=0pt]
  \path[draw=cd82626,fill=cd82626,line width=0.1588cm] (4.115, 4.4362) -- (4.115, 0.5601);

  \path[draw=cd82626,fill=cd82626,line width=0.1588cm] (1.381, 4.4362) -- (1.381, 0.5601);

  \path[draw=cd82626,fill=cd82626] (-0.8768, 2.4981) circle (0.0529cm);

  \path[draw=cd82626,fill=cd82626] (6.6476, 2.4981) circle (0.0529cm);

  \path[draw=cd82626,fill=cd82626,line width=0.1588cm] (6.6476, 2.4981) -- (-0.8768, 2.4981);

  \path[draw=black,fill=white,opacity=0.8,line width=0.0132cm] (5.3321, 2.4981) rectangle (6.5489, 1.5398);

  \tikzstyle{every node}=[font=\fontsize{10}{10}\selectfont]
  \node[text=black,anchor=south east,line width=0.0cm] (text27) at (6.6318, 1.5855){\begin{tabular}{l} $\delta(e)=$\\ $\epsilon\sigma^2(e)$ \end{tabular}};

  \path[draw=black,fill=white,opacity=0.8,line width=0.0132cm] (4.115, 2.2189) rectangle (4.9277, 1.6944);

  \node[text=black,line width=0.0cm,anchor=south west] (text32) at (4.1944, 1.7824){$\sigma(e)$};

  \path[draw=black,fill=white,opacity=0.8,line width=0.0132cm] (1.381, 2.2189) rectangle (3.0559, 1.2953);

  \node[text=black,line width=0.0cm,anchor=south west] (text34) at (1.2399, 1.3398){\begin{tabular}{l} $\varphi(e)=$\\ $(\sigma\epsilon)^{-1}(e)$ \end{tabular}};

  \path[draw=black,fill=white,opacity=0.8,line width=0.0132cm] (3.2948, 2.4981) rectangle (3.6767, 2.1542);

  \node[text=black,anchor=south east,line width=0.0cm] (text47) at (3.5476, 2.2667){$e$};

  \path[draw=black,fill=white,opacity=0.8,line width=0.0132cm] (4.368, 3.0131) rectangle (5.2785, 2.4981);

  \node[text=black,line width=0.0cm,anchor=south west] (text48) at (4.438, 2.5644){$\sigma^2(e)$};

  \path[draw=black,fill=white,opacity=0.8,line width=0.0132cm] (1.6808, 3.002) rectangle (2.4275, 2.4981);

  \node[text=black,line width=0.0cm,anchor=south west] (text58) at (1.7569, 2.5644){$\epsilon(e)$};

  \path[draw=cd82626,fill=cd82626] (4.115, 4.4362) circle (0.0529cm);

  \path[draw=cd82626,fill=cd82626] (4.115, 0.5601) circle (0.0529cm);

  \path[draw=cd82626,fill=cd82626] (1.381, 4.4362) circle (0.0529cm);

  \path[draw=cd82626,fill=cd82626] (1.381, 0.5601) circle (0.0529cm);

\end{tikzpicture}}
    	\scalebox{0.6}{\definecolor{cd82626}{RGB}{216,38,38}
\definecolor{c2626d8}{RGB}{38,38,216}

\def \globalscale {1.000000}
\begin{tikzpicture}[y=1cm, x=1cm, yscale=\globalscale,xscale=\globalscale, every node/.append style={scale=\globalscale}, inner sep=0pt, outer sep=0pt]

 \tikzstyle{every node}=[font=\fontsize{10}{10}\selectfont]
  \path[draw=cd82626,fill=cd82626] (6.6146, 8.296) circle (0.0529cm);

  \path[draw=cd82626,fill=cd82626] (10.4951, 8.296) circle (0.0529cm);

  \path[draw=cd82626,fill=cd82626] (10.4951, 0.5437) circle (0.0529cm);

  \path[draw=cd82626,fill=cd82626] (4.6743, 0.5437) circle (0.0529cm);

  \path[draw=cd82626,fill=cd82626] (4.6743, 8.296) circle (0.0529cm);

  \path[draw=cd82626,fill=cd82626] (0.7938, 8.296) circle (0.0529cm);

  \path[draw=cd82626,fill=cd82626] (0.7938, 2.4818) circle (0.0529cm);

  \path[draw=cd82626,fill=cd82626] (6.6146, 2.4818) circle (0.0529cm);

  \path[draw=c2626d8,fill=c2626d8] (2.734, 4.4199) circle (0.0529cm);

  \path[draw=c2626d8,fill=c2626d8] (2.734, 6.3579) circle (0.0529cm);

  \path[draw=c2626d8,fill=c2626d8] (8.5549, 6.3579) circle (0.0529cm);

  \path[draw=c2626d8,fill=c2626d8] (8.5549, 2.4818) circle (0.0529cm);

  \path[draw=c2626d8,fill=c2626d8] (12.4354, 2.4818) circle (0.0529cm);

  \path[draw=c2626d8,fill=c2626d8] (12.4354, 4.4199) circle (0.0529cm);

  \path[draw=cd82626,line width=0.1588cm] (6.6146, 8.296) -- (10.4951, 8.296);

  \path[draw=cd82626,line width=0.1588cm] (10.4951, 8.296) -- (10.4951, 0.5437);

  \path[draw=cd82626,line width=0.1588cm] (10.4951, 0.5437) -- (4.6743, 0.5437);

  \path[draw=cd82626,line width=0.1588cm] (4.6743, 0.5437) -- (4.6743, 8.296);

  \path[draw=cd82626,line width=0.1588cm] (4.6743, 8.296) -- (0.7938, 8.296);

  \path[draw=cd82626,line width=0.1588cm] (0.7938, 8.296) -- (0.7938, 2.4818);

  \path[draw=cd82626,line width=0.1588cm] (0.7938, 2.4818) -- (6.6146, 2.4818);

  \path[draw=cd82626,line width=0.1588cm] (6.6146, 2.4818) -- (6.6146, 8.296);

  \path[draw=c2626d8,line width=0.1588cm] (2.734, 4.4199) -- (2.734, 6.3579);

  \path[draw=c2626d8,line width=0.1588cm] (2.734, 6.3579) -- (8.5549, 6.3579);

  \path[draw=c2626d8,line width=0.1588cm] (8.5549, 6.3579) -- (8.5549, 2.4818);

  \path[draw=c2626d8,line width=0.1588cm] (8.5549, 2.4818) -- (12.4354, 2.4818);

  \path[draw=c2626d8,line width=0.1588cm] (12.4354, 2.4818) -- (12.4354, 4.4199);

  \path[draw=c2626d8,line width=0.1588cm] (12.4354, 4.4199) -- (2.734, 4.4199);

  \path[fill=white,opacity=0.8,line width=0.0cm] (7.7258, 4.4581) rectangle (8.255, 4.1141);

  \node[text=black,anchor=south east,line width=0.0cm] (text27) at (8.2318, 4.1804){-11};

  \path[fill=white,opacity=0.8,line width=0.0cm] (8.8635, 4.8021) rectangle (9.2869, 4.4581);

  \node[text=black,line width=0.0cm,anchor=south west] (text28) at (8.8779, 4.5244){10};

  \path[fill=white,opacity=0.8,line width=0.0cm] (10.4775, 2.2091) rectangle (10.8215, 1.8652);

  \node[text=black,line width=0.0cm,anchor=south west] (text29) at (10.4951, 1.9315){-2};

  \path[fill=white,opacity=0.8,line width=0.0cm] (10.2923, 3.1881) rectangle (10.5304, 2.8441);

  \node[text=black,anchor=south east,line width=0.0cm] (text30) at (10.4951, 2.9104){1};

  \path[fill=white,opacity=0.8,line width=0.0cm] (6.2971, 5.1196) rectangle (6.641, 4.7756);

  \node[text=black,anchor=south east,line width=0.0cm] (text31) at (6.6146, 4.8419){-7};

  \path[fill=white,opacity=0.8,line width=0.0cm] (6.5881, 4.1406) rectangle (6.8263, 3.7966);

  \node[text=black,line width=0.0cm,anchor=south west] (text32) at (6.6146, 3.8629){6};

  \path[fill=white,opacity=0.8,line width=0.0cm] (4.3656, 5.1196) rectangle (4.7096, 4.7756);

  \node[text=black,anchor=south east,line width=0.0cm] (text33) at (4.6743, 4.8419){-4};

  \path[fill=white,opacity=0.8,line width=0.0cm] (4.6567, 4.1406) rectangle (4.8948, 3.7966);

  \node[text=black,line width=0.0cm,anchor=south west] (text34) at (4.6743, 3.8629){3};

  \path[fill=white,opacity=0.8,line width=0.0cm] (4.3656, 7.0775) rectangle (4.7096, 6.7335);

  \node[text=black,anchor=south east,line width=0.0cm] (text35) at (4.6743, 6.7998){-5};

  \path[fill=white,opacity=0.8,line width=0.0cm] (4.6567, 6.0721) rectangle (4.8948, 5.7281);

  \node[text=black,line width=0.0cm,anchor=south west] (text36) at (4.6743, 5.7944){4};

  \path[fill=white,opacity=0.8,line width=0.0cm] (4.9742, 6.7335) rectangle (5.5033, 6.3896);

  \node[text=black,line width=0.0cm,anchor=south west] (text37) at (4.9973, 6.4558){-14};

  \path[fill=white,opacity=0.8,line width=0.0cm] (3.9423, 6.3896) rectangle (4.3656, 6.0456);

  \node[text=black,anchor=south east,line width=0.0cm] (text38) at (4.3513, 6.1119){13};

  \path[fill=white,opacity=0.8,line width=0.0cm] (6.9056, 6.7335) rectangle (7.4348, 6.3896);

  \node[text=black,line width=0.0cm,anchor=south west] (text39) at (6.9376, 6.4558){-15};

  \path[fill=white,opacity=0.8,line width=0.0cm] (5.9002, 6.3896) rectangle (6.3235, 6.0456);

  \node[text=black,anchor=south east,line width=0.0cm] (text40) at (6.2916, 6.1119){14};

  \path[fill=white,opacity=0.8,line width=0.0cm] (4.9742, 2.8706) rectangle (5.3181, 2.5266);

  \node[text=black,line width=0.0cm,anchor=south west] (text41) at (4.9973, 2.5929){-6};

  \path[fill=white,opacity=0.8,line width=0.0cm] (4.1275, 2.5266) rectangle (4.3656, 2.1827);

  \node[text=black,anchor=south east,line width=0.0cm] (text42) at (4.3513, 2.249){5};

  \path[fill=white,opacity=0.8,line width=0.0cm] (8.5196, 4.1406) rectangle (9.0488, 3.7966);

  \node[text=black,line width=0.0cm,anchor=south west] (text43) at (8.5549, 3.8629){-16};

  \path[fill=white,opacity=0.8,line width=0.0cm] (8.1492, 5.1196) rectangle (8.5725, 4.7756);

  \node[text=black,anchor=south east,line width=0.0cm] (text44) at (8.5549, 4.8419){15};

  \path[fill=white,opacity=0.8,line width=0.0cm] (10.4775, 4.1406) rectangle (10.8215, 3.7966);

  \node[text=black,line width=0.0cm,anchor=south west] (text45) at (10.4951, 3.8629){-1};

  \path[fill=white,opacity=0.8,line width=0.0cm] (10.2923, 5.1196) rectangle (10.5304, 4.7756);

  \node[text=black,anchor=south east,line width=0.0cm] (text46) at (10.4951, 4.8419){8};

  \path[fill=white,opacity=0.8,line width=0.0cm] (5.7944, 4.4581) rectangle (6.3235, 4.1141);

  \node[text=black,anchor=south east,line width=0.0cm] (text47) at (6.2916, 4.1804){-12};

  \path[fill=white,opacity=0.8,line width=0.0cm] (6.9056, 4.8021) rectangle (7.329, 4.4581);

  \node[text=black,line width=0.0cm,anchor=south west] (text48) at (6.9376, 4.5244){11};

  \path[fill=white,opacity=0.8,line width=0.0cm] (10.795, 2.8706) rectangle (11.139, 2.5266);

  \node[text=black,line width=0.0cm,anchor=south west] (text49) at (10.8182, 2.5929){-9};

  \path[fill=white,opacity=0.8,line width=0.0cm] (9.7631, 2.5266) rectangle (10.1865, 2.1827);

  \node[text=black,anchor=south east,line width=0.0cm] (text50) at (10.1721, 2.249){16};

  \path[fill=white,opacity=0.8,line width=0.0cm] (4.3656, 3.1881) rectangle (4.7096, 2.8441);

  \node[text=black,anchor=south east,line width=0.0cm] (text51) at (4.6743, 2.9104){-3};

  \path[fill=white,opacity=0.8,line width=0.0cm] (4.6567, 2.2091) rectangle (4.8948, 1.8652);

  \node[text=black,line width=0.0cm,anchor=south west] (text52) at (4.6743, 1.9315){2};

  \path[fill=white,opacity=0.8,line width=0.0cm] (9.6573, 4.4581) rectangle (10.1865, 4.1141);

  \node[text=black,anchor=south east,line width=0.0cm] (text53) at (10.1721, 4.1804){-10};

  \path[fill=white,opacity=0.8,line width=0.0cm] (10.795, 4.8021) rectangle (11.0331, 4.4581);

  \node[text=black,line width=0.0cm,anchor=south west] (text54) at (10.8182, 4.5244){9};

  \path[fill=white,opacity=0.8,line width=0.0cm] (6.2971, 7.0775) rectangle (6.641, 6.7335);

  \node[text=black,anchor=south east,line width=0.0cm] (text55) at (6.6146, 6.7998){-8};

  \path[fill=white,opacity=0.8,line width=0.0cm] (6.5881, 6.0721) rectangle (6.8263, 5.7281);

  \node[text=black,line width=0.0cm,anchor=south west] (text56) at (6.6146, 5.7944){7};

  \path[fill=white,opacity=0.8,line width=0.0cm] (3.8365, 4.4581) rectangle (4.3656, 4.1141);

  \node[text=black,anchor=south east,line width=0.0cm] (text57) at (4.3513, 4.1804){-13};

  \path[fill=white,opacity=0.8,line width=0.0cm] (4.9742, 4.8021) rectangle (5.3975, 4.4581);

  \node[text=black,line width=0.0cm,anchor=south west] (text58) at (4.9973, 4.5244){12};

\end{tikzpicture}}
    	\caption{Local action of the permutations $\epsilon, \sigma, \varphi, \delta$. A map with labeled half-edges.}
    	\label{fig:4-valent-map-permutations}
    \end{figure}
\end{example}

The following observation will serve only to find all (multi)loops in the sphere whose regions all have degree $\ge 3$, given their total number of regions.

\begin{lemma}
    \label{lem:average-degree-regions}
    For a $4$-regular map $\Gamma \subset \Surf$ in a surface of Euler characteristic $\chi$, the degree of the regions $R_i\in\Regions$ satisfy 
    \begin{equation*}\textstyle
        \sum_{R_i\in \Regions} (\deg(R_i)-4) = -4\chi.
    \end{equation*}
    In particular, a multiloop in the torus has average region-degree $4$, and a multiloop in the sphere whose regions have degree $\ge 3$ has at least $8$ regions of degree $3$.
\end{lemma}

\begin{proof}
    Denote the vertices by $V$, the edges by $E$ and the regions by $\Regions$.
    The $4$-valency of the map implies that $2 \rvert E\rvert = 4 \lvert V\rvert$.
    The surface embedding yields $\sum_\Regions \deg R_i = 2\lvert E\rvert = 4 \lvert V\rvert$.
    The filling property computes the Euler characteristic $\chi =\lvert V \rvert - \lvert E \rvert+ \lvert \Regions \rvert$ so that $\lvert \Regions \rvert = \chi+\lvert V \rvert$.
    Combine these, we find that for all $m\in \N$ we have $\sum_\Regions (\deg(R_i)-m) = (4-m)\lvert V\rvert-m\chi$, whence $\sum_\Regions (\deg R_i-4) = -4\chi$.
\end{proof}

\begin{remark}[Irreducible, indecomposible]
    Following \cite{Coquereaux-Zuber_maps_immersions_permutations_2024}, a plane map is \emph{irreducible} when it is not $2$-connected (disconnected by removing a vertex) and \emph{indecomposible} when it is not $3$-edge-connected (disconnected by removing two distinct edges). 
    
    These notions are the respective analogues of a link diagram having no nugatory crossings and being prime.
\end{remark}

\begin{figure}[h]
    \centering
    {\scalefont{1}\scalebox{0.35}{\definecolor{cd82626}{RGB}{216,38,38}
\definecolor{c2626d8}{RGB}{38,38,216}
\definecolor{c26d826}{RGB}{38,216,38}

\def \globalscale {1.000000}
\begin{tikzpicture}[y=1cm, x=1cm, yscale=\globalscale,xscale=\globalscale, every node/.append style={scale=\globalscale}, inner sep=0pt, outer sep=0pt]
  \path[draw=cd82626,fill=cd82626] (12.4301, 8.1809) circle (0.0529cm);

  \path[draw=cd82626,fill=cd82626] (12.4301, 1.1959) circle (0.0529cm);

  \path[draw=cd82626,fill=cd82626] (0.7885, 1.1959) circle (0.0529cm);

  \path[draw=cd82626,fill=cd82626] (0.7885, 8.1809) circle (0.0529cm);

  \path[draw=c2626d8,fill=c2626d8] (5.4451, 12.8376) circle (0.0529cm);

  \path[draw=c2626d8,fill=c2626d8] (10.1018, 12.8376) circle (0.0529cm);

  \path[draw=c2626d8,fill=c2626d8] (10.1018, 3.5243) circle (0.0529cm);

  \path[draw=c2626d8,fill=c2626d8] (5.4451, 3.5243) circle (0.0529cm);

  \path[draw=c26d826,fill=c26d826] (7.7735, 5.8526) circle (0.0529cm);

  \path[draw=c26d826,fill=c26d826] (3.1168, 5.8526) circle (0.0529cm);

  \path[draw=c26d826,fill=c26d826] (3.1168, 10.5092) circle (0.0529cm);

  \path[draw=c26d826,fill=c26d826] (7.7735, 10.5092) circle (0.0529cm);

  \path[draw=cd82626,line width=0.1588cm] (12.4301, 8.1809) -- (12.4301, 1.1959);

  \path[draw=cd82626,line width=0.1588cm] (12.4301, 1.1959) -- (0.7885, 1.1959);

  \path[draw=cd82626,line width=0.1588cm] (0.7885, 1.1959) -- (0.7885, 8.1809);

  \path[draw=cd82626,line width=0.1588cm] (0.7885, 8.1809) -- (12.4301, 8.1809);

  \path[draw=c2626d8,line width=0.1588cm] (5.4451, 12.8376) -- (10.1018, 12.8376);

  \path[draw=c2626d8,line width=0.1588cm] (10.1018, 12.8376) -- (10.1018, 3.5243);

  \path[draw=c2626d8,line width=0.1588cm] (10.1018, 3.5243) -- (5.4451, 3.5243);

  \path[draw=c2626d8,line width=0.1588cm] (5.4451, 3.5243) -- (5.4451, 12.8376);

  \path[draw=c26d826,line width=0.1588cm] (7.7735, 5.8526) -- (3.1168, 5.8526);

  \path[draw=c26d826,line width=0.1588cm] (3.1168, 5.8526) -- (3.1168, 10.5092);

  \path[draw=c26d826,line width=0.1588cm] (3.1168, 10.5092) -- (7.7735, 10.5092);

  \path[draw=c26d826,line width=0.1588cm] (7.7735, 10.5092) -- (7.7735, 5.8526);

  \tikzstyle{every node}=[font=\fontsize{30}{30}\selectfont]
  \node[text=black,line width=0.0cm,anchor=south west] (text23) at (0.9966, 0.2117){$3$};

  \node[text=black,line width=0.0cm,anchor=south west] (text24) at (0.9966, 7.1967){$3$};

  \node[text=black,line width=0.0cm,anchor=south west] (text25) at (3.3249, 9.525){$3$};

  \node[text=black,line width=0.0cm,anchor=south west] (text26) at (3.3249, 7.1967){$3$};

  \node[text=black,line width=0.0cm,anchor=south west] (text27) at (5.6533, 9.525){$3$};

  \node[text=black,line width=0.0cm,anchor=south west] (text28) at (5.6533, 7.1967){$3$};

  \node[text=black,line width=0.0cm,anchor=south west] (text29) at (5.6533, 11.8533){$3$};

  \node[text=black,line width=0.0cm,anchor=south west] (text30) at (5.6533, 4.8683){$3$};

\end{tikzpicture}}}
    {\scalefont{1}\scalebox{0.35}{\definecolor{cd82626}{RGB}{216,38,38}

\def \globalscale {1.000000}
\begin{tikzpicture}[y=1cm, x=1cm, yscale=\globalscale,xscale=\globalscale, every node/.append style={scale=\globalscale}, inner sep=0pt, outer sep=0pt]
  \path[draw=cd82626,fill=cd82626] (3.1168, 7.5195) circle (0.0529cm);

  \path[draw=cd82626,fill=cd82626] (3.1168, 0.5345) circle (0.0529cm);

  \path[draw=cd82626,fill=cd82626] (7.7735, 0.5345) circle (0.0529cm);

  \path[draw=cd82626,fill=cd82626] (7.7735, 12.1761) circle (0.0529cm);

  \path[draw=cd82626,fill=cd82626] (0.7885, 12.1761) circle (0.0529cm);

  \path[draw=cd82626,fill=cd82626] (0.7885, 5.1911) circle (0.0529cm);

  \path[draw=cd82626,fill=cd82626] (12.4301, 5.1911) circle (0.0529cm);

  \path[draw=cd82626,fill=cd82626] (12.4301, 9.8478) circle (0.0529cm);

  \path[draw=cd82626,fill=cd82626] (5.4451, 9.8478) circle (0.0529cm);

  \path[draw=cd82626,fill=cd82626] (5.4451, 2.8628) circle (0.0529cm);

  \path[draw=cd82626,fill=cd82626] (10.1018, 2.8628) circle (0.0529cm);

  \path[draw=cd82626,fill=cd82626] (10.1018, 7.5195) circle (0.0529cm);

  \path[draw=cd82626,line width=0.1588cm] (3.1168, 7.5195) -- (3.1168, 0.5345);

  \path[draw=cd82626,line width=0.1588cm] (3.1168, 0.5345) -- (7.7735, 0.5345);

  \path[draw=cd82626,line width=0.1588cm] (7.7735, 0.5345) -- (7.7735, 12.1761);

  \path[draw=cd82626,line width=0.1588cm] (7.7735, 12.1761) -- (0.7885, 12.1761);

  \path[draw=cd82626,line width=0.1588cm] (0.7885, 12.1761) -- (0.7885, 5.1911);

  \path[draw=cd82626,line width=0.1588cm] (0.7885, 5.1911) -- (12.4301, 5.1911);

  \path[draw=cd82626,line width=0.1588cm] (12.4301, 5.1911) -- (12.4301, 9.8478);

  \path[draw=cd82626,line width=0.1588cm] (12.4301, 9.8478) -- (5.4451, 9.8478);

  \path[draw=cd82626,line width=0.1588cm] (5.4451, 9.8478) -- (5.4451, 2.8628);

  \path[draw=cd82626,line width=0.1588cm] (5.4451, 2.8628) -- (10.1018, 2.8628);

  \path[draw=cd82626,line width=0.1588cm] (10.1018, 2.8628) -- (10.1018, 7.5195);

  \path[draw=cd82626,line width=0.1588cm] (10.1018, 7.5195) -- (3.1168, 7.5195);

  \tikzstyle{every node}=[font=\fontsize{30}{30}\selectfont]
  \node[text=black,line width=0.0cm,anchor=south west] (text23) at (3.3249, 6.5352){$3$};

  \node[text=black,line width=0.0cm,anchor=south west] (text24) at (0.9966, 11.1919){$3$};

  \node[text=black,line width=0.0cm,anchor=south west] (text25) at (3.3249, 4.2069){$3$};

  \node[text=black,line width=0.0cm,anchor=south west] (text26) at (0.9966, 4.2069){$4$};

  \node[text=black,line width=0.0cm,anchor=south west] (text27) at (5.6533, 8.8635){$3$};

  \node[text=black,line width=0.0cm,anchor=south west] (text28) at (7.9816, 8.8635){$3$};

  \node[text=black,line width=0.0cm,anchor=south west] (text29) at (5.6533, 6.5352){$4$};

  \node[text=black,line width=0.0cm,anchor=south west] (text30) at (7.9816, 6.5352){$3$};

  \node[text=black,line width=0.0cm,anchor=south west] (text31) at (7.9816, 4.2069){$3$};

  \node[text=black,line width=0.0cm,anchor=south west] (text32) at (5.6533, 4.2069){$3$};

\end{tikzpicture}}}
    {\scalefont{1}\scalebox{0.35}{\definecolor{cd82626}{RGB}{216,38,38}

\def \globalscale {1.000000}
\begin{tikzpicture}[y=1cm, x=1cm, yscale=\globalscale,xscale=\globalscale, every node/.append style={scale=\globalscale}, inner sep=0pt, outer sep=0pt]
  \path[draw=cd82626,fill=cd82626] (12.4354, 8.8194) circle (0.0529cm);

  \path[draw=cd82626,fill=cd82626] (12.4354, 2.9986) circle (0.0529cm);

  \path[draw=cd82626,fill=cd82626] (2.734, 2.9986) circle (0.0529cm);

  \path[draw=cd82626,fill=cd82626] (2.734, 6.8792) circle (0.0529cm);

  \path[draw=cd82626,fill=cd82626] (8.5549, 6.8792) circle (0.0529cm);

  \path[draw=cd82626,fill=cd82626] (8.5549, 12.7) circle (0.0529cm);

  \path[draw=cd82626,fill=cd82626] (4.6743, 12.7) circle (0.0529cm);

  \path[draw=cd82626,fill=cd82626] (4.6743, 4.9389) circle (0.0529cm);

  \path[draw=cd82626,fill=cd82626] (10.4951, 4.9389) circle (0.0529cm);

  \path[draw=cd82626,fill=cd82626] (10.4951, 10.7597) circle (0.0529cm);

  \path[draw=cd82626,fill=cd82626] (6.6146, 10.7597) circle (0.0529cm);

  \path[draw=cd82626,fill=cd82626] (6.6146, 1.0583) circle (0.0529cm);

  \path[draw=cd82626,fill=cd82626] (0.7938, 1.0583) circle (0.0529cm);

  \path[draw=cd82626,fill=cd82626] (0.7938, 8.8194) circle (0.0529cm);

  \path[draw=cd82626,line width=0.1588cm] (12.4354, 8.8194) -- (12.4354, 2.9986);

  \path[draw=cd82626,line width=0.1588cm] (12.4354, 2.9986) -- (2.734, 2.9986);

  \path[draw=cd82626,line width=0.1588cm] (2.734, 2.9986) -- (2.734, 6.8792);

  \path[draw=cd82626,line width=0.1588cm] (2.734, 6.8792) -- (8.5549, 6.8792);

  \path[draw=cd82626,line width=0.1588cm] (8.5549, 6.8792) -- (8.5549, 12.7);

  \path[draw=cd82626,line width=0.1588cm] (8.5549, 12.7) -- (4.6743, 12.7);

  \path[draw=cd82626,line width=0.1588cm] (4.6743, 12.7) -- (4.6743, 4.9389);

  \path[draw=cd82626,line width=0.1588cm] (4.6743, 4.9389) -- (10.4951, 4.9389);

  \path[draw=cd82626,line width=0.1588cm] (10.4951, 4.9389) -- (10.4951, 10.7597);

  \path[draw=cd82626,line width=0.1588cm] (10.4951, 10.7597) -- (6.6146, 10.7597);

  \path[draw=cd82626,line width=0.1588cm] (6.6146, 10.7597) -- (6.6146, 1.0583);

  \path[draw=cd82626,line width=0.1588cm] (6.6146, 1.0583) -- (0.7938, 1.0583);

  \path[draw=cd82626,line width=0.1588cm] (0.7938, 1.0583) -- (0.7938, 8.8194);

  \path[draw=cd82626,line width=0.1588cm] (0.7938, 8.8194) -- (12.4354, 8.8194);

  \tikzstyle{every node}=[font=\fontsize{30}{30}\selectfont]
  \node[text=black,line width=0.0cm,anchor=south west] (text27) at (1.0289, 0.0794){$4$};

  \node[text=black,line width=0.0cm,anchor=south west] (text28) at (6.8498, 3.9688){$3$};

  \node[text=black,line width=0.0cm,anchor=south west] (text29) at (1.0289, 7.8317){$3$};

  \node[text=black,line width=0.0cm,anchor=south west] (text30) at (2.9692, 5.9002){$3$};

  \node[text=black,line width=0.0cm,anchor=south west] (text31) at (6.8498, 7.8317){$3$};

  \node[text=black,line width=0.0cm,anchor=south west] (text32) at (6.8498, 5.9002){$4$};

  \node[text=black,line width=0.0cm,anchor=south west] (text33) at (4.9095, 11.721){$3$};

  \node[text=black,line width=0.0cm,anchor=south west] (text34) at (6.8498, 9.7896){$3$};

  \node[text=black,line width=0.0cm,anchor=south west] (text35) at (8.79, 9.7896){$3$};

  \node[text=black,line width=0.0cm,anchor=south west] (text36) at (4.9095, 7.8317){$4$};

  \node[text=black,line width=0.0cm,anchor=south west] (text37) at (4.9095, 5.9002){$3$};

\end{tikzpicture}}}
    {\scalefont{1}\scalebox{0.35}{\definecolor{cd82626}{RGB}{216,38,38}

\def \globalscale {1.000000}
\begin{tikzpicture}[y=1cm, x=1cm, yscale=\globalscale,xscale=\globalscale, every node/.append style={scale=\globalscale}, inner sep=0pt, outer sep=0pt]
  \path[draw=cd82626,fill=cd82626] (0.7938, 8.3167) circle (0.0529cm);

  \path[draw=cd82626,fill=cd82626] (0.7938, 4.4362) circle (0.0529cm);

  \path[draw=cd82626,fill=cd82626] (10.4951, 4.4362) circle (0.0529cm);

  \path[draw=cd82626,fill=cd82626] (10.4951, 10.257) circle (0.0529cm);

  \path[draw=cd82626,fill=cd82626] (4.6743, 10.257) circle (0.0529cm);

  \path[draw=cd82626,fill=cd82626] (4.6743, 0.5556) circle (0.0529cm);

  \path[draw=cd82626,fill=cd82626] (8.5549, 0.5556) circle (0.0529cm);

  \path[draw=cd82626,fill=cd82626] (8.5549, 12.1973) circle (0.0529cm);

  \path[draw=cd82626,fill=cd82626] (2.734, 12.1973) circle (0.0529cm);

  \path[draw=cd82626,fill=cd82626] (2.734, 6.3765) circle (0.0529cm);

  \path[draw=cd82626,fill=cd82626] (6.6146, 6.3765) circle (0.0529cm);

  \path[draw=cd82626,fill=cd82626] (6.6146, 2.4959) circle (0.0529cm);

  \path[draw=cd82626,fill=cd82626] (12.4354, 2.4959) circle (0.0529cm);

  \path[draw=cd82626,fill=cd82626] (12.4354, 8.3167) circle (0.0529cm);

  \path[draw=cd82626,line width=0.1588cm] (0.7938, 8.3167) -- (0.7938, 4.4362);

  \path[draw=cd82626,line width=0.1588cm] (0.7938, 4.4362) -- (10.4951, 4.4362);

  \path[draw=cd82626,line width=0.1588cm] (10.4951, 4.4362) -- (10.4951, 10.257);

  \path[draw=cd82626,line width=0.1588cm] (10.4951, 10.257) -- (4.6743, 10.257);

  \path[draw=cd82626,line width=0.1588cm] (4.6743, 10.257) -- (4.6743, 0.5556);

  \path[draw=cd82626,line width=0.1588cm] (4.6743, 0.5556) -- (8.5549, 0.5556);

  \path[draw=cd82626,line width=0.1588cm] (8.5549, 0.5556) -- (8.5549, 12.1973);

  \path[draw=cd82626,line width=0.1588cm] (8.5549, 12.1973) -- (2.734, 12.1973);

  \path[draw=cd82626,line width=0.1588cm] (2.734, 12.1973) -- (2.734, 6.3765);

  \path[draw=cd82626,line width=0.1588cm] (2.734, 6.3765) -- (6.6146, 6.3765);

  \path[draw=cd82626,line width=0.1588cm] (6.6146, 6.3765) -- (6.6146, 2.4959);

  \path[draw=cd82626,line width=0.1588cm] (6.6146, 2.4959) -- (12.4354, 2.4959);

  \path[draw=cd82626,line width=0.1588cm] (12.4354, 2.4959) -- (12.4354, 8.3167);

  \path[draw=cd82626,line width=0.1588cm] (12.4354, 8.3167) -- (0.7938, 8.3167);

  \tikzstyle{every node}=[font=\fontsize{30}{30}\selectfont]
  \node[text=black,line width=0.0cm,anchor=south west] (text27) at (1.0289, 7.276){$3$};

  \node[text=black,line width=0.0cm,anchor=south west] (text28) at (1.0289, 3.4131){$5$};

  \node[text=black,line width=0.0cm,anchor=south west] (text29) at (8.79, 7.276){$3$};

  \node[text=black,line width=0.0cm,anchor=south west] (text30) at (8.79, 3.4131){$3$};

  \node[text=black,line width=0.0cm,anchor=south west] (text31) at (4.9095, 7.276){$5$};

  \node[text=black,line width=0.0cm,anchor=south west] (text32) at (6.8498, 3.4131){$3$};

  \node[text=black,line width=0.0cm,anchor=south west] (text33) at (8.79, 9.234){$3$};

  \node[text=black,line width=0.0cm,anchor=south west] (text34) at (4.9095, 9.234){$3$};

  \node[text=black,line width=0.0cm,anchor=south west] (text35) at (2.9692, 11.1654){$3$};

  \node[text=black,line width=0.0cm,anchor=south west] (text36) at (4.9095, 5.3446){$3$};

  \node[text=black,line width=0.0cm,anchor=south west] (text37) at (4.9095, 3.4131){$3$};

  \node[text=black,line width=0.0cm,anchor=south west] (text38) at (2.9692, 7.276){$3$};

\end{tikzpicture}}}
    {\scalefont{1}\scalebox{0.35}{\definecolor{cd82626}{RGB}{216,38,38}
\definecolor{c2626d8}{RGB}{38,38,216}

\def \globalscale {1.000000}
\begin{tikzpicture}[y=1cm, x=1cm, yscale=\globalscale,xscale=\globalscale, every node/.append style={scale=\globalscale}, inner sep=0pt, outer sep=0pt]
  \path[draw=cd82626,fill=cd82626] (6.6146, 8.3167) circle (0.0529cm);

  \path[draw=cd82626,fill=cd82626] (10.4951, 8.3167) circle (0.0529cm);

  \path[draw=cd82626,fill=cd82626] (10.4951, 0.5556) circle (0.0529cm);

  \path[draw=cd82626,fill=cd82626] (4.6743, 0.5556) circle (0.0529cm);

  \path[draw=cd82626,fill=cd82626] (4.6743, 10.257) circle (0.0529cm);

  \path[draw=cd82626,fill=cd82626] (12.4354, 10.257) circle (0.0529cm);

  \path[draw=cd82626,fill=cd82626] (12.4354, 6.3765) circle (0.0529cm);

  \path[draw=cd82626,fill=cd82626] (2.734, 6.3765) circle (0.0529cm);

  \path[draw=cd82626,fill=cd82626] (2.734, 2.4959) circle (0.0529cm);

  \path[draw=cd82626,fill=cd82626] (6.6146, 2.4959) circle (0.0529cm);

  \path[draw=c2626d8,fill=c2626d8] (0.7938, 4.4362) circle (0.0529cm);

  \path[draw=c2626d8,fill=c2626d8] (0.7938, 12.1973) circle (0.0529cm);

  \path[draw=c2626d8,fill=c2626d8] (8.5549, 12.1973) circle (0.0529cm);

  \path[draw=c2626d8,fill=c2626d8] (8.5549, 4.4362) circle (0.0529cm);

  \path[draw=cd82626,line width=0.1588cm] (6.6146, 8.3167) -- (10.4951, 8.3167);

  \path[draw=cd82626,line width=0.1588cm] (10.4951, 8.3167) -- (10.4951, 0.5556);

  \path[draw=cd82626,line width=0.1588cm] (10.4951, 0.5556) -- (4.6743, 0.5556);

  \path[draw=cd82626,line width=0.1588cm] (4.6743, 0.5556) -- (4.6743, 10.257);

  \path[draw=cd82626,line width=0.1588cm] (4.6743, 10.257) -- (12.4354, 10.257);

  \path[draw=cd82626,line width=0.1588cm] (12.4354, 10.257) -- (12.4354, 6.3765);

  \path[draw=cd82626,line width=0.1588cm] (12.4354, 6.3765) -- (2.734, 6.3765);

  \path[draw=cd82626,line width=0.1588cm] (2.734, 6.3765) -- (2.734, 2.4959);

  \path[draw=cd82626,line width=0.1588cm] (2.734, 2.4959) -- (6.6146, 2.4959);

  \path[draw=cd82626,line width=0.1588cm] (6.6146, 2.4959) -- (6.6146, 8.3167);

  \path[draw=c2626d8,line width=0.1588cm] (0.7938, 4.4362) -- (0.7938, 12.1973);

  \path[draw=c2626d8,line width=0.1588cm] (0.7938, 12.1973) -- (8.5549, 12.1973);

  \path[draw=c2626d8,line width=0.1588cm] (8.5549, 12.1973) -- (8.5549, 4.4362);

  \path[draw=c2626d8,line width=0.1588cm] (8.5549, 4.4362) -- (0.7938, 4.4362);

  \tikzstyle{every node}=[font=\fontsize{30}{30}\selectfont]
  \node[text=black,line width=0.0cm,anchor=south west] (text27) at (4.9095, 9.2869){$4$};

  \node[text=black,line width=0.0cm,anchor=south west] (text28) at (6.8498, 7.329){$3$};

  \node[text=black,line width=0.0cm,anchor=south west] (text29) at (8.79, 9.2869){$3$};

  \node[text=black,line width=0.0cm,anchor=south west] (text30) at (8.79, 7.329){$3$};

  \node[text=black,line width=0.0cm,anchor=south west] (text31) at (4.9095, 1.5081){$4$};

  \node[text=black,line width=0.0cm,anchor=south west] (text32) at (1.0289, 3.466){$4$};

  \node[text=black,line width=0.0cm,anchor=south west] (text33) at (1.0289, 11.2183){$3$};

  \node[text=black,line width=0.0cm,anchor=south west] (text34) at (2.9692, 5.3975){$3$};

  \node[text=black,line width=0.0cm,anchor=south west] (text35) at (4.9095, 5.3975){$4$};

  \node[text=black,line width=0.0cm,anchor=south west] (text36) at (6.8498, 5.3975){$3$};

  \node[text=black,line width=0.0cm,anchor=south west] (text37) at (2.9692, 3.466){$3$};

  \node[text=black,line width=0.0cm,anchor=south west] (text38) at (4.9095, 3.466){$3$};

\end{tikzpicture}}}
    {\scalefont{1}\scalebox{0.35}{\definecolor{cd82626}{RGB}{216,38,38}
\definecolor{c2626d8}{RGB}{38,38,216}
\definecolor{c26d826}{RGB}{38,216,38}
\definecolor{c26d8d8}{RGB}{38,216,216}

\def \globalscale {1.000000}
\begin{tikzpicture}[y=1cm, x=1cm, yscale=\globalscale,xscale=\globalscale, every node/.append style={scale=\globalscale}, inner sep=0pt, outer sep=0pt]
  \path[draw=cd82626,fill=cd82626] (12.4279, 7.2002) circle (0.0529cm);

  \path[draw=cd82626,fill=cd82626] (12.4279, 2.2128) circle (0.0529cm);

  \path[draw=cd82626,fill=cd82626] (4.1124, 2.2128) circle (0.0529cm);

  \path[draw=cd82626,fill=cd82626] (4.1124, 7.2002) circle (0.0529cm);

  \path[draw=c2626d8,fill=c2626d8] (9.1017, 10.5251) circle (0.0529cm);

  \path[draw=c2626d8,fill=c2626d8] (0.7862, 10.5251) circle (0.0529cm);

  \path[draw=c2626d8,fill=c2626d8] (0.7862, 3.8753) circle (0.0529cm);

  \path[draw=c2626d8,fill=c2626d8] (9.1017, 3.8753) circle (0.0529cm);

  \path[draw=c26d826,fill=c26d826] (10.7648, 5.5377) circle (0.0529cm);

  \path[draw=c26d826,fill=c26d826] (7.4386, 5.5377) circle (0.0529cm);

  \path[draw=c26d826,fill=c26d826] (7.4386, 8.8627) circle (0.0529cm);

  \path[draw=c26d826,fill=c26d826] (10.7648, 8.8627) circle (0.0529cm);

  \path[draw=c26d8d8,fill=c26d8d8] (5.7755, 8.8627) circle (0.0529cm);

  \path[draw=c26d8d8,fill=c26d8d8] (2.4493, 8.8627) circle (0.0529cm);

  \path[draw=c26d8d8,fill=c26d8d8] (2.4493, 0.5503) circle (0.0529cm);

  \path[draw=c26d8d8,fill=c26d8d8] (5.7755, 0.5503) circle (0.0529cm);

  \path[draw=cd82626,line width=0.1588cm] (12.4279, 7.2002) -- (12.4279, 2.2128);

  \path[draw=cd82626,line width=0.1588cm] (12.4279, 2.2128) -- (4.1124, 2.2128);

  \path[draw=cd82626,line width=0.1588cm] (4.1124, 2.2128) -- (4.1124, 7.2002);

  \path[draw=cd82626,line width=0.1588cm] (4.1124, 7.2002) -- (12.4279, 7.2002);

  \path[draw=c2626d8,line width=0.1588cm] (9.1017, 10.5251) -- (0.7862, 10.5251);

  \path[draw=c2626d8,line width=0.1588cm] (0.7862, 10.5251) -- (0.7862, 3.8753);

  \path[draw=c2626d8,line width=0.1588cm] (0.7862, 3.8753) -- (9.1017, 3.8753);

  \path[draw=c2626d8,line width=0.1588cm] (9.1017, 3.8753) -- (9.1017, 10.5251);

  \path[draw=c26d826,line width=0.1588cm] (10.7648, 5.5377) -- (7.4386, 5.5377);

  \path[draw=c26d826,line width=0.1588cm] (7.4386, 5.5377) -- (7.4386, 8.8627);

  \path[draw=c26d826,line width=0.1588cm] (7.4386, 8.8627) -- (10.7648, 8.8627);

  \path[draw=c26d826,line width=0.1588cm] (10.7648, 8.8627) -- (10.7648, 5.5377);

  \path[draw=c26d8d8,line width=0.1588cm] (5.7755, 8.8627) -- (2.4493, 8.8627);

  \path[draw=c26d8d8,line width=0.1588cm] (2.4493, 8.8627) -- (2.4493, 0.5503);

  \path[draw=c26d8d8,line width=0.1588cm] (2.4493, 0.5503) -- (5.7755, 0.5503);

  \path[draw=c26d8d8,line width=0.1588cm] (5.7755, 0.5503) -- (5.7755, 8.8627);

  \tikzstyle{every node}=[font=\fontsize{30}{30}\selectfont]
  \node[text=black,line width=0.0cm,anchor=south west] (text31) at (1.0029, 2.8575){$4$};

  \node[text=black,line width=0.0cm,anchor=south west] (text32) at (5.9921, 2.8575){$4$};

  \node[text=black,line width=0.0cm,anchor=south west] (text33) at (2.6659, 2.8575){$3$};

  \node[text=black,line width=0.0cm,anchor=south west] (text34) at (4.329, 2.8575){$3$};

  \node[text=black,line width=0.0cm,anchor=south west] (text35) at (2.6659, 7.8581){$3$};

  \node[text=black,line width=0.0cm,anchor=south west] (text36) at (4.329, 6.1912){$3$};

  \node[text=black,line width=0.0cm,anchor=south west] (text37) at (1.0029, 9.525){$4$};

  \node[text=black,line width=0.0cm,anchor=south west] (text38) at (5.9921, 6.1912){$4$};

  \node[text=black,line width=0.0cm,anchor=south west] (text39) at (7.6552, 7.8581){$3$};

  \node[text=black,line width=0.0cm,anchor=south west] (text40) at (7.6552, 6.1912){$3$};

  \node[text=black,line width=0.0cm,anchor=south west] (text41) at (9.3183, 7.8581){$3$};

  \node[text=black,line width=0.0cm,anchor=south west] (text42) at (9.3183, 6.1912){$3$};

\end{tikzpicture}}}
     
    \caption{The irreducible indecomposible spherical multiloops (unoriented, up to reflection) with at most $12$ regions whose regions all have degrees $\geq 3$ (see \cite[A078666]{oeis}).}
    \label{fig:smallest_multiloops}
\end{figure}

\subsection{From multiloops and pin sets to words in the free group}
\label{subsec:MuLoops-to-FreeGroups}

In this subsection we restrict, for simplicity of the exposition, to the case $\Surf=\Sph^2$.

Consider a multiloop $\gamma \colon \sqcup_1^s \Sph^1 \looparrowright \Surf$ together with a subset $P\subset \Regions$ of cardinal $p+1\in \N_{>0}$.
We address the question of checking whether $P$ is a pinning-set, namely whether $\gamma$ is a taut loop in the surface $\Surf_P = \Surf \setminus P$.

We choose a base point $\infty\in \Surf \setminus \Gamma$ in a region belonging to $P$, called the unbounded region. Hence the fundamental group $\pi_1(\Surf_P,\infty)$ is free of rank $p$.

Recall that a filling multiloop $\gamma\colon \sqcup_1^s \Sph^1 \looparrowright \Surf$ is determined, modulo orientation of its strands, by its associated plane map $\Gamma \subset \Surf$. 
In its dual graph $\Gamma^* \subset \Surf$, we choose to include $\infty$ as a vertex which we call the root.

Fix a spanning tree $T^*$ for the dual map $\Gamma^* \subset \Surf$, orient its edges outwards from the root vertex $\infty$, and co-orient them using the orientation of the disc.
We also order its vertices and its edges by a depth first search from $\infty$.

\paragraph{Step 1: $\pi_1(\Surf_\Regions, \infty)$.}

We first compute an explicit presentation for the free group $\pi_1(\Surf_\Regions, \infty)$ where $\Regions = \pi_0(\Surf\setminus \Gamma)$ consists of all regions, from which we may easily express the conjugacy classes of the strands $\gamma_j$. The construction follows Figure \ref{fig:dual-tree-polygon-chordiag}.

Consider simple loops $\tau_1,\dots,\tau_r$ which intersect only at their base point $\infty$, such that $\tau_j$ intersects $T^*$ only at the edge with label $r$ and in the direction prescribed by its co-orientation. 
The free group $\pi_1(\Surf_\Regions,\infty)$ is freely generated by $\tau=\{\tau_1,\dots,\tau_r\}$.
The symmetric generating set $\tau^\pm =\{\tau_j^{\pm 1}\mid 1 \le j \le r\}$ inherits the cyclic order given by the emergence of paths from $\infty$.
%
The complement $\Surf\setminus T^*$ is a polygon whose $2r$ edges are labeled by $\tau^\pm$ in accordance with its cyclic order.
It lifts in the universal cover of $\Surf_\Regions$ to a fundamental domain under the action of $\pi_1(\Surf_\Regions,\infty)$, whose generators $\tau^\pm$ identify its edges two-by-two.

\begin{example}[from a multiloop to words in $\pi_1(\Surf_\Regions)$]
\label{ex:spanning-tree-pi1(S_R)}
Figure \ref{fig:dual-tree-polygon-chordiag} represents a multiloop $\gamma$ with $2$ strands and $1+9$ regions in $\Sph^2$.
In the dual map of $\Gamma$, we choose a co-oriented labeled spanning tree $T^*$ rooted at infinity.
The sphere punctured at the vertices of $T^*$ retracts by deformation onto the dual bouquet, so $\pi_1(\Sph^2_\Regions, \infty)\cong \Free_9=\langle a,b,c,d,e,f,g,h,i\rangle$. 
Up to conjugacy, we have $\gamma_1=aeH$, and $\gamma_2=bCDFIg$ (where capital letters denote inverses). Note the cyclic order induced on the symmetrized generating set: $hgFedDacCAEBbfGHIi$.
    \begin{figure}[h]
	\centering
	\scalebox{0.8}{\definecolor{cd82626}{RGB}{216,38,38}
\definecolor{c2626d8}{RGB}{38,38,216}

\def \globalscale {1.000000}
\begin{tikzpicture}[y=1cm, x=1cm, yscale=\globalscale,xscale=\globalscale, every node/.append style={scale=\globalscale}, inner sep=0pt, outer sep=0pt]
  \path[draw=cd82626,fill=cd82626] (5.8738, 7.964) circle (0.0529cm);

  \path[draw=cd82626,fill=cd82626] (9.7543, 7.964) circle (0.0529cm);

  \path[draw=cd82626,fill=cd82626] (9.7543, 0.2117) circle (0.0529cm);

  \path[draw=cd82626,fill=cd82626] (3.9335, 0.2117) circle (0.0529cm);

  \path[draw=cd82626,fill=cd82626] (3.9335, 7.964) circle (0.0529cm);

  \path[draw=cd82626,fill=cd82626] (0.0529, 7.964) circle (0.0529cm);

  \path[draw=cd82626,fill=cd82626] (0.0529, 2.1497) circle (0.0529cm);

  \path[draw=cd82626,fill=cd82626] (5.8738, 2.1497) circle (0.0529cm);

  \path[draw=c2626d8,fill=c2626d8] (1.9932, 4.0878) circle (0.0529cm);

  \path[draw=c2626d8,fill=c2626d8] (1.9932, 6.0259) circle (0.0529cm);

  \path[draw=c2626d8,fill=c2626d8] (7.814, 6.0259) circle (0.0529cm);

  \path[draw=c2626d8,fill=c2626d8] (7.814, 2.1497) circle (0.0529cm);

  \path[draw=c2626d8,fill=c2626d8] (11.6946, 2.1497) circle (0.0529cm);

  \path[draw=c2626d8,fill=c2626d8] (11.6946, 4.0878) circle (0.0529cm);

  \path[draw=cd82626,line width=0.1588cm] (5.8738, 7.964) -- (9.7543, 7.964);

  \path[draw=cd82626,line width=0.1588cm] (9.7543, 7.964) -- (9.7543, 0.2117);

  \path[draw=cd82626,line cap=butt,line join=miter,line width=0.1323cm,miter limit=4.0,<-] (3.9335, 4.8146) -- (3.9335, 4.0878);

  \path[draw=cd82626,line width=0.1588cm] (9.7543, 0.2117) -- (3.9335, 0.2117);

  \path[draw=cd82626,line width=0.1588cm] (3.9335, 0.2117) -- (3.9335, 7.964);

  \path[draw=cd82626,line width=0.1588cm] (3.9335, 7.964) -- (0.0529, 7.964);

  \path[draw=cd82626,line width=0.1588cm] (0.0529, 7.964) -- (0.0529, 2.1497);

  \path[draw=cd82626,line width=0.1588cm] (0.0529, 2.1497) -- (5.8738, 2.1497);

  \path[draw=cd82626,line width=0.1588cm] (5.8738, 2.1497) -- (5.8738, 7.964);

  \path[draw=c2626d8,line width=0.1588cm] (1.9932, 4.0878) -- (1.9932, 6.0259);

  \path[draw=c2626d8,line width=0.1588cm] (1.9932, 6.0259) -- (7.814, 6.0259);

  \path[draw=c2626d8,line width=0.1588cm] (7.814, 6.0259) -- (7.814, 2.1497);

  \path[draw=c2626d8,line width=0.1588cm] (7.814, 2.1497) -- (11.6946, 2.1497);

  \path[draw=c2626d8,line width=0.1588cm] (11.6946, 2.1497) -- (11.6946, 4.0878);

  \path[draw=c2626d8,line width=0.1588cm] (11.6946, 4.0878) -- (1.9932, 4.0878);

  \path[fill=black,even odd rule,line cap=round,line width=0.0353cm] (8.7842, 5.0568) circle (0.081cm);

  \path[fill=black,even odd rule,line cap=round,line width=0.0353cm] (10.7244, 5.0568) circle (0.081cm);

  \path[fill=black,even odd rule,line cap=round,line width=0.0353cm] (10.7244, 3.1188) circle (0.081cm);

  \path[fill=black,even odd rule,line cap=round,line width=0.0353cm] (10.7244, 3.1188) circle (0.081cm);

  \path[fill=black,even odd rule,line cap=round,line width=0.0353cm] (8.7842, 3.1188) circle (0.081cm);

  \path[fill=black,even odd rule,line cap=round,line width=0.0353cm] (6.8465, 3.1188) circle (0.081cm);

  \path[fill=black,even odd rule,line cap=round,line width=0.0353cm] (6.8439, 5.0568) circle (0.081cm);

  \path[fill=black,even odd rule,line cap=round,line width=0.0353cm] (4.9036, 5.0568) circle (0.081cm);

  \path[fill=black,even odd rule,line cap=round,line width=0.0353cm] (4.9036, 3.1188) circle (0.081cm);

  \path[fill=black,even odd rule,line cap=round,line width=0.0353cm] (2.9633, 5.0568) circle (0.081cm);

  \path[fill=black,even odd rule,line cap=round,line width=0.0353cm] (1.0231, 5.0568) circle (0.081cm);

  \path[draw=black,line cap=butt,line join=miter,line width=0.0529cm] (1.0231, 5.0568) -- (2.9633, 5.0568);

  \path[draw=black,line cap=butt,line join=miter,line width=0.0529cm] (2.9633, 5.0568) -- (4.9036, 5.0568);

  \path[draw=black,line cap=butt,line join=miter,line width=0.0529cm] (4.9036, 5.0568) -- (6.8439, 5.0568);

  \path[draw=black,line cap=butt,line join=miter,line width=0.0529cm] (6.8439, 5.0568) -- (8.7842, 5.0568);

  \path[draw=black,line cap=butt,line join=miter,line width=0.0529cm] (8.7842, 3.1188) -- (10.7244, 3.1188);

  \path[draw=black,line cap=butt,line join=miter,line width=0.0529cm] (6.8465, 3.1188) -- (8.7868, 3.1188);

  \path[draw=black,line cap=butt,line join=miter,line width=0.0529cm] (4.9036, 5.0568) -- (4.9036, 3.1188);

  \path[draw=black,line cap=butt,line join=miter,line width=0.0529cm] (6.8439, 5.0568) -- (6.8439, 3.1188);

  \path[draw=black,line cap=butt,line join=miter,line width=0.0529cm] (10.7244, 5.0568) -- (10.7244, 3.1188);

  \path[draw=black,line cap=butt,line join=miter,line width=0.0529cm,<-] (6.3594, 5.2991) -- (6.3594, 4.8146);

  \path[draw=black,line cap=butt,line join=miter,line width=0.0529cm,<-] (8.2985, 5.2991) -- (8.2985, 4.8146);

  \path[draw=black,line cap=butt,line join=miter,line width=0.0529cm,<-] (10.2388, 3.361) -- (10.2388, 2.8765);

  \path[draw=black,line cap=butt,line join=miter,line width=0.0529cm,<-] (8.2985, 3.361) -- (8.2985, 2.8765);

  \path[draw=black,line cap=butt,line join=miter,line width=0.0529cm,<-] (4.418, 5.2991) -- (4.418, 4.8146);

  \path[draw=black,line cap=butt,line join=miter,line width=0.0529cm,<-] (2.4777, 5.2991) -- (2.4777, 4.8146);

  \path[draw=black,line cap=butt,line join=miter,line width=0.0529cm,<-] (4.6614, 4.5723) -- (5.1459, 4.5723);

  \path[draw=black,line cap=butt,line join=miter,line width=0.0529cm,<-] (6.6016, 4.5723) -- (7.0861, 4.5723);

  \path[draw=black,line cap=butt,line join=miter,line width=0.0529cm,<-] (10.4822, 4.5723) -- (10.9667, 4.5723);

  \path[draw=c2626d8,line cap=butt,line join=miter,line width=0.1323cm,miter limit=4.0,<-] (7.814, 4.8146) -- (7.814, 4.0878);

  \tikzstyle{every node}=[font=\fontsize{12}{12}\selectfont]
  \node[anchor=south west] (text53) at (2.1962, 5.2536){$c$};

  \node[anchor=south west] (text54) at (4.0968, 5.2595){$a$};

  \node[anchor=south west] (text55) at (4.3776, 4.2745){$d$};

  \node[anchor=south west] (text56) at (6.0744, 5.2449){$e$};

  \node[anchor=south west] (text57) at (6.3178, 4.2275){$f$};

  \node[anchor=south west] (text58) at (7.9876, 5.2329){$b$};

  \node[anchor=south west] (text59) at (7.9876, 3.2664){$g$};

  \node[anchor=south west] (text60) at (9.9022, 3.3109){$h$};

  \node[anchor=south west] (text61) at (10.3047, 4.2796){$i$};

  \node[text=cd82626,anchor=south west] (text1) at (3.3334, 4.3364){$\gamma_1$};

  \node[text=c2626d8,anchor=south west] (text2) at (8.0919, 4.3135){$\gamma_2$};

\end{tikzpicture}}
    {\input{images/tikz/dualbouquet_v2.tex}}
	\caption{Accompanying Example \ref{ex:spanning-tree-pi1(S_R)}: from a multiloop to words in the free group.}
	\label{fig:dual-tree-polygon-chordiag}
\end{figure}
\end{example}

Note that each strand $\gamma_j$ of $\gamma$ corresponds to a conjugacy class in $\pi_1(\Surf_\Regions,\infty)$, hence to a unique reduced cyclic word over $\tau^\pm$ which may be computed by recording the sequence of edges of $T^*$ that it crosses, minding the co-orientation to determine the exponent of each generator.

\paragraph{Step 2: $\pi_1(\Surf_P,\infty)$.}

Now for a subset of regions $P\subset \Regions$, we deduce a presentation of the free group $\pi_1(\Surf_P,\infty)$, and the new reduced cyclic words associated to the conjugacy classes of the strands $\gamma_j$. The construction follows Figure \ref{fig:filling-punctures}.

By the van Kampen Theorem, the group $\pi_1(\Surf_P,\infty)$ is the quotient of $\pi_1(\Surf_{\Regions}, \infty)$ by the normal subgroup generated by small loops only surrounding the punctures to be filled, associated to the set of regions $Q=\Regions \setminus P$.
%

A region $R_j \in \Regions$ corresponds to a non-rooted vertex of $T^*$, which has one incoming edge and a certain number of outgoing edges.
Let $\rho_j$ be the word over $\tau^\pm$ obtained by reading the cutting sequence of $T^*$ of a loop winding once around $R_j$, starting at the incoming edge.
Thus $\rho_j$ is a product of $\deg(R_j)$ generators in $\tau^\pm$ which are all distinct even up to inversion.
We deduce the presentation:
\begin{equation*}
\pi_1(\Surf_P,\infty)=\pi_1(\Surf_{\Regions}, \infty) \mod{\vartriangleleft \{ \rho_j \mid R_j \in Q\} \vartriangleright}
\end{equation*}
in which every relator $\rho(R_j)$ yields a rewriting rule of the form $\tau_{j_1} \mapsto \tau^{\mp}_{j_d}\cdots \tau^{\mp}_{j_2}$ which we obtain by solving the equation $\rho(R_j)=\Id$.

Moreover, for an arbitrary word $\alpha$ over $\tau^\pm$, we may apply these rewriting rules in the order prescribed by a depth-first search from $\infty$ to obtain a canonical representative of $\alpha$ in $\pi_1(\Surf_P,\infty)$ which uses none of the generators associated to regions in $Q$.
This holds in particular for the words associated to the strands of $\gamma$.

Observe that the remaining generators $\tau^\pm_P=\{\tau_j^\pm \mid R_j\in P\}$ of $\pi_1(\Surf_P,\infty)$ inherit a cyclic order from $\tau^\pm$, which is compatible with their action on a lift of $\Surf_P\setminus T^*$ to a polygonal fundamental domain in the universal cover of $\Surf_P$. See Figure \ref{fig:filling-punctures} and \ref{fig:lifts_intersecting_domain}.

\begin{example}[filling punctures and adding relations]
\label{ex:filling-punctures}
Following Example \ref{ex:spanning-tree-pi1(S_R)} in Figure \ref{fig:filling-punctures}, we choose to pin $P=\{r,p_1,p_2\}$ and fill the other punctures $Q = \Regions\setminus P$.
The new root $r$ is the closest vertex in $P$ to the old root $\infty$.
In $T^*$, we keep one edge outgoing from each vertex of $P$, and the other edge (in light green) yields gets collapsed and yields a relation a rewriting rule in the fundamental group.
The cyclic order on the remaining generators which correspond to the surviving edges is inherited from the original one: $hgFedDacCAEBbfGHIi\mapsto aABb$.
Thus \[\pi_1(\Sph^2_P,\infty)\cong \langle a,b,c,d,e,f,g,h,i \mid i=1, h=g, g=F, f=bE, e=da, c=1, d=1 \rangle\]
and applying the rewriting rules in DFS order from the root $r$, the loops $\gamma_1$ and $\gamma_2$ become:
\begin{align*}
    &\gamma_1=aeH\mapsto_2 aeG\mapsto_3 aef\mapsto_4 aebE\mapsto_5 adabAD\mapsto_7 aabA \sim ab \\
    &\gamma_2=bCDFIg\mapsto_{1,3}bCDFF\mapsto_4 bCDeBeB\mapsto_5 bCDdaBdaB \mapsto_{6,7}baBaB\sim aBa
\end{align*}
\vspace{-0.5cm}
\begin{figure}[h]
	\centering
	\scalebox{0.9}{\definecolor{c018000}{RGB}{1,128,0}
\definecolor{lime}{RGB}{0,255,0}
\definecolor{c00fa00}{RGB}{0,250,0}

\def \globalscale {1.000000}
\begin{tikzpicture}[y=1cm, x=1cm, yscale=\globalscale,xscale=\globalscale, every node/.append style={scale=\globalscale}, inner sep=0pt, outer sep=0pt]
  \path[draw=black,line cap=butt,line join=miter,line width=0.0529cm] (1.0231, 5.0568) -- (2.9633, 5.0568);

  \path[draw=black,line cap=butt,line join=miter,line width=0.0529cm] (2.9633, 5.0568) -- (4.9036, 5.0568);

  \path[draw=black,line cap=butt,line join=miter,line width=0.0529cm] (4.9036, 5.0568) -- (6.8439, 5.0568);

  \path[draw=black,line cap=butt,line join=miter,line width=0.0529cm] (6.8439, 5.0568) -- (8.7842, 5.0568);

  \path[draw=black,line cap=butt,line join=miter,line width=0.0529cm] (8.7842, 3.1188) -- (10.7244, 3.1188);

  \path[draw=black,line cap=butt,line join=miter,line width=0.0529cm] (6.8465, 3.1188) -- (8.7868, 3.1188);

  \path[draw=black,line cap=butt,line join=miter,line width=0.0529cm] (4.9036, 5.0568) -- (4.9036, 3.1188);

  \path[draw=black,line cap=butt,line join=miter,line width=0.0529cm] (6.8439, 5.0568) -- (6.8439, 3.1188);

  \path[draw=black,line cap=butt,line join=miter,line width=0.0529cm] (10.7244, 5.0568) -- (10.7244, 3.1188);

  \path[draw=black,line cap=butt,line join=miter,line width=0.0529cm,<-] (6.3594, 5.2991) -- (6.3594, 4.8146);

  \path[draw=black,line cap=butt,line join=miter,line width=0.0529cm,<-] (8.2985, 5.2991) -- (8.2985, 4.8146);

  \path[draw=black,line cap=butt,line join=miter,line width=0.0529cm,<-] (10.2388, 3.361) -- (10.2388, 2.8765);

  \path[draw=black,line cap=butt,line join=miter,line width=0.0529cm,<-] (8.2985, 3.361) -- (8.2985, 2.8765);

  \path[draw=black,line cap=butt,line join=miter,line width=0.0529cm,<-] (4.418, 5.2991) -- (4.418, 4.8146);

  \path[draw=black,line cap=butt,line join=miter,line width=0.0529cm,<-] (2.4777, 5.2991) -- (2.4777, 4.8146);

  \path[draw=black,line cap=butt,line join=miter,line width=0.0529cm,<-] (4.6614, 4.5723) -- (5.1459, 4.5723);

  \path[draw=black,line cap=butt,line join=miter,line width=0.0529cm,<-] (6.6016, 4.5723) -- (7.0861, 4.5723);

  \path[draw=black,line cap=butt,line join=miter,line width=0.0529cm,<-] (10.4822, 4.5723) -- (10.9667, 4.5723);

  \tikzstyle{every node}=[font=\fontsize{12}{12}\selectfont]
  \node[anchor=south west] (text53) at (2.1962, 5.2536){$c$};

  \node[anchor=south west] (text54) at (4.0968, 5.2595){$a$};

  \node[anchor=south west] (text55) at (4.3776, 4.2745){$d$};

  \node[anchor=south west] (text56) at (6.0744, 5.2449){$e$};

  \node[anchor=south west] (text57) at (6.3178, 4.2275){$f$};

  \node[anchor=south west] (text58) at (7.9876, 5.2329){$b$};

  \node[anchor=south west] (text59) at (7.9876, 3.2664){$g$};

  \node[anchor=south west] (text60) at (9.9022, 3.3109){$h$};

  \node[anchor=south west] (text61) at (10.3047, 4.2796){$i$};

  \path[fill=black,even odd rule,line cap=round,line width=0.0353cm] (10.7244, 5.0568) circle (0.081cm);

  \path[fill=black,even odd rule,line cap=round,line width=0.0353cm] (10.7244, 3.1188) circle (0.081cm);

  \path[fill=black,even odd rule,line cap=round,line width=0.0353cm] (8.7842, 3.1188) circle (0.081cm);

  \path[fill=black,even odd rule,line cap=round,line width=0.0353cm] (6.8465, 3.1188) circle (0.081cm);

  \path[fill=black,even odd rule,line cap=round,line width=0.0353cm] (6.8439, 5.0568) circle (0.081cm);

  \path[fill=black,even odd rule,line cap=round,line width=0.0353cm] (4.9036, 5.0568) circle (0.081cm);

  \path[fill=black,even odd rule,line cap=round,line width=0.0353cm] (4.9036, 3.1188) circle (0.081cm);

  \path[fill=c018000,even odd rule,line cap=round,line join=miter,line width=0.0353cm,miter limit=4.0] (2.9633, 5.0568) circle (0.1995cm);

  \path[fill=black,even odd rule,line cap=round,line width=0.0353cm] (1.0231, 5.0568) circle (0.081cm);

  \path[fill=c018000,even odd rule,line cap=round,line width=0.0353cm] (8.7842, 5.0568) circle (0.1995cm);

  \path[fill=c018000,even odd rule,line cap=round,line join=miter,line width=0.0353cm,miter limit=4.0] (10.7244, 3.1166) circle (0.1995cm);

  \node[text=c018000,line cap=round,line join=miter,line width=0.0353cm,miter limit=4.0,anchor=south west] (text26) at (3.0678, 5.2963){$p_2$};

  \node[text=c018000,line cap=round,line join=miter,line width=0.0353cm,miter limit=4.0,anchor=south west] (text27) at (8.8927, 5.3007){$p_1$};

  \node[text=c018000,line cap=round,line join=miter,line width=0.0353cm,miter limit=4.0,anchor=south west] (text28) at (10.8921, 2.7704){$r$};

  \path[fill=lime,opacity=0.5,even odd rule,line cap=round,line width=0.0353cm,rotate around={45.0:(0.0, 7.911)}] (2.8211, 1.0535) circle (0.3445cm);

  \path (7.0505, 5.3352) -- (7.0635, 5.3221);

  \path[draw=lime,opacity=0.5,line cap=butt,line join=miter,line width=0.1852cm] (6.8523, 5.4013) -- (6.8523, 5.6776);

  \path[draw=black,even odd rule,line cap=round,line width=0.0353cm,rotate around={45.0:(0.0, 7.911)}] (2.8211, 1.0535) circle (0.3445cm);

  \path[fill=c00fa00,opacity=0.5,even odd rule,line cap=round,line width=0.0353cm] (6.0434, 6.1936) rectangle (7.9484, 5.6776);

  \node[text=black,anchor=south west] (text1) at (6.0803, 5.6937){4) $f=bE$};

  \path[fill=lime,opacity=0.5,even odd rule,line cap=round,line width=0.0353cm,rotate around={45.0:(0.0, 7.911)}] (1.4492, 2.4254) circle (0.3445cm);

  \path[draw=lime,opacity=0.5,line cap=butt,line join=miter,line width=0.1852cm] (4.9121, 5.4013) -- (4.9121, 5.6776);

  \path[draw=black,even odd rule,line cap=round,line width=0.0353cm,rotate around={45.0:(0.0, 7.911)}] (1.4492, 2.4254) circle (0.3445cm);

  \path[fill=c00fa00,opacity=0.5,even odd rule,line cap=round,line width=0.0353cm] (4.1032, 6.1936) rectangle (5.9222, 5.6776);

  \node[text=black,anchor=south west] (text4) at (4.14, 5.6937){5) $e=da$};

  \path[fill=lime,opacity=0.5,even odd rule,line cap=round,line width=0.0353cm,rotate around={45.0:(0.0, 7.911)}] (-1.2948, 5.1694) circle (0.3445cm);

  \path[draw=lime,opacity=0.5,line cap=butt,line join=miter,line width=0.1852cm] (1.0315, 5.4013) -- (1.0315, 5.6776);

  \path[draw=black,even odd rule,line cap=round,line width=0.0353cm,rotate around={45.0:(0.0, 7.911)}] (-1.2948, 5.1694) circle (0.3445cm);

  \path[fill=c00fa00,opacity=0.5,even odd rule,line cap=round,line width=0.0353cm] (0.2226, 6.1936) rectangle (1.8521, 5.6776);

  \node[text=black,anchor=south west] (text6) at (0.2595, 5.6937){6) $c=1$};

  \path[fill=lime,opacity=0.5,even odd rule,line cap=round,line width=0.0353cm,rotate around={45.0:(0.0, 7.911)}] (5.5651, -1.6905) circle (0.3445cm);

  \path[draw=lime,opacity=0.5,line cap=butt,line join=miter,line width=0.1852cm] (10.7329, 5.4013) -- (10.7329, 5.6776);

  \path[draw=black,even odd rule,line cap=round,line width=0.0353cm,rotate around={45.0:(0.0, 7.911)}] (5.5651, -1.6905) circle (0.3445cm);

  \path[fill=c00fa00,opacity=0.5,even odd rule,line cap=round,line width=0.0353cm] (9.924, 6.1936) rectangle (11.5418, 5.6776);

  \node[text=black,anchor=south west] (text8) at (9.9609, 5.6937){1) $i=1$};

  \path[fill=lime,opacity=0.5,even odd rule,line cap=round,line width=0.0353cm,rotate around={45.0:(0.0, 7.911)}] (2.8246, -1.6908) circle (0.3445cm);

  \path[draw=lime,opacity=0.5,line cap=butt,line join=miter,line width=0.1852cm] (8.7953, 2.4982) -- (8.7953, 2.7743);

  \path[draw=black,even odd rule,line cap=round,line width=0.0353cm,rotate around={45.0:(0.0, 7.911)}] (2.8246, -1.6908) circle (0.3445cm);

  \path[fill=c00fa00,opacity=0.5,even odd rule,line cap=round,line width=0.0353cm] (7.9864, 2.4982) rectangle (9.6928, 1.9821);

  \node[text=black,anchor=south west] (text10) at (8.0232, 1.9982){2) $h=g$};

  \path[fill=lime,opacity=0.5,even odd rule,line cap=round,line width=0.0353cm,rotate around={45.0:(0.0, 7.911)}] (1.4507, -0.317) circle (0.3445cm);

  \path[draw=lime,opacity=0.5,line cap=butt,line join=miter,line width=0.1852cm] (6.8523, 2.4982) -- (6.8523, 2.7743);

  \path[draw=black,even odd rule,line cap=round,line width=0.0353cm,rotate around={45.0:(0.0, 7.911)}] (1.4507, -0.317) circle (0.3445cm);

  \path[fill=c00fa00,opacity=0.5,even odd rule,line cap=round,line width=0.0353cm] (6.0434, 2.4982) rectangle (7.7762, 1.9821);

  \node[text=black,anchor=south west] (text12) at (6.0803, 1.9982){3) $g=F$};

  \path[fill=lime,opacity=0.5,even odd rule,line cap=round,line width=0.0353cm,rotate around={45.0:(0.0, 7.911)}] (0.0787, 1.055) circle (0.3445cm);

  \path[draw=lime,opacity=0.5,line cap=butt,line join=miter,line width=0.1852cm] (4.9121, 2.4982) -- (4.9121, 2.7743);

  \path[draw=black,even odd rule,line cap=round,line width=0.0353cm,rotate around={45.0:(0.0, 7.911)}] (0.0787, 1.055) circle (0.3445cm);

  \path[fill=c00fa00,opacity=0.5,even odd rule,line cap=round,line width=0.0353cm] (4.1032, 2.4982) rectangle (5.7819, 1.9821);

  \node[text=black,anchor=south west] (text14) at (4.14, 1.9982){7) $d=1$};

\end{tikzpicture}}
    \caption*{}
\end{figure}
\vspace{-1cm}
\begin{figure}[h]
    \centering
    \scalebox{0.9}{\definecolor{cd82626}{RGB}{216,38,38}
\definecolor{c2626d8}{RGB}{38,38,216}
\definecolor{lime}{RGB}{0,255,0}
\definecolor{c018000}{RGB}{1,128,0}

\def \globalscale {1.000000}
\begin{tikzpicture}[y=1cm, x=1cm, yscale=\globalscale,xscale=\globalscale, every node/.append style={scale=\globalscale}, inner sep=0pt, outer sep=0pt]
  \path[draw=cd82626,fill=cd82626] (5.8738, 7.964) circle (0.0529cm);

  \path[draw=cd82626,fill=cd82626] (9.7543, 7.964) circle (0.0529cm);

  \path[draw=cd82626,fill=cd82626] (9.7543, 0.2117) circle (0.0529cm);

  \path[draw=cd82626,fill=cd82626] (3.9335, 0.2117) circle (0.0529cm);

  \path[draw=cd82626,fill=cd82626] (3.9335, 7.964) circle (0.0529cm);

  \path[draw=cd82626,fill=cd82626] (0.0529, 7.964) circle (0.0529cm);

  \path[draw=cd82626,fill=cd82626] (0.0529, 2.1497) circle (0.0529cm);

  \path[draw=cd82626,fill=cd82626] (5.8738, 2.1497) circle (0.0529cm);

  \path[draw=c2626d8,fill=c2626d8] (1.9932, 4.0878) circle (0.0529cm);

  \path[draw=c2626d8,fill=c2626d8] (1.9932, 6.0259) circle (0.0529cm);

  \path[draw=c2626d8,fill=c2626d8] (7.814, 6.0259) circle (0.0529cm);

  \path[draw=c2626d8,fill=c2626d8] (7.814, 2.1497) circle (0.0529cm);

  \path[draw=c2626d8,fill=c2626d8] (11.6946, 2.1497) circle (0.0529cm);

  \path[draw=c2626d8,fill=c2626d8] (11.6946, 4.0878) circle (0.0529cm);

  \path[draw=cd82626,line width=0.1588cm] (5.8738, 7.964) -- (9.7543, 7.964);

  \path[draw=cd82626,line width=0.1588cm] (9.7543, 7.964) -- (9.7543, 0.2117);

  \path[draw=cd82626,line cap=butt,line join=miter,line width=0.1323cm,miter limit=4.0,<-] (3.9335, 4.8146) -- (3.9335, 4.0878);

  \path[draw=cd82626,line width=0.1588cm] (9.7543, 0.2117) -- (3.9335, 0.2117);

  \path[draw=cd82626,line width=0.1588cm] (3.9335, 0.2117) -- (3.9335, 7.964);

  \path[draw=cd82626,line width=0.1588cm] (3.9335, 7.964) -- (0.0529, 7.964);

  \path[draw=cd82626,line width=0.1588cm] (0.0529, 7.964) -- (0.0529, 2.1497);

  \path[draw=cd82626,line width=0.1588cm] (0.0529, 2.1497) -- (5.8738, 2.1497);

  \path[draw=cd82626,line width=0.1588cm] (5.8738, 2.1497) -- (5.8738, 7.964);

  \path[draw=c2626d8,line width=0.1588cm] (1.9932, 4.0878) -- (1.9932, 6.0259);

  \path[draw=c2626d8,line width=0.1588cm] (1.9932, 6.0259) -- (7.814, 6.0259);

  \path[draw=c2626d8,line width=0.1588cm] (7.814, 6.0259) -- (7.814, 2.1497);

  \path[draw=c2626d8,line width=0.1588cm] (7.814, 2.1497) -- (11.6946, 2.1497);

  \path[draw=c2626d8,line width=0.1588cm] (11.6946, 2.1497) -- (11.6946, 4.0878);

  \path[draw=lime,line cap=butt,line join=miter,line width=0.1588cm,<-] (1.0231, 5.0568) -- (2.9633, 5.0568);

  \path[draw=lime,line cap=butt,line join=miter,line width=0.1588cm,<-] (4.9036, 5.0568) -- (6.8439, 5.0568);

  \path[draw=c2626d8,line width=0.1588cm] (11.6946, 4.0878) -- (1.9932, 4.0878);

  \path[draw=lime,line cap=butt,line join=miter,line width=0.1588cm,<-] (8.7842, 3.1188) -- (10.7244, 3.1188);

  \path[draw=lime,line cap=butt,line join=miter,line width=0.1588cm,<-] (6.8439, 5.058) -- (6.8439, 3.1177);

  \path[draw=lime,line cap=butt,line join=miter,line width=0.1588cm,<-] (6.8439, 3.1188) -- (8.7842, 3.1188);

  \path[draw=lime,line cap=butt,line join=miter,line width=0.1588cm,<-] (10.7244, 5.0568) -- (10.7244, 3.1166);

  \path[draw=lime,line cap=butt,line join=miter,line width=0.1588cm,<-] (4.9036, 3.1177) -- (4.9036, 5.058);

  \path[draw=black,line cap=butt,line join=miter,line width=0.0529cm] (1.0231, 5.0568) -- (2.9633, 5.0568);

  \path[draw=black,line cap=butt,line join=miter,line width=0.0529cm] (2.9633, 5.0568) -- (4.9036, 5.0568);

  \path[draw=black,line cap=butt,line join=miter,line width=0.0529cm] (4.9036, 5.0568) -- (6.8439, 5.0568);

  \path[draw=black,line cap=butt,line join=miter,line width=0.0529cm] (6.8439, 5.0568) -- (8.7842, 5.0568);

  \path[draw=black,line cap=butt,line join=miter,line width=0.0529cm] (8.7842, 3.1188) -- (10.7244, 3.1188);

  \path[draw=black,line cap=butt,line join=miter,line width=0.0529cm] (6.8465, 3.1188) -- (8.7868, 3.1188);

  \path[draw=black,line cap=butt,line join=miter,line width=0.0529cm] (4.9036, 5.0568) -- (4.9036, 3.1188);

  \path[draw=black,line cap=butt,line join=miter,line width=0.0529cm] (6.8439, 5.0568) -- (6.8439, 3.1188);

  \path[draw=black,line cap=butt,line join=miter,line width=0.0529cm] (10.7244, 5.0568) -- (10.7244, 3.1188);

  \path[draw=black,line cap=butt,line join=miter,line width=0.0529cm,<-] (6.3594, 5.2991) -- (6.3594, 4.8146);

  \path[draw=black,line cap=butt,line join=miter,line width=0.0529cm,<-] (8.2985, 5.2991) -- (8.2985, 4.8146);

  \path[draw=black,line cap=butt,line join=miter,line width=0.0529cm,<-] (10.2388, 3.361) -- (10.2388, 2.8765);

  \path[draw=black,line cap=butt,line join=miter,line width=0.0529cm,<-] (8.2985, 3.361) -- (8.2985, 2.8765);

  \path[draw=black,line cap=butt,line join=miter,line width=0.0529cm,<-] (4.418, 5.2991) -- (4.418, 4.8146);

  \path[draw=black,line cap=butt,line join=miter,line width=0.0529cm,<-] (2.4777, 5.2991) -- (2.4777, 4.8146);

  \path[draw=black,line cap=butt,line join=miter,line width=0.0529cm,<-] (4.6614, 4.5723) -- (5.1459, 4.5723);

  \path[draw=black,line cap=butt,line join=miter,line width=0.0529cm,<-] (6.6016, 4.5723) -- (7.0861, 4.5723);

  \path[draw=black,line cap=butt,line join=miter,line width=0.0529cm,<-] (10.4822, 4.5723) -- (10.9667, 4.5723);

  \path[draw=c2626d8,line cap=butt,line join=miter,line width=0.1323cm,miter limit=4.0,<-] (7.814, 4.8146) -- (7.814, 4.0878);

  \tikzstyle{every node}=[font=\fontsize{12}{12}\selectfont]
  \node[anchor=south west] (text53) at (2.1962, 5.2536){$c$};

  \node[anchor=south west] (text54) at (4.0968, 5.2595){$a$};

  \node[anchor=south west] (text55) at (4.3776, 4.2745){$d$};

  \node[anchor=south west] (text56) at (6.0744, 5.2449){$e$};

  \node[anchor=south west] (text57) at (6.3178, 4.2275){$f$};

  \node[anchor=south west] (text58) at (7.9876, 5.2329){$b$};

  \node[anchor=south west] (text59) at (7.9876, 3.2664){$g$};

  \node[anchor=south west] (text60) at (9.9022, 3.3109){$h$};

  \node[anchor=south west] (text61) at (10.3047, 4.2796){$i$};

  \node[text=cd82626,anchor=south west] (text1) at (3.3169, 4.3232){$\gamma_1$};

  \node[text=c2626d8,anchor=south west] (text2) at (8.0919, 4.3135){$\gamma_2$};

  \path[fill=black,even odd rule,line cap=round,line width=0.0353cm] (10.7244, 5.0568) circle (0.081cm);

  \path[fill=black,even odd rule,line cap=round,line width=0.0353cm] (10.7244, 3.1188) circle (0.081cm);

  \path[fill=black,even odd rule,line cap=round,line width=0.0353cm] (8.7842, 3.1188) circle (0.081cm);

  \path[fill=black,even odd rule,line cap=round,line width=0.0353cm] (6.8465, 3.1188) circle (0.081cm);

  \path[fill=black,even odd rule,line cap=round,line width=0.0353cm] (6.8439, 5.0568) circle (0.081cm);

  \path[fill=black,even odd rule,line cap=round,line width=0.0353cm] (4.9036, 5.0568) circle (0.081cm);

  \path[fill=black,even odd rule,line cap=round,line width=0.0353cm] (4.9036, 3.1188) circle (0.081cm);

  \path[fill=c018000,even odd rule,line cap=round,line join=miter,line width=0.0353cm,miter limit=4.0] (2.9633, 5.0568) circle (0.1995cm);

  \path[fill=black,even odd rule,line cap=round,line width=0.0353cm] (1.0231, 5.0568) circle (0.081cm);

  \path[fill=c018000,even odd rule,line cap=round,line width=0.0353cm] (8.7842, 5.0568) circle (0.1995cm);

  \path[fill=c018000,even odd rule,line cap=round,line join=miter,line width=0.0353cm,miter limit=4.0] (10.7244, 3.1166) circle (0.1995cm);

  \node[text=c018000,line cap=round,line join=miter,line width=0.0353cm,miter limit=4.0,anchor=south west] (text26) at (3.0678, 5.2963){$p_2$};

  \node[text=c018000,line cap=round,line join=miter,line width=0.0353cm,miter limit=4.0,anchor=south west] (text27) at (8.8927, 5.3007){$p_1$};

  \node[text=c018000,line cap=round,line join=miter,line width=0.0353cm,miter limit=4.0,anchor=south west] (text28) at (10.8921, 2.7704){$r$};

\end{tikzpicture}}
	\caption{Accompanying Example \ref{ex:filling-punctures}: filling punctures yields relations and rewriting rules.}
	\label{fig:filling-punctures}
\end{figure}

\end{example}

\paragraph{Conclusion.}

Let us summarize and emphasize some details of the whole procedure.

In short, given a multiloop $\gamma\colon \sqcup_1^s \Sph^1 \looparrowright \Sph^2$ and a subset of $1+p$ regions $P\subset \Regions$, we compute in time $O(\Card(\Regions))$ a presentation of the cyclically ordered group $\pi_1(\Surf_P)$ and a set of words representing the conjugacy classes of the strands $\gamma_i$.

The cyclically ordered presentation consists of the generating set $\tau_P=\{\tau_j\}_1^p$ for $ \pi_1(\Surf_P)$ labeled by the bounded regions in $P$, which in our case is free from relations, together with a cyclic order on the symmetric set $(\tau_P)^\pm= \{\tau_j^{\pm 1}\}_1^p$.

Moreover for $P_1\subset P_2 \subset \Regions$, the natural quotient map $\pi_1(\Surf_{P_2})\to \pi_1(\Surf_{P_1})$ is order preserving, and our
presentations for these cyclically ordered groups (which we described by the relations $\rho_j$ corresponding to elements in $Q = P_2\setminus P_1$) are functorial in the sense that the cyclic order on $(\tau_{P_1})^\pm$ is inherited by that of $(\tau_{P_2})^\pm$.

\begin{remark}[winding numbers]
    For a loop $\alpha\colon \Sph^1 \looparrowright \Surf_P$ with $\infty\in P$ as above, every puncture $o\in P$ has a winding number (defined as the linking number between the $1$-cycle $\alpha$ and the relative $0$-chain $[o]-[\infty]$ in the oriented surface $\Surf_P$).
    
    This yields a winding number function $w_\infty\colon \pi_1(\Surf_P)\to \Z^P$ which corresponds to the abelianization $\Free_p \to \Z^p$.
    Given our presentation of $\pi_1(\Surf_P)$ as above, this winding number function may be computed in time $\Card(P)^2$ as the matrix of winding numbers between the generators $\tau_j$ and the punctures $P_i$.
\end{remark}

\begin{remark}[higher genus]
    This algorithm adapts to filling multiloops $\gamma \colon \sqcup_1^s \Sph^1 \subset \Surf_P$ in punctured surfaces of higher genus.
    
    For this one may choose a base point $\infty \in \Surf$ which is close to a puncture, and perform cuts along disjoint simple arcs between the punctures to represent $\Surf$ as the identification of the sides of a $2n$-gon whose vertices are removed. 
    This leads to a presentation for $\pi_1(\Surf\setminus \{\infty\})=\Free_n$ by a free set of generators, together with a cyclic order on the corresponding symmetric generating set.
    The complexity is still $O(\Card(\Regions))$. 
\end{remark}

\subsection{Computing self-intersection of homotopy classes}
\label{subsec:si-algorithm}

In this subsection we explain an adaptation of the algorithm proposed in \cite{Cohen-Lustig_algo-geometric-intersection-number_1987}, extending \cite{Birman-Series_algo-simple-closed-curves_1984}, to compute the self-intersection number of a multicurve $\gamma$ in $\Surf_P$. 
%
%
For clarity and completeness, we include a description of this algorithm in relation with our previous presentation so as to explain some details not covered in those references.

The algorithm takes as input a free group $\Free_p$ of rank $p$ over a set $\tau = \{\tau_j\mid 1\le j \le p\}$ with a cyclic order on the symmetric set $\tau^\pm$, and a finite set $\gamma$ of reduced words in these generators.
It returns the self-intersection number: 
\begin{equation*} \textstyle
    \si_P(\gamma) = \sum_k \si_P(\gamma_k) + \sum_{i \ne j} \ti_P(\gamma_i,\gamma_j).
\end{equation*}

\begin{remark}
\label{rem:non_primitive_si}
For curves $\alpha,\beta$ in $\Surf_P$ and $a,b \in \N$ we have $\ti_P(\alpha^a, \beta^b)= ab\cdot \ti_P(\alpha,\beta)$.
Note that the self-intersection number of a curve satisfies $\si_P(\alpha)= \tfrac{1}{2}\ti_P(\alpha,\alpha)$. 

If $\alpha$ is primitive then for all $n\in \N^*$ we have $\si_P(\alpha^n)=n^2\si_P(\alpha)+(n-1)$. Conversely, if there exists $n\in \N_{>1}$ such that $\si_P(\alpha^n)=n^2\si_P(\alpha)+(n-1)$ then $\alpha$ is primitive (because otherwise $\alpha = \beta^m$ for some primitive $\beta$ and $m\in \N_{>1}$, and applying the previous formula to $\alpha^n=\beta^{mn}$ would lead to $m\le 1$, a contradiction).
\end{remark}

\subsubsection*{Description and justification of the algorithm}

\paragraph{Combinatorial action.}

The elements of $\Free_p$ correspond to words on the alphabet $\tau^\pm$ that are \emph{reduced} in the sense that consecutive letters are not mutually inverse.

The Cayley graph of $(\Free_p,\tau^\pm)$ is the infinite $2p$-regular tree $\Tree_{p}$ with a base vertex $\Id$, whose edges around each vertex are labeled by $\tau^\pm$.
We identify the boundary $\partial \Tree_{p}$ with the subset of one sided sequences $(\tau^\pm)^\N$ which are reduced.
The cyclic order on $\tau^\pm$ extends to a cyclic order $\cord(x,y,z)$ on the boundary $\partial \Tree_{p}$.

The group $\Free_p$ acts by left translation on its Cayley graph $\Tree_p$, thus on its boundary $\partial \Tree_p$.
The action of $\gamma\in \Free_p\setminus\{\Id\}$ has attractive and repulsive fixed points $\gamma^+,\gamma^-\in \partial \Tree_p$, which correspond to the periodization of the words $\gamma^{+1}$ and $\gamma^{-1}$.

\paragraph{Geometric action.}

Consider a collection of $2p$ disjoint geodesics in the hyperbolic plane $\mathbb{HP}$ with distinct endpoints, labeled in cyclic order by $\tau^\pm$ (as we would obtain by lifting the spanning tree $T^*\subset \Surf_P$ from Figure \ref{fig:dual-tree-polygon-chordiag} to the universal cover of $\Surf_P$).
Let $\tau_j$ act as the hyperbolic isometry sending the geodesic $\tau_j^{-1}$ to the geodesic $\tau_j^{+1}$ to obtain a free action of $\Free_p$ on $\mathbb{HP}$, with fundamental domain an open polygon bounded by the geodesics $\tau_j^\pm$.

This faithfully represents $\Free_p$ as a subgroup of $\operatorname{Isom}^+(\mathbb{HP})$, its Cayley graph $\Tree_{p}$ as the dual tree of the tessellation by translates of the fundamental polygon in $\mathbb{HP}$, and its boundary $\partial \Tree_{p}$ as a Cantor subspace in the circle boundary $ \partial \mathbb{HP}$.

A nontrivial $\gamma\in \Free_p$ acts by hyperbolic translation on $\Tree_p\subset \mathbb{HP}$ along an axis $(\gamma^-,\gamma^+)$ joining its repulsive and attractive fixed points in the boundary $\partial \Tree_p \subset \mathbb{HP}$.

\begin{figure}[h]
	\centering
	\scalebox{1}{\input{images/tikz/lifts_crossing_domain4.tex}}
	\caption{Accompanying Example \ref{ex:counting_cross_val}: The red and blue loops in Figure \ref{fig:filling-punctures} are represented by the conjugacy classes of $\gamma_1=ab$ and $\gamma_2=aBa$, and their lifts in the universal cover of $\Sph^2_P$ which intersect the fundamental domain are the axes of the cyclically reduced representatives.}
    \label{fig:lifts_intersecting_domain}
\end{figure}

\paragraph{Cyclically reduced words.}

Recall that elements of $\Free_p$ correspond to reduced words on the alphabet $\tau^\pm$.
A word on $\tau_\pm$ is called \emph{cyclically reduced} when all its cyclic permutations are reduced.
Geometrically, the reduced word associated to a hyperbolic $\gamma\in \Free_p$ is cyclically reduced if and only if its hyperbolic axis intersects the fundamental region.
Hence a (nontrivial) conjugacy class corresponds in $\Free_p$ to the cyclic permutations of a cyclically reduced word on $\tau^\pm$, namely the representatives of its hyperbolic translations whose axis intersects the fundamental domain. 
We will denote by $\sigma$ the cyclic shift of words which moves the first letter at the end. 
The following is straightforward.

\textsc{Algorithm 1:}
For $\gamma\in \Free_p$ represented as a word on $\tau^\pm$ of length $\len(\gamma)$, we can find a cyclically reduced representative for its conjugacy class in $O(\len(\gamma))$ steps.

\paragraph{Primitive roots.}

For an element $\gamma\in \Free_p$, its primitivity exponent is the largest $n\in \N_{>0}$ for which there exists $\gamma_0\in \Free_p$ such that $\gamma=\gamma_0^n$: such a $\gamma_0$ is unique and called the primitive root of $\gamma$.
A nontrivial conjugacy class in $\Free_p$ is primitive when the cyclically reduced words that represent it are primitive words.
The following is classical.

\textsc{Algorithm 2:}
For a cyclically reduced word $\gamma$ on $\tau^\pm$ of length $\len(\gamma)$, we can compute its primitive root and primitivity exponent in time $O(\len(\gamma))$. 

\paragraph{Intersection number.}

For nontrivial $\alpha,\beta \in \Free_p$ acting on $\Tree_p\subset \mathbb{HP}$ with axes $(\alpha^-,\alpha^+)$ and $(\beta^-, \beta^+)$, define their crossing number $\cross(\alpha,\beta)\in \{-1,0,+1\}$ by:
\begin{equation*}
    \cross(\alpha,\beta)=\tfrac{1}{2}\left(\cord(\alpha^+,\beta^+, \alpha^-)-\cord(\alpha^+, \beta^-, \alpha^-)\right).
\end{equation*} 
This can be computed from the reduced words $\alpha,\beta \in \Free_n$ in time $O(\len(\alpha)+\len(\beta))$ as follows. 
First conjugate $\alpha,\beta$ by the longest $w\in \Free_n$ such that $\alpha,\beta$ have $w$ as a common prefix and $w$ (thus $w$ is trivial if either $\alpha$ or $\beta$ is cyclically reduced).
Next determine $\cord(\alpha^{+},\beta^{\pm}, \alpha^{-})$ by observing the first index where $\beta^\pm$ differs from both $\alpha^{+}$ and $\alpha^-$, which by periodicity must be smaller than $\len(\beta)$ unless $\beta^\pm$ coincides with either $\alpha^+$ or $\alpha^-$.

For nontrivial primitive $\alpha, \beta\in \Free_p$, define their \emph{order of contact} $\operatorname{val}(\alpha,\beta)$ as the number of translates of the fundamental domain that they both intersect (the union of which is convex hence contractible).
If $\alpha, \beta$ correspond to cyclically reduced primitive words, then denoting $w^\pm$ the longest common prefix of the infinite words $\alpha^{\pm}, \beta^{\pm}\in \partial \Free_p$
\begin{equation*}
    \operatorname{val}(\alpha,\beta)=1+\len(w^+)+\len(w^-)
\end{equation*} 
This can be computed from the finite words $\alpha,\beta$ in time $O(\len(\alpha)+\len(\beta))$. 
Indeed, we have $\len(w^\pm) \ge \len(\alpha)+\len(\beta)$ if and only if $\alpha, \beta$ have the same primitive roots, hence are equal by the primitivity assumption.

\begin{proposition}[intersection]
\label{prop:crossval}
For nontrivial primitive conjugacy classes in the cyclically ordered group $(\Free_p,\tau^\pm)$ represented by cyclically reduced words $\alpha,\beta$ on $\tau^\pm$, the intersection number of the corresponding curves is given by:
\begin{equation}
\label{eq:inter-crossval}
\tag{cross-val}
    \ti_P(\alpha,\beta) = 
    \sum_{i=1}^{\len(\alpha)}
    \sum_{j=1}^{\len(\beta)}
    \frac{\left|\cross(\sigma^i\alpha,\sigma^j\beta)\right|}{\val(\sigma^i\alpha,\sigma^j\beta)} 
\end{equation}
counting the pairs of cyclically reduced representatives of the conjugacy classes up to simultaneous conjugacy, whose endpoints are linked.
\end{proposition}

\begin{proof}
This is a reformulation of \cite[§1, §2]{Cohen-Lustig_algo-geometric-intersection-number_1987} and \cite[§3]{Birman-Series_algo-simple-closed-curves_1984}, so we briefly sketch the argument.

The topological intersection number $\ti_P(\alpha,\beta)$ counts the number of intersection points in the fundamental domain between representatives for the conjugacy classes (see Figure \ref{fig:lifts_intersecting_domain}). 
It is equal to the sum over pairs of representatives whose axes both enter the fundamental domain (given by cyclically reduced words $\sigma^i\alpha, \sigma^j\beta$) and which intersect (that is with $\lvert \cross(\sigma^i\alpha,\sigma^j\beta)\rvert = 1$) the inverse of the number of fundamental regions that they traverse in common (equal to $\operatorname{val}(\sigma^i\alpha, \sigma^j \beta)$).
This formula holds even when $\alpha=\beta$, and yields $\ti_P(\alpha,\alpha)=2\si_P(\alpha)$.
\end{proof}

\begin{example}[computing intersections]
\label{ex:counting_cross_val}
In Figure \ref{fig:lifts_intersecting_domain}, the red and blue loops from Figure \ref{fig:filling-punctures} are represented by the conjugacy classes of $\alpha=\gamma_1=ab$ and $\beta=\gamma_2=aBa$ with $\len(\alpha)=2$ and $\len(\beta)=3$, and their lifts in the universal cover of $\Sph^2_P$ which intersect the fundamental domain are the axes of the cyclically reduced representatives. Those satisfy:
\begin{itemize}[noitemsep]
    \item[$\alpha\beta$:] $\left|\cross(\sigma^i\alpha,\sigma^j\beta)\right|=1$ for $(i,j)\in \{(0,0),(0,2),(1,0)\} \subset \Z/2\times \Z/3$ and $0$ otherwise,
    \item[] $\val(\sigma^i\alpha,\sigma^j\beta)=2$ for all $(i,j)\in \Z/2\times \Z/3$.
\end{itemize}
\begin{itemize}[noitemsep]
    \item[$\alpha\alpha$:] $\left|\cross(\sigma^i\alpha,\sigma^j\alpha)\right|=1$ when $i\neq j$ in $\Z/2$ and $0$ otherwise,
    \item[] $\val(\sigma^i\alpha,\sigma^j\alpha)=1$ when $i\neq j$ in $\Z/2$ and $\infty$ otherwise.
\end{itemize}
\begin{itemize}[noitemsep]
    \item[$\beta\beta$:] $\left|\cross(\sigma^i\beta,\sigma^j\beta)\right|=0$ when $i=j$ or $(i,j)\in \{(1,2),(2,1)\} \subset \Z/3\times \Z/3$ and $1$ otherwise,
    \item[] $\val(\sigma^i\beta,\sigma^j\beta)=2$ when $i\neq j$ in $\Z/3$ and $\infty$ otherwise.
\end{itemize}
As expected, Equation \eqref{eq:inter-crossval} yields $\ti_P(\alpha,\beta)=\ti_P(\alpha,\alpha)=\ti_P(\beta,\beta)=2$. \end{example}

\textsc{Algorithm 3:}
For primitive $\alpha,\beta\in \Free_p\setminus\{\Id\}$ represented as cyclically reduced words on the symmetric set of generators $\tau^\pm$, we can compute the topological intersection number $\ti_P(\alpha,\beta)$ of the conjugacy classes in time $O(\len(\alpha)\len(\beta))$.

\subsubsection*{Consequences and summaries of algorithms for pinning multiloops}

We may now apply the algorithms from subsection \ref{subsec:MuLoops-to-FreeGroups} and \ref{subsec:si-algorithm} to the pinning problem.

\begin{algorithm}[self-intersection]
\label{algo:self-intersection}
    We have described an algorithm which has:
    \begin{itemize}[noitemsep, align=left]
        \item[Input:] A filling multiloop $\gamma\colon \sqcup_1^s \Sph^1 \looparrowright \Surf$ together with a subset of regions $P\subset \Regions$.
        \item[Output:] The self-intersection number of the multicurve $\gamma \colon \sqcup_1^s \Sph^1 \looparrowright \Surf_P$.
        \item[Complexity:] Running time $O(\Card(\Regions)^2)$.
    \end{itemize}
\end{algorithm}

\begin{remark}[The case $P=\emptyset$]

We presented Algorithm \ref{algo:Intro:self-intersection} so that it holds for $P\ne \emptyset$, but it can be adapted to compute the self-intersection number of multicurves in closed surfaces.

In the closed sphere, any multicurve has self-intersection $0$.
In the closed torus, every loop is homotopic to a power of a simple loop, and the self-intersection number of a loops follows as in Remark \ref{rem:non_primitive_si} whereas the intersection number of loops is the absolute value of their algebraic intersection number.

The case of closed surfaces of genus $g\ge 2$ is more subtle, but our description of Algorithm \ref{algo:self-intersection} can be extended by adapting \cite{Birman-Series_Dehn-algo-simple-closed-curves_1987, Lustig_Paths-geodesics-geometric-intersection-II_1987}.
We refer to \cite{Despres-Lazarus_Computing-intersection_2019} for a discussion of their work, and another approach to Algorithm \ref{algo:self-intersection}.
\end{remark}

From now on, we assume that we can deal with all cases, but we do insist that all algorithms should take care of the special cases individually and appropriately.
In particular, the following cautionary example shows that one cannot reduce the to the case $P\ne \emptyset$ so easily.

\begin{example}[loops that increase $\si$ with any puncture]
    \label{ex:increasing-si-any-puncture}
    Note that in $\Sph^2$, every multiloop $\gamma$ is homotopically trivial and remains so even after placing a puncture in any region, so for any one of its regions $\infty \in \Regions$ we have $\si_{\{\infty\}}(\gamma)=0$.
    
    In a surface $\Surf$ which is not simply connected, for all $n\in \N^*$ we can construct a loop $\gamma \subset \Surf$ which is homotopically trivial but such that for any region $\infty \in \Regions$ we have $\si_{\{\infty\}}(\gamma)\ge n$ (in particular $\gamma$ is non-homotopically trivial in $\Surf_{\infty}$).

    Let us describe our construction step by step in the torus.
    
    We first depict in the left of Figure \ref{fig:torus_commutator} a loop in the punctured torus which represent the commutator: the punctured torus is constructed as the quotient of the plane minus the interiors of two discs by identifying their circular boundaries labeled $C$ with the indicated orientations, and the purple loop represents the commutator. This loop becomes homotopically trivial in the torus (after filling the puncture at infinity). The right of Figure \ref{fig:torus_commutator} depicts that same loop in the square torus (where $C$ is now vertical) punctured at the large central square.
    
    Now cover the torus by $4$ closed rectangles $\{R_1,R_2,R_3,R_4\}$ as in the leftmost Figure \ref{fig:torus_multiloop}.
    For each $R_i$, construct a loop $\gamma_i$ in the complement which represents the commutator of that punctured torus as in the middle of Figure \ref{fig:torus_multiloop}. Choosing $\epsilon$ small enough ensures that their intersections are controlled in small annular neighborhoods of each rectangle.
    Finally, join them into one loop by resolving and clasping some intersection points in some of these annular neighborhoods as indicated at right in Figure \ref{fig:torus_multiloop} (for example, once between $\gamma_1$ and $\gamma_2$, once between $\gamma_1$ and $\gamma_3$, and once between $\gamma_1$ and $\gamma_4$). The resolutions perform the product of the loops based at the intersection points, and the clasping avoids any merging of regions around the intersection points.

    The regions of the resulting loop may belong to $1$ or $2$ or $3$ rectangles, and we may check in each case that after placing a puncture, we obtain a loop whose self-intersection is at least $3$.

    This example can be modified by taking powers of commutators to construct loops $\gamma$ in the torus with $\si_\emptyset(\gamma)=0$ but $\si_P(\gamma)$ arbitrarily large for every $P\neq\emptyset$.
    It can also be adapted to yield such examples in higher genus.
\begin{figure}[H]
    \centering
    \hspace{-0.5cm}
    \scalebox{0.3}{\definecolor{ca700f1}{RGB}{167,0,241}

\def \globalscale {1.000000}
\begin{tikzpicture}[y=1cm, x=1cm, yscale=\globalscale,xscale=\globalscale, every node/.append style={scale=\globalscale}, inner sep=0pt, outer sep=0pt]
  \begin{scope}[shift={(2.8547, 5.9196)}]
    \path[draw=ca700f1,line cap=butt,line join=miter,line width=0.0865cm] (19.4733, 3.1087).. controls (21.3891, 1.1929) and (20.5425, -0.8281) .. (18.6267, -2.7439).. controls (17.1898, -4.1807) and (16.1013, -4.3682) .. (14.224, -3.5906).. controls (9.7991, -1.7577) and (8.4889, 3.3078) .. (4.064, 5.1407).. controls (2.1867, 5.9183) and (0.7807, 5.7309) .. (-0.6562, 4.294).. controls (-3.5298, 1.4203) and (-3.5298, -0.2089) .. (-0.6562, -3.0826).. controls (2.2175, -5.9562) and (5.7654, -5.5939) .. (10.16, -5.5939).. controls (14.5546, -5.5939) and (18.1025, -5.9562) .. (20.9762, -3.0826).. controls (23.8498, -0.2089) and (23.8498, 1.4203) .. (20.9762, 4.294).. controls (19.5393, 5.7309) and (18.1333, 5.9183) .. (16.256, 5.1407).. controls (11.8311, 3.3078) and (10.5209, -1.7577) .. (6.096, -3.5906).. controls (4.2187, -4.3682) and (3.1302, -4.1807) .. (1.6933, -2.7439).. controls (-0.2225, -0.8281) and (-1.0691, 1.1929) .. (0.8467, 3.1087).. controls (5.6361, 7.8982) and (14.6839, 7.8982) .. (19.4733, 3.1087) -- cycle;

    \path[draw=black,even odd rule,line width=0.14cm] (4.064, 0.4734) ellipse (1.962cm and 1.962cm);

    \path[draw=black,line cap=round,line join=miter,line width=0.14cm,miter limit=4.0] (3.8574, 2.6899) -- (4.1114, 2.4359) -- (3.8574, 2.1819);

    \path[draw=blue,line cap=butt,line join=miter,line width=0.0865cm] (10.16, -3.5906).. controls (10.16, -3.5906) and (8.6869, 1.6404) .. (5.8629, 2.7582).. controls (3.0389, 3.8759) and (-0.3436, 1.8362) .. (2.4054, -1.4827).. controls (4.572, -4.0986) and (10.16, -3.5906) .. (10.16, -3.5906) -- cycle;

    \path[draw=blue,line cap=round,line join=miter,line width=0.0865cm,miter limit=4.0] (8.4107, -3.8842) -- (8.1887, -3.6018) -- (8.4711, -3.3798);

    \path[draw=ca700f1,line cap=round,line join=miter,line width=0.0865cm,miter limit=4.0] (9.9492, -5.3687) -- (10.2125, -5.6131) -- (9.9682, -5.8764);

    \tikzstyle{every node}=[font=\fontsize{30}{30}\selectfont]
    \node[text=red,line width=0.0265cm,anchor=south west] (text10) at (11.2115, -3.1318){$\alpha$};

    \node[text=blue,line width=0.0265cm,anchor=south west] (text11) at (8.5593, -3.2824){$\beta$};

    \node[text=ca700f1,line width=0.0265cm,anchor=south west] (text13) at (8.7837, -5.3206){$[\alpha,\beta]$};

    \path[draw=black,line cap=round,line join=miter,line width=0.14cm,miter limit=4.0] (16.4626, 2.6899) -- (16.2086, 2.4359) -- (16.4626, 2.1819);

    \path[draw=black,even odd rule,line width=0.14cm] (16.256, 0.4734) ellipse (1.962cm and 1.962cm);

    \node[line width=0.0265cm,anchor=south west] (text14) at (3.5418, 1.1439){$C$};

    \node[line width=0.0265cm,anchor=south west] (text15) at (15.9348, 1.1626){$C$};

    \path[draw=red,line cap=butt,line join=miter,line width=0.0865cm] (10.16, -3.5906).. controls (10.16, -3.5906) and (11.6331, 1.6404) .. (14.4571, 2.7582).. controls (17.2811, 3.8759) and (20.6636, 1.8362) .. (17.9146, -1.4827).. controls (15.748, -4.0986) and (10.16, -3.5906) .. (10.16, -3.5906) -- cycle;

    \path[draw=red,line cap=round,line join=miter,line width=0.0865cm,miter limit=4.0] (11.9093, -3.8842) -- (12.1313, -3.6018) -- (11.8489, -3.3798);

    \path[fill=black,even odd rule,line cap=butt,line join=miter,line width=0.0565cm,miter limit=4.0] (10.16, -3.5906) circle (0.3519cm);

  \end{scope}

\end{tikzpicture}}
    \hspace{-0.3cm}
    \scalebox{0.3}{\definecolor{ca700f1}{RGB}{167,0,241}
\definecolor{c3100bb}{RGB}{49,0,187}
\definecolor{navy}{RGB}{0,0,128}

\def \globalscale {1.000000}
\begin{tikzpicture}[y=1cm, x=1cm, yscale=\globalscale,xscale=\globalscale, every node/.append style={scale=\globalscale}, inner sep=0pt, outer sep=0pt]
  \begin{scope}[shift={(0.3728, -0.3278)}]
    \path[draw=ca700f1,line cap=round,line join=miter,line width=0.0865cm,miter limit=4.0] (2.54, 0.6555).. controls (2.54, 0.6555) and (0.5943, 5.6948) .. (2.54, 7.6405).. controls (4.3361, 9.4366) and (8.3639, 9.4366) .. (10.16, 7.6405).. controls (12.1057, 5.6948) and (10.16, 0.6555) .. (10.16, 0.6555);

    \path[draw=ca700f1,line cap=round,line join=miter,line width=0.0865cm,miter limit=4.0] (2.54, 13.3555).. controls (2.54, 13.3555) and (2.6255, 12.688) .. (3.5983, 11.7151).. controls (5.1341, 10.0847) and (11.7026, 7.8291) .. (8.89, 4.8889).. controls (7.2436, 3.2425) and (5.4564, 3.2425) .. (3.81, 4.8889).. controls (1.2129, 7.8868) and (7.3285, 10.1184) .. (9.1017, 11.7151).. controls (10.0745, 12.688) and (10.16, 13.3555) .. (10.16, 13.3555);

    \path[draw=black,line cap=round,line join=miter,line width=0.1065cm] (0.0, 13.3555) -- (0.0, 7.0055) -- (0.0, 0.6555);

    \path[draw=black,line cap=round,line join=miter,line width=0.1065cm] (12.7, 13.3555) -- (12.7, 7.0055) -- (12.7, 0.6555);

    \path[draw=black,line cap=round,line join=miter,line width=0.1065cm] (0.0, 13.3555) -- (6.35, 13.3555) -- (12.7, 13.3555);

    \path[draw=black,line cap=round,line join=miter,line width=0.1065cm] (0.0, 0.6555) -- (6.35, 0.6555) -- (12.7, 0.6555);

    \path[draw=ca700f1,line cap=round,line join=miter,line width=0.0865cm,miter limit=4.0] (6.2345, 3.9893) -- (6.552, 3.6718) -- (6.2345, 3.3543);

    \path[draw=black,line cap=round,line join=miter,line width=0.1065cm,miter limit=4.0] (6.1802, 13.673) -- (6.4977, 13.3555) -- (6.1802, 13.038);

    \path[draw=black,line cap=round,line join=miter,line width=0.1065cm,miter limit=4.0] (6.1802, 0.338) -- (6.4977, 0.6555) -- (6.1802, 0.973);

    \path[draw=black,line cap=round,line join=miter,line width=0.1065cm,miter limit=4.0] (-0.3285, 7.0055) -- (-0.011, 7.323) -- (0.3065, 7.0055);

    \path[draw=black,line cap=round,line join=miter,line width=0.1065cm,miter limit=4.0] (-0.3285, 6.688) -- (-0.011, 7.0055) -- (0.3065, 6.688);

    \path[draw=black,line cap=round,line join=miter,line width=0.1065cm,miter limit=4.0] (12.3715, 6.9945) -- (12.689, 7.312) -- (13.0065, 6.9945);

    \path[draw=black,line cap=round,line join=miter,line width=0.1065cm,miter limit=4.0] (12.3715, 6.677) -- (12.689, 6.9945) -- (13.0065, 6.677);

    \path[draw=c3100bb,fill=navy,fill opacity=0.2289,even odd rule,line cap=round,line join=miter,line width=0.2165cm,miter limit=4.0] (6.35, 7.0055) circle (0.368cm);

    \tikzstyle{every node}=[font=\fontsize{30}{30}\selectfont]
    \node[text=ca700f1,line width=0.0265cm,anchor=south west] (text13) at (5.1324, 4.1536){$[\alpha,\beta]$};

  \end{scope}

\end{tikzpicture}}
    \caption{A loop in a punctured torus representing the commutator $[\alpha,\beta]=\alpha\beta\alpha^{-1}\beta^{-1}$.}
    \label{fig:torus_commutator}
\end{figure}
\begin{figure}[H]
    \centering
    \hspace{-0.5cm}
    \scalebox{0.5}{\definecolor{cc8c8ff}{RGB}{200,200,255}
\definecolor{cffc8c8}{RGB}{255,200,200}
\definecolor{cc8ffc8}{RGB}{200,255,200}
\definecolor{cc8c8c8}{RGB}{200,200,200}
\definecolor{lime}{RGB}{0,255,0}

\def \globalscale {1.000000}
\begin{tikzpicture}[y=1cm, x=1cm, yscale=\globalscale,xscale=\globalscale, every node/.append style={scale=\globalscale}, inner sep=0pt, outer sep=0pt]
  \begin{scope}[shift={(0.3728, -0.3278)}]
    \path[fill=cc8c8ff,even odd rule,line cap=butt,line join=round,line width=0.5cm,miter limit=4.0] (0.0, 13.3555) rectangle (3.81, 9.5455);

    \path[fill=cc8c8ff,even odd rule,line cap=butt,line join=round,line width=0.5cm,miter limit=4.0] (8.89, 13.3555) rectangle (12.7, 9.5455);

    \path[fill=cc8c8ff,even odd rule,line cap=butt,line join=round,line width=0.5cm,miter limit=4.0] (0.0, 4.4655) rectangle (3.81, 0.6555);

    \path[fill=cc8c8ff,even odd rule,line cap=butt,line join=round,line width=0.5cm,miter limit=4.0] (8.89, 4.4655) rectangle (12.7, 0.6555);

    \path[fill=cffc8c8,even odd rule,line cap=butt,line join=round,line width=0.5cm,miter limit=4.0] (3.81, 13.3555) rectangle (8.89, 10.8155);

    \path[fill=cffc8c8,even odd rule,line cap=butt,line join=round,line width=0.5cm,miter limit=4.0] (3.81, 3.1955) rectangle (8.89, 0.6555);

    \path[fill=cc8ffc8,even odd rule,line cap=butt,line join=round,line width=0.5cm,miter limit=4.0,cm={ 0.0,1.0,1.0,0.0,(-13.3555, 13.3555)}] (-8.89, 15.8955) rectangle (-3.81, 13.3555);

    \path[fill=cc8c8c8,fill opacity=0.3793,even odd rule,line join=round,line width=0.5cm] (2.54, 10.8155) rectangle (10.16, 3.1955);

    \path[fill=cc8ffc8,even odd rule,line cap=butt,line join=round,line width=0.5cm,miter limit=4.0,cm={ 0.0,1.0,1.0,0.0,(-13.3555, 13.3555)}] (-8.89, 26.0555) rectangle (-3.81, 23.5155);

    \path[draw=blue,even odd rule,line cap=butt,line join=round,line width=0.0465cm,miter limit=4.0,dash pattern=on 0.0465cm off 0.093cm] (8.5725, 13.3555) rectangle (12.7, 9.228);

    \path[draw=black,even odd rule,line join=round,line width=0.0465cm,dash pattern=on 0.0465cm off 0.093cm] (2.2225, 11.133) rectangle (10.4775, 2.878);

    \path[draw=blue,even odd rule,line cap=butt,line join=round,line width=0.0465cm,miter limit=4.0,dash pattern=on 0.0465cm off 0.093cm] (0.0, 13.3555) rectangle (4.1275, 9.228);

    \path[draw=blue,even odd rule,line cap=butt,line join=round,line width=0.0465cm,miter limit=4.0,dash pattern=on 0.0465cm off 0.093cm] (0.0233, 4.783) rectangle (4.1275, 0.6788);

    \path[draw=blue,even odd rule,line cap=butt,line join=round,line width=0.0465cm,miter limit=4.0,dash pattern=on 0.0465cm off 0.093cm] (8.5725, 4.783) rectangle (12.6768, 0.6788);

    \path[draw=red,even odd rule,line cap=butt,line join=round,line width=0.0465cm,miter limit=4.0,dash pattern=on 0.0465cm off 0.093cm] (3.4925, 13.3323) rectangle (9.2075, 10.498);

    \path[draw=red,even odd rule,line cap=butt,line join=round,line width=0.0465cm,miter limit=4.0,dash pattern=on 0.0465cm off 0.093cm] (3.5125, 3.513) rectangle (9.2275, 0.6555);

    \path[draw=lime,even odd rule,line cap=butt,line join=round,line width=0.0465cm,miter limit=4.0,dash pattern=on 0.0465cm off 0.093cm,rotate around={-90.0:(0.0, 13.3555)}] (3.4925, 16.213) rectangle (9.2075, 13.3555);

    \path[draw=lime,even odd rule,line cap=butt,line join=round,line width=0.0465cm,miter limit=4.0,dash pattern=on 0.0465cm off 0.093cm,rotate around={-90.0:(0.0, 13.3555)}] (3.4925, 26.0323) rectangle (9.2075, 23.198);

    \path[draw=black,line cap=round,line join=miter,line width=0.1065cm] (0.0, 13.3555) -- (0.0, 7.0055) -- (0.0, 0.6555);

    \path[draw=black,line cap=round,line join=miter,line width=0.1065cm] (12.7, 13.3555) -- (12.7, 7.0055) -- (12.7, 0.6555);

    \path[draw=black,line cap=round,line join=miter,line width=0.1065cm] (0.0, 13.3555) -- (6.35, 13.3555) -- (12.7, 13.3555);

    \path[draw=black,line cap=round,line join=miter,line width=0.1065cm] (0.0, 0.6555) -- (6.35, 0.6555) -- (12.7, 0.6555);

    \begin{scope}[shift={(-0.1588, -0.1588)}]
      \tikzstyle{every node}=[font=\fontsize{20}{20}\selectfont]
      \node[text=black,line width=0.0265cm,anchor=south west] (text5) at (6.1579, 6.7963){$R_1$};

      \node[text=black,line width=0.0265cm,anchor=south west] (text22) at (6.1379, 12.1938){$R_4$};

      \node[text=black,line width=0.0265cm,anchor=south west] (text23) at (6.1579, 1.3988){$R_4$};

      \node[text=black,line width=0.0265cm,anchor=south west] (text24) at (1.0779, 1.3826){$R_2$};

      \node[text=black,line width=0.0265cm,anchor=south west] (text25) at (1.0779, 6.7801){$R_3$};

      \node[text=black,line width=0.0265cm,anchor=south west] (text26) at (1.0779, 12.1938){$R_2$};

      \node[text=black,line width=0.0265cm,anchor=south west] (text27) at (11.2579, 12.1938){$R_2$};

      \node[text=black,line width=0.0265cm,anchor=south west] (text28) at (11.2379, 6.7963){$R_3$};

      \node[text=black,line width=0.0265cm,anchor=south west] (text29) at (11.2379, 1.3988){$R_2$};

    \end{scope}
    \path[draw=black,line cap=round,line join=miter,line width=0.1065cm,miter limit=4.0] (6.1912, 13.673) -- (6.5087, 13.3555) -- (6.1912, 13.038);

    \path[draw=black,line cap=round,line join=miter,line width=0.1065cm,miter limit=4.0] (6.1912, 0.338) -- (6.5087, 0.6555) -- (6.1912, 0.973);

    \path[draw=black,line cap=round,line join=miter,line width=0.1065cm,miter limit=4.0] (-0.3175, 7.0055) -- (0.0, 7.323) -- (0.3175, 7.0055);

    \path[draw=black,line cap=round,line join=miter,line width=0.1065cm,miter limit=4.0] (-0.3175, 6.688) -- (0.0, 7.0055) -- (0.3175, 6.688);

    \path[draw=black,line cap=round,line join=miter,line width=0.1065cm,miter limit=4.0] (12.3825, 6.9945) -- (12.7, 7.312) -- (13.0175, 6.9945);

    \path[draw=black,line cap=round,line join=miter,line width=0.1065cm,miter limit=4.0] (12.3825, 6.677) -- (12.7, 6.9945) -- (13.0175, 6.677);

  \end{scope}

\end{tikzpicture}}
    \hspace{-0.3cm}
    \scalebox{0.5}{\definecolor{cffc8c8}{RGB}{255,200,200}

\def \globalscale {1.000000}
\begin{tikzpicture}[y=1cm, x=1cm, yscale=\globalscale,xscale=\globalscale, every node/.append style={scale=\globalscale}, inner sep=0pt, outer sep=0pt]
  \begin{scope}[shift={(0.3728, -0.3278)}]
    \path[draw=red,even odd rule,line join=round,line width=0.0565cm,dash pattern=on 0.0565cm off 0.1695cm] (1.27, 9.5455) rectangle (10.16, 3.1955);

    \path[draw=black,line cap=round,line join=miter,line width=0.1065cm] (0.0, 13.3555) -- (0.0, 7.0055) -- (0.0, 0.6555);

    \path[draw=black,line cap=round,line join=miter,line width=0.1065cm] (12.7, 13.3555) -- (12.7, 7.0055) -- (12.7, 0.6555);

    \path[draw=black,line cap=round,line join=miter,line width=0.1065cm] (0.0, 13.3555) -- (6.35, 13.3555) -- (12.7, 13.3555);

    \path[draw=black,line cap=round,line join=miter,line width=0.1065cm] (0.0, 0.6555) -- (6.35, 0.6555) -- (12.7, 0.6555);

    \path[fill=cffc8c8,even odd rule,line cap=butt,line join=round,line width=0.1065cm,miter limit=4.0] (2.54, 8.2755) rectangle (8.89, 4.4655);

    \tikzstyle{every node}=[font=\fontsize{20}{20}\selectfont]
    \node[text=black,line cap=butt,line join=miter,line width=0.0265cm,miter limit=4.0,anchor=south west] (text5) at (5.4891, 6.2231){$R_i$};

    \path[draw=red,line cap=butt,line join=miter,line width=0.0565cm] (1.905, 0.6555) -- (1.905, 8.9105) -- (9.525, 8.9105) -- (9.525, 0.6555);

    \path[draw=red,line cap=butt,line join=miter,line width=0.0565cm] (1.905, 13.3555) -- (1.905, 9.228) -- (5.08, 9.228) -- (9.2075, 8.593) -- (9.2075, 3.8305) -- (2.2225, 3.8305) -- (2.2225, 8.593) -- (6.35, 9.228) -- (9.525, 9.228) -- (9.525, 13.3555);

    \path[draw=red,line cap=butt,line join=miter,line width=0.0565cm,miter limit=4.0] (5.5562, 4.148) -- (5.8737, 3.8305) -- (5.5562, 3.513);

    \path[draw=black,line cap=round,line join=miter,line width=0.1065cm,miter limit=4.0] (6.1802, 13.673) -- (6.4977, 13.3555) -- (6.1802, 13.038);

    \path[draw=black,line cap=round,line join=miter,line width=0.1065cm,miter limit=4.0] (6.1802, 0.338) -- (6.4977, 0.6555) -- (6.1802, 0.973);

    \path[draw=black,line cap=round,line join=miter,line width=0.1065cm,miter limit=4.0] (-0.3285, 7.0055) -- (-0.011, 7.323) -- (0.3065, 7.0055);

    \path[draw=black,line cap=round,line join=miter,line width=0.1065cm,miter limit=4.0] (-0.3285, 6.688) -- (-0.011, 7.0055) -- (0.3065, 6.688);

    \path[draw=black,line cap=round,line join=miter,line width=0.1065cm,miter limit=4.0] (12.3715, 6.9945) -- (12.689, 7.312) -- (13.0065, 6.9945);

    \path[draw=black,line cap=round,line join=miter,line width=0.1065cm,miter limit=4.0] (12.3715, 6.677) -- (12.689, 6.9945) -- (13.0065, 6.677);

    \node[text=black,line width=0.0265cm,anchor=south west] (text5-5) at (8.2562, 3.372){$\epsilon$};

    \path[draw=black,line cap=butt,line join=miter,line width=0.0165cm,dash pattern=on 0.033cm off 0.0165cm] (8.723, 3.1955) -- (8.723, 4.4655);

    \path[draw=black,line cap=butt,line join=miter,line width=0.0165cm] (8.5642, 4.3068) -- (8.723, 4.4655) -- (8.8817, 4.3068);

    \path[draw=black,line cap=butt,line join=miter,line width=0.0165cm] (8.5642, 3.3543) -- (8.723, 3.1955) -- (8.8817, 3.3543);

    \node[text=red,line width=0.0265cm,anchor=south west] (text5-5-2) at (4.5986, 3.9023){$\gamma_i$};

  \end{scope}

\end{tikzpicture}}
    \scalebox{1}{\def \globalscale {1.000000}
\begin{tikzpicture}[y=1cm, x=1cm, yscale=\globalscale,xscale=\globalscale, every node/.append style={scale=\globalscale}, inner sep=0pt, outer sep=0pt]
  \begin{scope}[shift={(-5.08, -2.58)}]
    \path[draw=blue,line cap=butt,line join=miter,line width=0.0565cm] (5.08, 7.8187) -- (7.62, 7.8187);

    \path[draw=red,line cap=butt,line join=miter,line width=0.0565cm] (6.35, 6.5487) -- (6.35, 9.0887);

    \path[draw=red,line cap=butt,line join=miter,line width=0.0565cm,miter limit=4.0] (6.1913, 8.93) -- (6.35, 9.0887) -- (6.5088, 8.93);

    \path[draw=blue,line cap=butt,line join=miter,line width=0.0565cm,miter limit=4.0] (7.4613, 7.66) -- (7.62, 7.8187) -- (7.4613, 7.9775);

    \path[draw=black,line cap=butt,line join=miter,line width=0.0565cm] (6.35, 2.7387) -- (6.35, 3.0562) -- (6.0325, 3.3737) -- (6.0325, 4.0087) -- (6.35, 4.3262) -- (6.985, 4.3262) -- (7.3025, 4.0087) -- (7.62, 4.0087);

    \path[draw=black,line cap=butt,line join=miter,line width=0.0565cm] (6.35, 5.2787) -- (6.35, 4.9612) -- (6.6675, 4.6437) -- (6.6675, 4.0087) -- (6.35, 3.6912) -- (5.715, 3.6912) -- (5.3975, 4.0087) -- (5.08, 4.0087);

    \path[draw=black,line cap=butt,line join=miter,line width=0.0565cm] (6.1913, 5.12) -- (6.35, 5.2787) -- (6.5088, 5.12);

    \path[draw=black,line cap=butt,line join=miter,line width=0.0565cm] (7.4613, 3.85) -- (7.62, 4.0087) -- (7.4613, 4.1675);

    \path[draw=black,line cap=round,line join=miter,line width=0.0165cm] (6.35, 6.223) -- (6.35, 5.588);

    \path[draw=black,line cap=round,line join=miter,line width=0.0165cm] (6.5088, 5.7467) -- (6.35, 5.588) -- (6.1913, 5.7467);

    \path[draw=white,draw opacity=0.0,line cap=butt,line join=miter,line width=0.0565cm,miter limit=4.0] (6.35, 2.58) -- (6.35, 2.7387);

  \end{scope}

\end{tikzpicture}}
    \caption{Building a loop $\gamma$ in the torus with $\si_\emptyset(\gamma)=0$ but $\si_P(\gamma)\geq 3$ for any $P\neq\emptyset$.}
    \label{fig:torus_multiloop}
\end{figure}
\end{example}

\begin{corollary}[\textsc{Pinning} and \textsc{PinMin} are \textsf{P}]
\label{cor:MultiLooPinMin-P-easy}
Given a multiloop $\gamma$ and a subset of regions $P\subset \Regions$, we can successively:
    \begin{itemize}[noitemsep, align=left]
        \item[\textsc{Pinning}:] certify that $P$ is pinning in time $O(\Card(\Regions)^2)$, and if so
        \item[\textsc{PinMin}:] check if $P$ is a minimal pinning set in time $O(\Card(P)\Card(\Regions)^2)$, and if not 
        \item[\textsc{PinMinC}:] construct a minimal pinning set $M\subset P$ in time $O(\Card(P)\Card(\Regions)^2)$.
    \end{itemize}
\end{corollary}

\begin{proof}
    Consider a multiloop $\gamma\colon \sqcup_1^s \Sph^1 \looparrowright \Surf$ and a set of regions $P\subset \Regions$. 
    
    First, we may compute $\si_P(\gamma)$ using the previous algorithm in time $O(\Card(\Regions)^2)$ and compare it with $\#\gamma$. This checks whether $P$ is pinning and if so, confirms that $\varpi(\gamma)\le \Card(P)$. 

    Now suppose that we have checked that $P$ is pinning. To a linear order of $P$ we associate a minimal pinning subset computed as follows. 
    We construct a decreasing chain of pinning sets, starting with $P$, by trying to remove them one by one.
    We try removing the first region and check if the new set is pinning in time $O(\Card(\Regions)^2)$: if it is we continue with that new set, otherwise we move to the next region.
    By doing so, we never need to try removing a pin twice: either it was already removed, or we know it cannot be removed.
    Thus in the worst case we need to call at most $\Card(P)$ times to our $O(\Card(\Regions)^2)$ checking algorithm.
    Note that this construction defines a surjective function from the set of linear orders of $\Regions$ to the set of minimal pinning sets of $\gamma$.
\end{proof}

\begin{corollary}[\textsc{MultiLooPinNum} is \textsf{NP}]
\label{cor:MultiLooPinNum-is-NP}
    The \textsc{MultiLooPinNum} problem, defined by 
    \begin{itemize}[align=left, noitemsep]
        \item[Instance:] A filling multiloop $\gamma \colon \sqcup_1^s \Sph^1 \looparrowright \Surf$, and an integer $p\in \N$.
        \item[Question:] Does $\gamma$ have pinning number $\varpi(\gamma)\le p$?
    \end{itemize}
    belongs to the class \textsf{NP}.
\end{corollary}

\newpage

\section{Geometric topology: computing immersed discs}

This section focuses on loops, namely multiloops with one strand, maybe non-filling.

In the first subsection we define various notions of monogons and bigons, provide a sufficient condition for a subset of regions to be pinning which will be useful to justify the Algorithm \ref{algo:planar-vertex-cover-to-LooPinNum}, and prove the main Theorem \ref{thm:taut-no-immonorbigon}.

The second subsection applies Theorem \ref{thm:taut-no-immonorbigon} to put the pinning sets of a loop in correspondence with the vertex covers of a hypergraph or equivalently the satisfying assignments of a positive conjunctive normal form. 
Hence it is possible that the pinning number and minimal pinning sets of many loops can be found efficiently with certain \textsc{SAT}-solvers.

The third subsection will reduce the vertex cover problem for planar graphs to our \textsc{LooPinNum} problem, proving it is \textsf{NP}-hard (for any genus).

\subsection{Taut loops, their monorbigons and mobidiscs}
\label{subsec:monorbigons}

Fix an oriented surface $\Surf$, and recall that for a subset $P\subset \Surf$ we denote by $\Surf_P = \Surf\setminus P$. 
Consider a loop $\gamma \colon \Sph^1 \looparrowright \Surf_P$ which may be filling or not; denote by $\Regions=\Surf\setminus \gamma$ its set of regions and by $V\subset \Surf_P$ its set of double points.

\begin{definition}[singular monorbigons]
\label{def:singular-monorbigon}


Consider a loop $\gamma \colon \Sph^1 \looparrowright \Surf_P$.

A \emph{singular monogon} is a nontrivial closed interval $I\subset \Sph^1$ such that $\gamma(\partial I)=\{x\}$ for some $x\in V$ and $\gamma(I)$ is null-homotopic.
This singular monogon is \emph{embedded} when $\gamma$ is injective on the interior of $I$.

A \emph{singular bigon} is a disjoint union of two nontrivial closed intervals $I \sqcup J \subset \Sph^1$ such that $\gamma(\partial I)=\{x,y\}=\gamma(\partial J)$ for some $x,y\in V$ and $\gamma(I\sqcup J) \subset \Surf_P$ is null-homotopic.
This singular bigon is \emph{embedded} when $\gamma$ is injective on the interior of $I\sqcup J$.

A \emph{singular monorbigon} $K$ refers to a singular monogon $I$ or singular bigon $I\sqcup J$.
The restriction $\gamma(K)$ is a subloop of $\gamma$ and we call $\gamma(\partial K)$ its \emph{marked points}.
\end{definition}

\begin{figure}[h]
    \centering
     \scalebox{0.26}{\input{images/tikz/monorbigon_cases.tex}}
    \caption{Examples and non-examples of different kinds of monorbigons.}
    \label{fig:monorbigons}
\end{figure}

\begin{algorithm}[computing singular monorbigons]
\label{algo:singular-monorbigons}
For a filling loop $\gamma \colon \Sph^1 \looparrowright \Surf_P$, we can list all its singular monorbigons in time $O(\Card(\Regions)^3)$.
\end{algorithm}

\begin{proof}
First list the $2\Card(V)$ intervals $I$ and $< 2\Card(V)^2$ interval unions $I\sqcup J$ satisfying the conditions in definition \ref{def:singular-monorbigon}, except the homotopic triviality.
(Recall $\Card(V)=O(\Card(\Regions))$.)

To check if the subloop $\gamma(I)$ or $\gamma(I\sqcup J)$ is homotopically trivial: first compute the associated word in the fundamental group $\pi_1(\Surf_P)$, then apply a linear-time algorithm solving the word problem for that surface group. When $P\neq\emptyset$ this group is free, and when $P=\emptyset$ we add a relation which is a product of commutators and apply Dehn's algorithm \cite{Dehn_papers-group-topology_1987}. \end{proof}

Our interest in the concept of embedded and singular monorbigons lies in the following Theorems that were proven by Hass--Scott in \cite{Hass-Scott_Intersection-curves-surfaces_1985}.

\begin{theorem}[taut loops have no singular monorbigons]
\label{thm:Hass-Scott-2.7-4.2}
    Fix a loop $\gamma \colon \Sph^1 \looparrowright \Surf_P$, with $\#\gamma$ double points and whose homotopy class in $\Surf_P$ has self-intersection number $\si_P(\gamma)$.
    \begin{itemize}[noitemsep]
        \item[(2.7)] If $\# \gamma >\si_P(\gamma)=0$, then $\gamma$ has an embedded monorbigon in $\Surf_P$.
        \item[(4.2)] If $\#\gamma > \si_P(\gamma)$, then $\gamma$ has a singular monorbigon in $\Surf_P$.
    \end{itemize}
\end{theorem}

\begin{remark}[embedded monorbigons do not suffice]
 The ``Immersed'' column in Figure \ref{fig:monorbigons} makes it clear that one may not replace ``singular'' by ``embedded'' in part (4.2) of Theorem \ref{thm:Hass-Scott-2.7-4.2}.
\end{remark}

Now we wish to understand how to pin a singular monorbigon: which minimal collections of regions should be punctured in order to remove that singular monorbigon?

\begin{definition}[Linking function]
Consider a loop $\alpha\colon \Sph^1 \looparrowright \Surf_P$ with regions $\Regions$, and assume that it yields a trivial homology class in $H_1(\Surf_P;\Z)=H_1(\Surf;\Z)$.

For any two regions $o,\infty \in \Regions$, the linking number $\lk(\alpha, [o]-[\infty]) \in \Z$ is defined as the signed intersection number between $\alpha$ and any arc $\beta$ from $\infty$ to $o$.

The \emph{linking function} $\lk(\alpha,\cdot) \colon \Regions \times \Regions \to \Z$ defined by $(o,\infty)\mapsto \lk(\alpha, [o]-[\infty])$ changes sign under change of orientation of the $1$-cycle $\alpha$ or the $0$-chain $[o]-[\infty]$, in particular its absolute value and vanishing locus are well defined independently of these orientations.
\end{definition}

\begin{remark}[homology]
    \label{rem:link(singular-monorbigons)=0}
    For a loop $\gamma\colon \Sph^1\looparrowright \Surf_P$, a singular monorbigon $K$ yields a subloop $\alpha = \gamma(K)$ which is homotopically trivial in $\Surf_P$, hence trivial in homology. 
    
    In particular $\alpha$ is trivial in $H_1(\Surf_P;\Z/2)$, so its complement is properly $2$-colorable.
    Moreover $\alpha$ is trivial in $H_1(\Surf_P;\Z)$, so for all $p_0,p_1\in P$ we have $\lk(\alpha ; [p_0]-[p_1])=0$.
\end{remark}

\begin{lemma}[singular monorbigons link trivially]
\label{lem:link(singular-monorbigons)=0}
Consider a loop $\gamma \colon \Sph^1 \looparrowright\Surf$ and a set of regions $P\subset \Regions$.
If $P$ is not pinning $\gamma$, then there exists a singular monorbigon $\alpha$ of $\gamma$ such that for every pair of regions $o,\infty \in P$ we have $\lk(\alpha, [o]-[\infty]) = 0$.
\end{lemma}

\begin{proof}
If $P$ is not pinning $\gamma$ then the loop $\gamma_P \colon \Sph^1 \looparrowright \Surf_P\setminus \{\infty\}$ is not taut, so by \ref{thm:Hass-Scott-2.7-4.2} it has a singular monorbigon $\alpha$, so the statement follows from Remark \ref{rem:link(singular-monorbigons)=0}.
\end{proof}

\begin{remark}[converse is false]
There are loops $\gamma\colon \Sph^1 \looparrowright \Surf$ admitting pinning sets $P\subset \Regions$ which contradict the converse to the previous lemma: for every monorbigon $\alpha$ of $\gamma$ and every pair of pinned regions $o,\infty \in P$, we have $\lk(\alpha, [o]-[\infty])=0$.

One example is obtained as follows. 
The thrice punctured sphere $\Sph^2\setminus \{0,1,\infty\}$ with fundamental group $\langle \tau_0,\tau_1 \rangle$ has a unique hyperbolic metric: consider the geodesic representative $\gamma$ for the commutator $[\tau_0,\tau_1]$ as depicted in Figure \ref{fig:trefoil-commutator}.
As commutators have trivial abelianization, for all $p,q\in \{0,1,\infty\}$ we get $\lk(\gamma, [p]-[q])=0$.
\begin{figure}[h]
    \centering
    \includegraphics[width=0.4\linewidth]{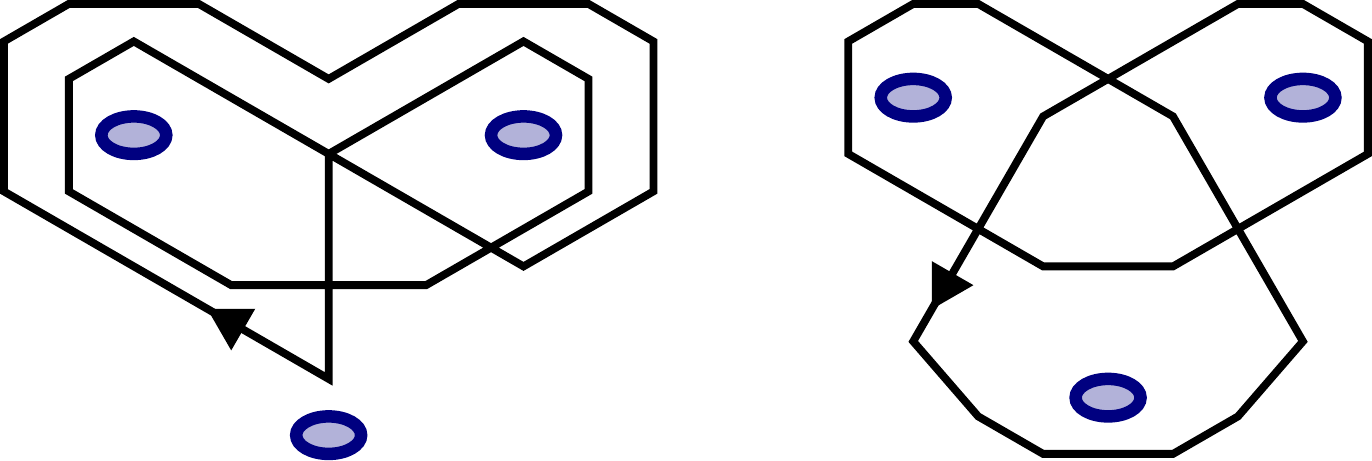}
    \caption{The trefoil in $\Sph^2\setminus\{0,1,\infty\}$ represents a commutator in $\pi_1(\Sph^2\setminus\{0,1,\infty\})$.}
    \label{fig:trefoil-commutator}
\end{figure}
\end{remark}

To characterize the pinning sets of loops in terms of the linking numbers of their singular monorbigons, we will thus restrict to the class of immersed monorbigons.

\begin{definition}[immersed monorbigons and mobidiscs]
\label{def:immersed-monorbigon-mobidiscs}
    Consider a loop $\gamma \colon \Sph^1 \looparrowright \Surf_P$. 
    
    We say that a singular monorbigon $K\subset \Sph^1$ is \emph{immersed} when there exists an immersion $\iota\colon \Disc\looparrowright \Surf_P$ such that the restriction $\gamma \colon K \looparrowright \Surf_P$ factors through $\iota \colon \partial \Disc \looparrowright \Surf_P$.
    A \emph{mobidisc} consists of the regions of $\gamma$ in the image $\iota(\Disc)$ of such an immersion.
    
    Denote by $\MoB(\gamma)\subset \Parts(\Regions)$ the set of all its mobidiscs, and its subset of proper mobidiscs by $\MoB(\gamma)\setminus \{\Regions\}$ (which may or may not coincide with $\MoB(\gamma)$).
\end{definition}

\begin{remark}[precisions]
    According to the definition of immersed monorbigons, the immersion $\iota\colon \partial \Disc\looparrowright \gamma(K)$ is surjective, and its only double points are those of $\gamma(K)$.
    
    In particular $\iota \colon \partial \Disc \looparrowright \gamma(K)$ cannot overlap itself (even with opposite directions), so the immersion suggested on the right of Figure \ref{fig:weak-bigon} does not satisfy the conditions in Definition \ref{def:immersed-monorbigon-mobidiscs}.
\begin{figure}[h]
    \centering
    \includegraphics[scale=0.43, angle=-90]{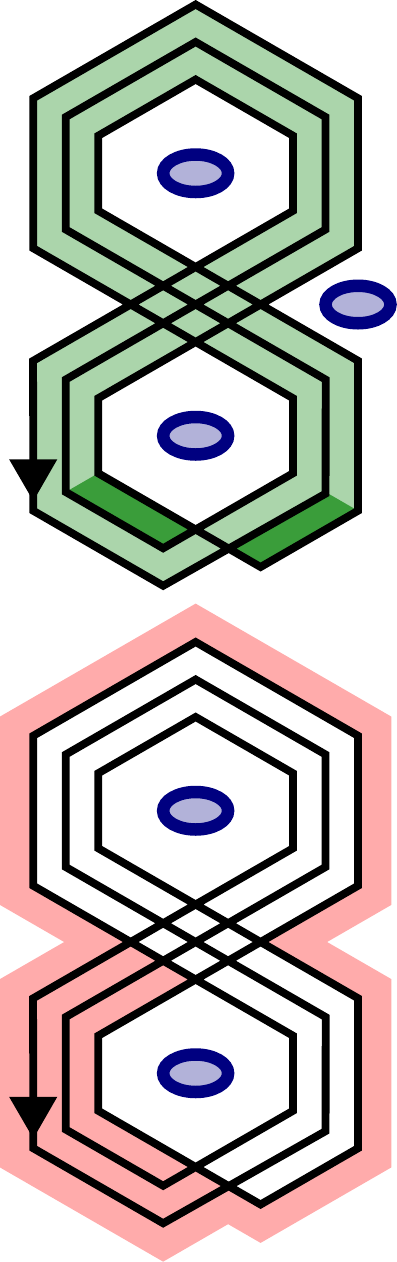}
    \caption{Left: The loop is not taut, since it has a singular bigon (which is immersed, and even embedded).
    Right: Placing an additional puncture at infinity makes the loop taut. 
    Note that the immersed disc covering the green portions does not satisfy the conditions in Definition \ref{def:immersed-monorbigon-mobidiscs}.}
    \label{fig:weak-bigon}
\end{figure}
\end{remark}

\begin{lemma}[uniqueness of mobidiscs]
    \label{lem:unique-monorbidisc}
    For a loop $\gamma \colon \Sph^1 \looparrowright \Surf_P$ with $\infty \in P$, every immersed monorbigon $K$ gives rise to exactly one mobidisc: it consists of the regions $o\in \Regions\setminus P$ such that $\lvert \lk(\gamma(K), [o]-[\infty]) \rvert \ne 0$.
\end{lemma}

\begin{proof}
    One may deduce this from the more general \cite[Theorem 3.2 and Corollary 3.3]{Ezell-Marx_Branched-extensions-curves-orientable-surfaces_1980} or the results in \cite{Pappas_Extensions-codim-1-immersions_1996}, but here is a simple explanation.
    
    Consider a region $o$ and an immersed disc $\iota\colon \Disc \mapsto \Surf_P$ such that the immersed monorbigon $\gamma\colon K\looparrowright \Surf_P$ factors through $\iota\colon \partial\Disc \looparrowright \Surf_P$.
    
    If $\lk(\gamma(K), [o]-[\infty]) \ne 0$ then $o$ must belong to $\im(\iota)$ since $\gamma(K)$ is null-homotopic inside $\im(\iota)$ (this holds even for a singular monorbigon bounding a singular disc).

    If $o$ belongs to $\im(\gamma)$ and $\lk(\gamma(K),[o]-[\infty])=0$, then there must exist two points in the preimage $\iota^{-1}(o)$ at which the Jacobian of $\iota$ has opposite determinants.
    Hence the Jacobian determinant of $\iota$ must vanish along any arc of $\Disc$ joining these two points, contradicting the immersion property.
\end{proof}

\begin{remark}[non-uniqueness of extensions] 
    Note however that an immersed monorbigon $K$ of $\gamma$ can extend to several non-isotopic immersions $\iota \colon \Disc \looparrowright \Surf$ which define the same mobidisc.
    The Figure \ref{fig:Milnor-paysley} depicts the example of such a plane loop which was famously discovered by Milnor.

    Given a loop in a surface, the description and computation of all its extensions to immersions of a disc was raised by Hopf and Thom \cite{Thom_inequadiff-globales_1959}.
    It has been solved in the plane, first by Blank \cite{Blank_extending-immersions-circle_1967, Poenaru_extension-immersions-codim-1_1995} and then by Shor and Van Wyk in \cite{Shor-VanWyk_Detecting-decomposing-self-overlapping-curves_1992}.
    The status for other surfaces remains unclear: Frisch worked on generalizations to oriented surfaces in his Thesis  \cite{Frisch_classification-immersions-curves-in-surfaces_2010}, but restricted to the case of the sphere in his preprint \cite{Frisch_extending-immersions-sphere_2010}, and none of those have been published. 
    More general problems are broached in \cite{Ezell-Marx_Branched-extensions-curves-orientable-surfaces_1980} and \cite{Pappas_Extensions-codim-1-immersions_1996} using different methods.
\end{remark}

\begin{figure}[h]
    \centering
    \scalebox{0.38}{\input{images/tikz/milnordoodle.tex}}
    \caption{The ``Milnor doodle'' bounds an immersed disc in two different ways. Viewed from above, the violet portions are the closest whereas the red portions are the farthest. The number in each region correspond to the local degree of the immersion (they agree in both cases).}
    \label{fig:Milnor-paysley}
\end{figure}

We may now formulate the following improvement of Theorem \ref{thm:Hass-Scott-2.7-4.2}(4.2).

\begin{theorem}[non-taut loops have immersed monorbigons]
    \label{thm:taut-no-immonorbigon}
    Consider a closed oriented surface $\Surf$ with a (possibly empty) set of punctures $P\subset \Surf$. 
    
    A loop $\gamma \colon \Sph^1 \looparrowright \Surf_P$ is not taut if and only if it has an singular monorbigon that bounds an immersed disc which avoids at least one region.
\end{theorem}

\begin{proof}
    While the statement in \cite[Theorem 4.2]{Hass-Scott_Intersection-curves-surfaces_1985} only mentions singular monorbigons, its proof actually finds singular monorbigons whose lifts in the universal cover are embedded.
    This was also observed in \cite[Lemma 4.2]{Chang-Mesmay-al_Tightening-curves-local-moves_2018}, where such monorbigons are called \emph{basic}.
    In particular these bound immersed discs, by the projection map of the universal cover.

    The theorem is then clear if $P\neq\emptyset$. We now explain why the proof of \cite[Theorem 4.2]{Hass-Scott_Intersection-curves-surfaces_1985} actually finds such an immersed disc which avoids at least one region, namely that a non-taut loop in a closed surface remains non-taut after placing some puncture. 

    In the case of the closed sphere or the closed torus, their proof yields an embedded monorbigon, which must therefore avoid some region.
    
    In the case of closed surfaces with higher genus, it reduces to the case of punctured surfaces as in the proof of \cite[Theorem 3.5]{Hass-Scott_Intersection-curves-surfaces_1985}: they construct a surface $\Surf_1$ with at least one puncture (as its fundamental group is generated by $2$ elements) which is a subsurface in some cover of $\Surf$ that contains a lift of the loop, and pursue to find an immersed disc inside there.
\end{proof}

We finally reformulate Theorem \ref{thm:taut-no-immonorbigon} in terms of pinning sets and mobidiscs.

\begin{theorem}[pinning mobidiscs]
    \label{thm:pinning-mobidiscs}
    For a loop $\gamma \colon \Sph^1 \looparrowright \Surf$ with regions $\Regions$, a subset of regions $P\in \Parts(\Regions)$ is pinning if and only if it intersects every proper mobidisc $D\in \MoB(\gamma)\setminus \{\Regions\}$.
\end{theorem}

\begin{remark}[pinning ideal]
    \label{rem:non-pinning-ideal}
    Theorem \ref{thm:pinning-mobidiscs} amounts to saying that a subset of regions $N\subset \Regions$ is non-pinning if and only if it is contained in the complement $\Regions \setminus D$ of a mobidisc $D\in \MoB(\gamma)$. 
    Equivalently, in the dual lattice $(\Parts(\Regions), \supset)$, the ideal of non-pinning sets $(\Parts(\Regions)\setminus \mathcal{PI}(\gamma), \supset)$ is generated by the complements of mobidiscs $\{\Regions\setminus D \colon D\in \MoB(\gamma)\}$.
\end{remark}

\begin{remark}[the advantage of proper mobidiscs]
    The fact that one may always find a mobidisc that avoids at least one region has the following consequences.
    
    Topologically, it means that there does not exists a non-taut loop in a closed surface whose only mobidiscs consist of the whole collection of regions, namely that would become taut in the complement of any puncture.
    
    Computationally, it may be useful to compute the pinning ideal, as one may range over the set of regions of loop, and look for mobidiscs of the loop in the surface punctured at that region.
\end{remark}

\begin{remark}[just for loops]
    Theorems \ref{thm:taut-no-immonorbigon} and \ref{thm:pinning-mobidiscs} only hold for loops.
    
    Indeed, \cite[Figure 0.1]{Hass-Scott_Intersection-curves-surfaces_1985} exhibits a multiloop with $2$ strands inside an annulus which is not taut, but that does not bound any singular bigon involving both strands.
    
    After adding a puncture as in Figure \ref{fig:no_singular_monorbigon_non_taut_multiloop}, we obtain a multiloop in a thrice-punctured sphere which is not taut, but now has no singular monorbigons at all. 
    \begin{figure}[h]
        \centering
        \scalebox{0.36}{\definecolor{c010000}{RGB}{1,0,0}

\def \globalscale {1.000000}
\begin{tikzpicture}[y=1cm, x=1cm, yscale=\globalscale,xscale=\globalscale, every node/.append style={scale=\globalscale}, inner sep=0pt, outer sep=0pt]
  \path[draw=red,fill=red,fill opacity=0.3,line cap=butt,line join=miter,line width=0.0265cm,miter limit=4.0] (5.0132, 10.172).. controls (2.5132, 10.172) and (0.0132, 8.172) .. (0.0132, 5.172).. controls (0.0132, 2.672) and (2.5132, 0.172) .. (5.0132, 0.172).. controls (7.5132, 0.172) and (10.0132, 2.672) .. (10.0132, 5.172).. controls (10.0132, 7.672) and (7.5132, 10.172) .. (5.0132, 10.172) -- cycle(5.0132, 7.672).. controls (6.1917, 7.672) and (6.8525, 6.422) .. (7.5132, 5.172).. controls (6.8525, 3.922) and (6.1917, 2.672) .. (5.0132, 2.672).. controls (3.8347, 2.672) and (2.5132, 3.9935) .. (2.5132, 5.172).. controls (2.5132, 6.3505) and (3.8347, 7.672) .. (5.0132, 7.672) -- cycle;

  \path[draw=black,line cap=butt,line join=miter,line width=0.1058cm] (5.0132, 10.172).. controls (2.5132, 10.172) and (0.0132, 8.172) .. (0.0132, 5.172).. controls (0.0132, 2.672) and (2.5132, 0.172) .. (5.0132, 0.172).. controls (7.5132, 0.172) and (10.0132, 2.672) .. (10.0132, 5.172).. controls (10.0132, 7.672) and (7.5132, 10.172) .. (5.0132, 10.172) -- cycle;

  \path[draw=black,line cap=butt,line join=miter,line width=0.1058cm] (5.0132, 7.672).. controls (7.3703, 7.672) and (7.6562, 2.672) .. (10.0132, 2.672).. controls (12.3703, 2.672) and (12.6562, 7.672) .. (15.0132, 7.672).. controls (16.1917, 7.672) and (17.5132, 6.3505) .. (17.5132, 5.172).. controls (17.5132, 3.9935) and (16.1917, 2.672) .. (15.0132, 2.672).. controls (12.6562, 2.672) and (12.3703, 7.672) .. (10.0132, 7.672).. controls (7.6562, 7.672) and (7.3703, 2.672) .. (5.0132, 2.672).. controls (3.8347, 2.672) and (2.5132, 3.9935) .. (2.5132, 5.172).. controls (2.5132, 6.3505) and (3.8347, 7.672) .. (5.0132, 7.672) -- cycle;

  \path[fill=c010000,even odd rule,line cap=round,line width=0.0873cm] (5.0132, 5.172) circle (0.3044cm);

  \path[fill=c010000,even odd rule,line cap=round,line width=0.0873cm] (11.0132, 5.172) circle (0.3044cm);

  \path[fill=c010000,even odd rule,line cap=round,line width=0.0873cm] (15.0132, 5.172) circle (0.3044cm);

  \tikzstyle{every node}=[font=\fontsize{35}{35}\selectfont]
  \node[anchor=south west] (text54) at (4.4802, 5.7647){$p_1$};

  \node[anchor=south west] (text57) at (10.4598, 5.7647){$p_2$};

  \node[anchor=south west] (text59) at (14.4285, 5.7647){$p_3$};

\end{tikzpicture}}
        \caption{This multiloop $\Sph^1 \sqcup \Sph^1 \looparrowright \Sph^2\setminus \{p_1,p_2,p_3\}$ is not taut yet has no singular monorbigons.}
        \label{fig:no_singular_monorbigon_non_taut_multiloop}
    \end{figure}
\end{remark}

\subsection{Reducing pin the loop to boolean formula}
\label{subsec:PIN-to-SAT}

In this subsection, we fix a loop $\gamma\colon \Sph^1 \looparrowright \Surf$ with regions $\Regions$ and double points $V$.

Let us now reformulate the \textsc{LooPinNum} problem as a satisfiability problem.

\begin{definition}[formulae and satisfying assignments]
Consider a set of variables $x=\{x_1,\dots,x_r\}$ which may take boolean values $\{\true,\false\}$.
A \emph{formula} $\phi$ is an expression in the variables $x$ involving negations $\neg$, disjunctions $\wedge$, conjunctions $\vee$ (and parentheses). 

Evaluating the variables $x$ at $\xi \in \{\true,\false\}^r$ yields a value $\phi(\xi)\in \{\true,\false\}$ for the formula $\phi$.
Such an assignment $\xi$ of the variables $x$ is said to \emph{satisfy} $\phi$ when $\phi(\xi)=\true$.

Two formulae $\phi,\phi'$ on the set of variables $x$ are called \emph{logically equivalent} when their evaluations on every assignment $\xi$ of the variables $x$ are equal $\phi(\xi)=\phi'(\xi)$. 
\end{definition}

Every formula is equivalent to a conjunctive normal form.

\begin{definition}[conjunctive normal form]
Consider a set of boolean variables $x=\{x_1,\dots,x_r\}$.

A \emph{clause} is a formula of the form $l_1\vee \dots \vee l_k$ where each $l_j$ is a \emph{literal}, namely a variable $x_t$ or its negation $\neg x_i$.
A \emph{conjunctive normal form} (abbreviated CNF) is a formula which is conjunction $c_1\wedge \dots\wedge c_m$ of clauses $c_j$.
A conjunctive normal form is called \emph{positive} when none of the clauses feature negated variables.
\end{definition}

\begin{definition}[solution ideal]
    Consider a positive conjunctive normal form $\phi$ over a set of boolean variables $x=\{x_1,\dots,x_r\}$.
    
    For an assignment $\xi$ of $x$ satisfying $\phi$, we call $\{x_i \mid \xi_i=\true\}\in \Parts(x)$ a \emph{solution} of $\phi$. 
    
    The solutions of $\phi$ form an ideal $\mathcal{S}(\phi) \subset \Parts(x)$: a sub-poset which is absorbing under union, and it contains the whole set of all elements $\{x_1,\dots,x_r\}$.
\end{definition}

\begin{remark}[hypergraph vertex cover]
    The solutions of a positive conjunctive normal form correspond to the vertex covers of the hypergraph whose vertices are indexed by the variables and hyperedges correspond to the clauses.
\end{remark}

We will only need to consider conjunctive normal forms which are positive.

\begin{definition}[mobidisc formula]
To a loop $\gamma \colon \Sph^1 \looparrowright \Surf$ we associate its \emph{mobidisc formula}: it is the positive conjunctive normal form on the set of boolean variables indexed by its set of regions whose clauses correspond to $\MoB(\gamma)\subset \Parts(\Regions)$.
\end{definition}

\begin{question}
    Can we describe, among all positive conjunctive normal forms, the mobidisc formulae arising from filling loops, and from those in genus $g$?
    
    The next subsection will show that this class is complex even when $g=0$: the task of finding (the cardinal of) an optimal solution is \textsf{NP}-complete.
\end{question}

We may now reformulate Theorem \ref{thm:pinning-mobidiscs} as follows.

\begin{theorem}[pinning a loop = satisfying its mobidisc formula]
\label{thm:LooPin-to-mobidisc-CNF}
The pinning sets of a loop correspond to the solutions of its mobidisc formula.
In particular, the pinning ideal of a loop is isomorphic to the solution ideal of its mobidisc formula.
\end{theorem}

\begin{remark}[solving the mobidisc formula]
Theorem \ref{thm:pinning-mobidiscs} implies that the CNF associated to $\MoB(\gamma)$ and $\MoB(\gamma)\setminus\{\Regions\}$ are logically equivalent. 
In any given case, one may prune the mobidisc formula by removing all clauses which contain another one (which arise topologically as the innermost mobidiscs).

The pinning problem is solved by simplifying such a formula into a minimal DNF: the disjunction of the conjunctions corresponding to the minimal pinning sets. These minimal pinning sets form the smallest basis for the pinning semi-lattice and pinning ideal.
\end{remark}

\begin{remark}[\textsc{SAT}-solver]
    Once the mobidisc formula has been computed, one may send it to a \textsc{SAT}-solver designed to compute optimal solutions to positive conjunctive normal forms.
    
    In practice, the current \textsc{SAT}-solvers for these problems are surprisingly fast.
\end{remark}

We finally consider the case of a loop in the sphere $\gamma\colon \Sph^1 \looparrowright \Sph^2$.
Recall that the subset of proper mobidiscs $\MoB(\gamma)\setminus \{\Regions\}$ consists of those that avoid at least one region.

\begin{algorithm}[computing mobidiscs]
    \label{algo:mobidisc}
    For a loop in the sphere $\gamma\colon \Sph^1 \looparrowright \Sph^2$, we can compute its collection of proper mobidiscs $\MoB(\gamma)\setminus\{\Regions\} \subset \Parts(\Regions)$ in time $O(\Card(\Regions)^5\log(\Card(\Regions)))$.
\end{algorithm}

\begin{proof}
    For each choice of puncture $\infty \in \Regions$, we apply the following algorithm to find the mobidiscs which do not contain $\infty$, namely of the plane loop $\gamma\in \Sph_\infty$.
    
    We first apply algorithm \ref{algo:singular-monorbigons} to find all singular monorbigons of $\gamma$ in time $O(\Card(\Regions)^3)$.
    For each singular monorbigon $K$, we check if $\gamma(K)$ bounds an immersed disc.
    This may be achieved with the algorithm of \cite{Blank_extending-immersions-circle_1967, Poenaru_extension-immersions-codim-1_1995} in time $O(\Card(\Regions)^6)$ (according to \cite{Shor-VanWyk_Detecting-decomposing-self-overlapping-curves_1992}).
    Instead, we may apply the algorithm \cite{Shor-VanWyk_Detecting-decomposing-self-overlapping-curves_1992} in time $O(\Card(\Regions)^3\log(\Card(\Regions)))$ after computing an embedding of $\gamma$ as a planar graph with straight edges in time $O(\Card(\Regions))$ using \cite{Hopcroft-Tarjan_planar-embedding_1974}.

    Finally, for each of these $O(\Card(\Regions)^2)$ immersed monorbigons, we compute the corresponding mobidiscs using winding numbers as in Lemma \ref{lem:unique-monorbidisc}, in time $O(\Card(\Regions))$.
    
    All this yields the desired collection of proper mobidiscs $\MoB(\gamma)\setminus \{\Regions\}$ in the announced time.
\end{proof}

\begin{remark}[in \textsc{NP}]
This recovers Corollary \ref{cor:MultiLooPinNum-is-NP} for the special case of loops in the sphere.
\end{remark}

\begin{example}[a mobidisc formula]
    The mobidisc formula of the \href{https://christopherlloyd.github.io/LooPindex/multiloops/11\%5E1_97.html}{loop $11^1_{97}$} in Figure \ref{fig:mobidisc_example} (which we think of as living in $\Sph^2$) is equivalent, after taking the set of minimal clauses with respect to implication, to:
    \begin{align*}
        & (1\vee 2) \wedge
        (2\vee 3) \wedge 
        (2\vee 10) \wedge 
        (4\vee 6)\wedge
        (6\vee 7) \wedge 
        (6\vee 8) \wedge \\
        & (1\vee 4 \vee 5) \wedge (1\vee 8\vee 9) \wedge
        (3\vee5\vee 7) \wedge (3\vee 8 \vee 11) \wedge
        (4\vee10\vee11) \wedge (7\vee9\vee 10) \wedge \\
        & (1\vee 5\vee7\vee 9) \wedge
        (3\vee4\vee 5\vee 11)  \wedge (8\vee9\vee 10\vee 11)
    \end{align*}
    Its minimal satisfying assignments correspond to optimal pinning sets ($A,B$) and minimal pinning sets ($A,B,a,b,c,d,e,f,g,h,i,j,k$), which we may think of as the conjunctions of the variables labeling regions (for instance $A=2\wedge 4\wedge 7 \wedge 8$), hence the mobidisc formula is equivalent to the disjunction $A \vee B \vee  a \vee b \vee c \vee d \vee e \vee  f \vee g\vee h\vee i\vee j\vee k$.
    
    The mobidiscs giving rise to the clauses $(2\vee 10)$ and $(8\vee 9\vee 10\vee 11)$ are the bigons shaded in orange and blue.
\begin{figure}[H]
    \centering
    \scalebox{0.8}{\input{images/tikz/11_1_97pindata_with_highlighted_bigons.tex}}		
    \caption{The \href{https://christopherlloyd.github.io/LooPindex/multiloops/11\%5E1_97.html}{loop $11^1_{97}$}, its minimal pinning sets, and some highlighted mobidiscs.}
    \label{fig:mobidisc_example}
\end{figure}
\end{example}

\subsection{Reducing planar vertex cover to pin-the-loop}

\label{subsec:PlanarVertexCover-to-LooPinNum}

We now show that the \textsc{LooPinNum} problem is \textsf{NP}-hard by reduction from the following decision problem, which is known to be \textsf{NP}-hard by \cite{Garey-Johnson-Stockmeyer_simplified-NP-complete-graph_1976}[Theorem 2.7].

\begin{definition}[Planar vertex cover]
The \textsc{planar vertex cover} problem has:
\begin{itemize}[align=left, noitemsep]
    \item[Instance:] A planar graph $G=(V,E)$ and a positive integer $k\leq \Card(V)$.
    \item[Question:] Is there $U\subset V$ with $\Card(U)\leq k$ such that for all $\{v_1,v_2\}\in E$, at least one of $v_1$ and $v_2$ belongs to $U$?
\end{itemize}
\end{definition}

\begin{remark}[Plane connected]
\label{rem:plane-connected-vertex-cover}
For our purposes, the \textsc{planar vertex cover} problem has the same complexity as its restriction to connected plane graphs.

On the one hand, the vertex covers of a graph correspond to the disjoint union of the vertex covers of each connected components.

On the other hand, there are linear time algorithms \cite{Hopcroft-Tarjan_planar-embedding_1974} which given a graph, decide whether it is planar and when so compute a plane embedding (encoded as a cyclic order of the edges around each vertex).
\end{remark}

\begin{algorithm}[\textsc{planar vertex cover} reduces to \textsc{LooPinNum}]
\label{algo:planar-vertex-cover-to-LooPinNum}
    The \textsc{planar vertex cover} problem admits a polynomial reduction to the \textsc{LooPinNum} problem.
    
    More precisely, there exists a polynomial time algorithm which to a planar graph $G=(V,E)$ associates a plane loop $\gamma \colon \Sph^1 \looparrowright \Disc$ with $O(\Card(E)^2)$ double points, such that the $k$-vertex-covers of $G$ are in one to one correspondence with the $(6\Card(E)+k)$-pinning-sets of $\gamma$.
\end{algorithm}

\begin{proof}[Proof outline]
Given an instance $(G=(V,E),k\in\N)$ of \textsc{planar vertex cover}, we will construct an instance $(\gamma,6\Card(E)+k)$ of \textsc{LooPinNum}, and show that they yield equivalent problems.
As mentioned in Remark \ref{rem:plane-connected-vertex-cover} we may assume that $G$ is connected and endowed with a plane embedding in $\Ball(0,1)$.

\begin{figure}[h]
    \centering
    \hspace{-0.5cm}
    \scalebox{0.8}{\definecolor{c00cd00}{RGB}{0,205,0}
\definecolor{ce82df3}{RGB}{232,45,243}
\definecolor{cefcd00}{RGB}{239,205,0}
\definecolor{ce92cf3}{RGB}{233,44,243}
\definecolor{c00cdff}{RGB}{0,205,255}
\definecolor{ce72ef3}{RGB}{231,46,243}
\definecolor{cff7700}{RGB}{255,119,0}

\def \globalscale {1.000000}
\begin{tikzpicture}[y=1cm, x=1cm, yscale=\globalscale,xscale=\globalscale, every node/.append style={scale=\globalscale}, inner sep=0pt, outer sep=0pt]
  \begin{scope}[shift={(11.3808, 8.372)}]
    \path[draw=c00cd00,fill opacity=0.9961,even odd rule,line cap=round,line width=0.035cm] (4.0996, -7.4952) -- (3.2143, -7.703) -- (-3.0477, -6.5988) -- (-3.8137, -6.0974);

    \path[draw=c00cd00,fill opacity=0.9961,even odd rule,line cap=round,line width=0.035cm] (-3.8653, -6.4073) -- (-2.9772, -6.1988) -- (3.2848, -7.303) -- (4.0507, -7.8043);

    \path[draw=red,fill opacity=0.9961,even odd rule,line cap=round,line width=0.035cm] (-0.9992, -5.6364) -- (-0.0873, -5.6649) -- (1.2674, -6.2967) -- (1.8774, -6.9792);

    \path[draw=red,fill opacity=0.9961,even odd rule,line cap=round,line width=0.035cm] (2.0054, -6.6932) -- (1.0958, -6.6648) -- (-0.259, -6.033) -- (-0.8691, -5.3504);

    \path[draw=ce82df3,fill opacity=0.9961,even odd rule,line cap=round,line width=0.035cm] (3.4765, -8.3225) -- (3.0465, -7.5208) -- (2.5925, -2.3315) -- (2.8786, -1.4617);

    \path[draw=ce82df3,fill opacity=0.9961,even odd rule,line cap=round,line width=0.035cm] (2.5658, -1.492) -- (2.9971, -2.2961) -- (3.4511, -7.4854) -- (3.1652, -8.3545);

    \path[draw=cefcd00,fill opacity=0.9961,even odd rule,line cap=round,line width=0.035cm] (-0.5753, -6.5356) -- (-0.3288, -5.7186) -- (2.6392, -2.1833) -- (3.4564, -1.7706);

    \path[draw=cefcd00,fill opacity=0.9961,even odd rule,line cap=round,line width=0.035cm] (3.2148, -1.5679) -- (2.9503, -2.4443) -- (-0.0177, -5.9797) -- (-0.7836, -6.3664);

    \path[fill=c00cd00,opacity=0.25,even odd rule,line cap=round,line width=0.035cm] (-3.4788, -6.3166) -- (-3.0477, -6.5988) -- (3.2143, -7.703) -- (3.716, -7.5852) -- (3.2848, -7.303) -- (-2.9772, -6.1988) -- cycle;

    \path[fill=red,opacity=0.25,even odd rule,line cap=round,line width=0.035cm] (1.6108, -6.6809) -- (1.2674, -6.2966) -- (-0.0873, -5.6649) -- (-0.6024, -5.6488) -- (-0.259, -6.033) -- (1.0958, -6.6648) -- cycle;

    \path[fill=cefcd00,opacity=0.25,even odd rule,line cap=round,line width=0.035cm] (3.0992, -1.951) -- (2.6392, -2.1833) -- (-0.3288, -5.7186) -- (-0.4777, -6.2119) -- (-0.0177, -5.9797) -- (2.9503, -2.4443) -- cycle;

    \path[fill=ce92cf3,opacity=0.25,even odd rule,line cap=round,line width=0.035cm] (2.7535, -1.842) -- (2.5925, -2.3315) -- (3.0465, -7.5208) -- (3.2901, -7.9749) -- (3.4511, -7.4854) -- (2.9971, -2.2961) -- cycle;

    \path[fill=c00cdff,opacity=0.25,even odd rule,line cap=round,line width=0.035cm] (3.1827, -2.0421) -- (2.6783, -2.1474) -- (-3.1533, -6.2282) -- (-3.4247, -6.6662) -- (-2.9203, -6.5609) -- (2.9113, -2.4801) -- cycle;

    \path[draw=c00cdff,fill opacity=0.9961,even odd rule,line cap=round,line width=0.035cm] (-3.634, -7.0039) -- (-3.1532, -6.2282) -- (2.6783, -2.1474) -- (3.5745, -1.9604);

    \path[draw=c00cdff,fill opacity=0.9961,even odd rule,line cap=round,line width=0.035cm] (3.3936, -1.7018) -- (2.9113, -2.4801) -- (-2.9203, -6.5609) -- (-3.8166, -6.748);

    \path[draw=c00cd00,fill=c00cd00,fill opacity=0.9961,line cap=round,line join=miter,line width=0.15cm] (-11.2552, -6.3945) -- (-4.9688, -7.503);

    \path[draw=c00cdff,line cap=round,line join=miter,line width=0.15cm] (-11.2552, -6.3946) -- (-5.4236, -2.3138);

    \path[draw=ce72ef3,draw opacity=0.9451,line cap=round,line join=miter,line width=0.15cm] (-5.4236, -2.3137) -- (-4.9688, -7.503);

    \path[draw=cefcd00,line cap=round,line join=miter,line width=0.15cm] (-8.3916, -5.8491) -- (-5.4236, -2.3138);

    \path[draw=red,line cap=round,line join=miter,line width=0.15cm] (-7.0368, -6.4807) -- (-8.3916, -5.8491);

    \tikzstyle{every node}=[font=\fontsize{16}{16}\selectfont]
    \node[line width=0.0265cm,anchor=south west] (text48) at (-9.2026, -7.8131){$G$};

    \node[line width=0.0265cm,anchor=south west] (text49) at (-1.2804, -7.9279){$\gamma$};

    \path[fill=black,even odd rule,line cap=round,line width=0.08cm] (-11.2551, -6.3945) circle (0.1257cm);

    \path[draw=cff7700,fill=cff7700,even odd rule,line cap=round,line width=0.15cm] (-8.3916, -5.8491) circle (0.086cm);

    \path[draw=cff7700,fill=cff7700,even odd rule,line cap=round,line width=0.15cm] (-5.4236, -2.3138) circle (0.086cm);

    \path[draw=cff7700,fill=cff7700,even odd rule,line cap=round,line width=0.15cm] (-4.9688, -7.503) circle (0.086cm);

    \path[draw=cff7700,fill=cff7700,even odd rule,line cap=round,line width=0.15cm] (-0.1733, -5.8491) circle (0.086cm);

    \path[draw=cff7700,fill=cff7700,even odd rule,line cap=round,line width=0.15cm] (2.7947, -2.3138) circle (0.086cm);

    \path[draw=cff7700,fill=cff7700,even odd rule,line cap=round,line width=0.15cm] (3.2495, -7.503) circle (0.086cm);

    \path[fill=black,even odd rule,line cap=round,line width=0.08cm] (-7.0368, -6.4807) circle (0.1257cm);

  \end{scope}

\end{tikzpicture}}		
    \caption{Rough idea of the reduction. To a plane graph $G$ is associated a loop $\gamma$ with a bigon for each edge. A vertex cover yields a choice of pins in the local pictures. The technical part of the construction is to join these bigons to form a single loop with a predictable pinning number.}
    \label{fig:reduction_proof_idea}
\end{figure}

Following the example in Figure \ref{fig:precise-gadget-assembly}, we will associate to every edge $e\in E$ an edge-gadget $g_e$, and connect them all together through the boundary circle $\partial \Ball(0,1)$ to form a loop $\gamma$.
The technical part of the construction is to separate the gadgets adequately so as to ensure to ensure the genericity of $\gamma$ and the equivalence of the two problems.
In particular, we will have a correspondence (by the \eqref{eqn:bigons-meet-vertices} property) between the vertices of $G$ and the maximal nonempty intersections of certain edge-bigons of $\gamma$.
Finally we will prove the equivalence of the two problems. One direction will make extensive use of the Lemma \ref{lem:link(singular-monorbigons)=0} on linking numbers which we derived from the characterization of taut loops in Theorem \ref{thm:Hass-Scott-2.7-4.2}.

It is worth mentioning that the metric detail and combinatorial intricacies involved in the manufacturing of the gadgets will ensure the validity of the proof. \end{proof}


\begin{proof}[Construction of the loop]
Let us now perform the construction in detail.

We work in the Euclidean plane $\R^2$ with distance function $d$. For a subset $X\subset \R^2$ and $\epsilon \in \R_+$ we denote by $\Ball(X,\epsilon)$ the $\epsilon$-neighborhood of $X$.

\begin{definition}[Edge gadget]

Fix an oriented segment $e\subset \Ball(0,1)$ and denote by $L_e$ be the oriented line containing $e$.

For a small $\epsilon>0$, the edge gadget $g_e=g_e(\epsilon)$ is the pair of piecewise-smooth curves depicted in Figure \ref{fig:edge-gadget}. 
It is piecewise-linear except for portions of arcs in $\partial\Ball(0,1)$.

We assume that $\epsilon$ is small enough so that $\partial\Ball(0,1) \setminus\Ball(L_e,3\epsilon\tan(\frac{\epsilon}{2}))$ consists of two disjoint arcs which together with $g_e(\epsilon)$ form a single piecewise-smooth loop. 

The shaded bigon $M_e$ is contained in $\overline{\Ball(e,\epsilon)}$, and we assume that $\epsilon$ is small enough so that $\overline{\Ball(e,3\epsilon)}\subset\Ball(0,1)$. 

To such a gadget is associated the six \emph{forced pins} $\{p_1,p_2,p_3,p_4,q_1,q_2\}$, all of which necessarily exist for the piecewise-smooth curve to be taut after surgery. 

The double points $\{u_i\}_{i=1}^4\cup\{v_i\}_{i=1}^2\cup\{w_i\}_{i=1}^6$ will be referred to in the proof justifying the validity of the reduction.

\begin{figure}[h]
    \centering
    \scalebox{0.12}{\definecolor{c00cd00}{RGB}{0,205,0}
\definecolor{cf50600}{RGB}{245,6,0}
\definecolor{c000600}{RGB}{0,6,0}
\definecolor{c000200}{RGB}{0,2,0}

\def \globalscale {1.000000}
\begin{tikzpicture}[y=1cm, x=1cm, yscale=\globalscale,xscale=\globalscale, every node/.append style={scale=\globalscale}, inner sep=0pt, outer sep=0pt]
  \begin{scope}[shift={(23.05, -3.05)}]
    \path[fill=c00cd00,opacity=0.25,even odd rule,line cap=round,line width=0.35cm] (30.0, 15.1) -- (38.0, 11.1) -- (46.0, 11.1) -- (54.0, 15.1) -- (46.0, 19.1) -- (38.0, 19.1) -- cycle;

    \path[draw=black,line cap=round,line join=miter,line width=0.35cm] (56.0614, 27.1021).. controls (57.0957, 26.4865) and (58.0978, 25.8185) .. (59.0644, 25.1003) -- (67.0, 25.1) -- (75.0, 21.1) -- (79.0, 21.1) -- (79.0, 25.1) -- (75.0, 25.1) -- (67.0, 21.1) -- (63.5802, 21.1037).. controls (66.8658, 17.6539) and (69.4219, 13.5705) .. (71.0916, 9.1015) -- (75.0, 9.1) -- (77.0, 10.1) -- (92.0, 10.1) -- (82.0, 5.1) -- (78.0, 5.1) -- (78.0, 9.1) -- (82.0, 9.1) -- (92.0, 4.1) -- (77.0, 4.1) -- (75.0, 5.1) -- (72.3216, 5.1005).. controls (72.4828, 4.4388) and (72.625, 3.772) .. (72.7478, 3.1016);

    \path[draw=red,fill=red,line cap=butt,line join=miter,line width=0.6cm,->] (38.0, 15.1) -- (46.0, 15.1);

    \path[draw=cf50600,fill=red,even odd rule,line width=0.0cm] (38.0, 15.1) circle (0.5cm);

    \path[draw=cf50600,fill=red,even odd rule,line width=0.0cm] (46.0, 15.1) circle (0.5cm);

    \path[draw=black,line cap=butt,line join=miter,line width=0.1cm,dash pattern=on 0.4cm off 0.4cm] (30.0187, 16.1) -- (30.0, 24.1);

    \path[draw=black,line cap=butt,line join=miter,line width=0.1cm,dash pattern=on 0.4cm off 0.4cm] (15.0, 13.1) -- (7.0, 13.1);

    \path[draw=black,line cap=butt,line join=miter,line width=0.1cm,dash pattern=on 0.4cm off 0.4cm] (15.0, 17.1) -- (9.5, 17.1);

    \path[draw=black,line cap=butt,line join=miter,line width=0.1cm,dash pattern=on 0.4cm off 0.4cm] (75.0001, 19.1103) -- (83.0001, 19.129);

    \path[draw=black,line cap=butt,line join=miter,line width=0.1cm,dash pattern=on 0.4cm off 0.4cm] (75.0, 11.1) -- (83.0, 11.1);

    \path[draw=black,line cap=butt,line join=miter,line width=0.1cm,dash pattern=on 0.4cm off 0.4cm] (87.0, 27.1) -- (-23.0, 27.1);

    \path[draw=black,line cap=butt,line join=miter,line width=0.1cm,dash pattern=on 0.4cm off 0.4cm] (87.0, 3.1) -- (-23.0, 3.1);

    \path[draw=black,line cap=butt,line join=miter,line width=0.1cm,<-] (14.0, 17.1) -- (14.0056, 16.6);

    \path[draw=black,line cap=butt,line join=miter,line width=0.1cm,->] (14.016, 13.6) -- (14.0217, 13.1);

    \path[draw=black,line cap=butt,line join=miter,line width=0.1cm,<-] (82.0, 19.1) -- (82.0028, 16.5994);

    \path[draw=black,line cap=butt,line join=miter,line width=0.1cm,->] (81.9972, 13.6006) -- (82.0, 11.1);

    \path[draw=black,line cap=butt,line join=miter,line width=0.1cm,<-] (-22.2334, 27.1) -- (-22.232, 16.5994);

    \path[draw=black,line cap=butt,line join=miter,line width=0.1cm,->] (-22.2376, 13.6005) -- (-22.2334, 3.1);

    \path[draw=c000600,fill=black,even odd rule,line width=0.0cm] (-7.0, 15.1) circle (0.5cm);

    \path[draw=c000600,fill=black,even odd rule,line width=0.0cm] (5.0, 23.1) circle (0.5cm);

    \path[draw=c000600,fill=black,even odd rule,line width=0.0cm] (-9.0, 7.1) circle (0.5cm);

    \path[draw=c000600,fill=black,even odd rule,line width=0.0cm] (77.0, 23.1) circle (0.5cm);

    \path[draw=c000600,fill=black,even odd rule,line width=0.0cm] (80.0, 7.1) circle (0.5cm);

    \path[draw=c000600,fill=black,even odd rule,line width=0.0cm] (71.961, 15.1) circle (0.5cm);

    \path[draw=black,fill opacity=0.9961,even odd rule,line cap=round,line width=0.35cm] (10.653, 19.1342).. controls (11.2075, 19.8109) and (11.787, 20.4668) .. (12.39, 21.1002) -- (10.0, 21.1) -- (8.0, 20.1) -- (-7.0, 20.1) -- (3.0, 25.1) -- (7.0, 25.1) -- (7.0, 21.1) -- (3.0, 21.1) -- (-7.0, 26.1) -- (8.0, 26.1) -- (10.0, 25.1) -- (16.9087, 25.1) -- (17.0053, 25.1715).. controls (17.9616, 25.8771) and (18.9528, 26.534) .. (19.9753, 27.1396);

    \path[draw=black,fill opacity=0.9961,even odd rule,line cap=round,line width=0.35cm] (3.2281, 3.1144).. controls (3.3503, 3.7802) and (3.4916, 4.4423) .. (3.6518, 5.1) -- (1.0, 5.1) -- (-7.0, 9.1) -- (-11.0, 9.1) -- (-11.0, 5.1) -- (-7.0, 5.1) -- (1.0, 9.1) -- (4.8822, 9.1024) -- (4.8931, 9.1315).. controls (5.143, 9.7984) and (5.4129, 10.4576) .. (5.7026, 11.1083);

    \path[draw=c00cd00,fill opacity=0.9961,even odd rule,line cap=round,line width=0.35cm] (5.7026, 11.1083).. controls (6.0026, 11.7821) and (6.3236, 12.4464) .. (6.6651, 13.1001) -- (3.0, 13.1) -- (-5.0, 17.1) -- (-9.0, 17.1) -- (-21.0, 11.1) -- (22.0, 11.1) -- (38.0, 19.1) -- (46.0, 19.1) -- (62.0, 11.1) -- (74.0, 11.1) -- (74.0, 19.1) -- (62.0, 19.1) -- (46.0, 11.1) -- (38.0, 11.1) -- (22.0, 19.1) -- (-21.0, 19.1) -- (-9.0, 13.1) -- (-5.0, 13.1) -- (3.0, 17.1) -- (9.1217, 17.1229).. controls (9.6193, 17.8263) and (10.1247, 18.4895) .. (10.653, 19.1342);

    \path[draw=black,line cap=butt,line join=miter,line width=0.1cm,dash pattern=on 0.4cm off 0.4cm] (38.0, 16.1) -- (38.0, 24.1);

    \path[draw=black,fill opacity=0.9961,even odd rule,line cap=butt,line width=0.1cm,<->,dash pattern=on 0.4cm off 0.4cm] (24.895, 12.8431)arc(156.095:180.0:5.5695 and -5.5695)arc(180.0:203.9054:5.5695 and -5.5695);

    \path[draw=black,line cap=butt,line join=miter,line width=0.1cm,<-] (37.9764, 23.1796) -- (35.4758, 23.1768);

    \path[draw=black,line cap=butt,line join=miter,line width=0.1cm,->] (32.4769, 23.1824) -- (29.9764, 23.1796);

  \end{scope}
  \begin{scope}[shift={(4.05, -0.05)}]
    \tikzstyle{every node}=[font=\fontsize{72}{72}\selectfont]
    \node[text=black,line width=0.0265cm,anchor=south west] (text14) at (29.9495, 10.8174){$\epsilon\tan(\frac{\epsilon}{2})$};

    \node[text=black,line width=0.0265cm,anchor=south west] (text14-7) at (52.8403, 19.5479){$\epsilon$};

    \node[text=black,line width=0.0265cm,anchor=south west] (text14-2) at (42.1248, 11.5856){$\epsilon$};

    \node[text=black,line width=0.0265cm,anchor=south west] (text14-1) at (96.8817, 10.9813){$2\epsilon\tan(\frac{\epsilon}{2})$};

    \node[text=black,line width=0.0265cm,anchor=south west] (text14-1-1) at (-7.3266, 11.0077){$6\epsilon\tan(\frac{\epsilon}{2})$};

    \node[text=red,line width=0.0265cm,anchor=south west] (text14-7-4) at (60.6206, 13.3495){$e$};

    \node[text=c000200,line width=0.0265cm,anchor=south west] (text14-7-4-2) at (59.2199, 9.0245){$M_e$};

    \node[text=black,line width=0.0265cm,shift={(19.0, -3.0)},anchor=south west] (text14-7-4-2-6) at (2.2671, 22.3782){$p_1$};

    \node[text=black,line width=0.0265cm,anchor=south west] (text14-7-4-2-6-9) at (13.892, 19.3466){$u_1$};

    \node[text=black,line width=0.0265cm,anchor=south west] (text14-7-4-2-6-3) at (10.8252, 3.2724){$p_2$};

    \node[text=black,line width=0.0265cm,anchor=south west] (text14-7-4-2-6-7) at (93.1612, 19.4522){$p_3$};

    \node[text=black,line width=0.0265cm,anchor=south west] (text14-7-4-2-6-7-0) at (85.7944, 19.3927){$u_3$};

    \node[text=black,line width=0.0265cm,shift={(19.0, -3.0)},anchor=south west] (text14-7-4-2-6-7-1) at (80.7969, 6.3782){$p_4$};

    \node[text=black,line width=0.0265cm,anchor=south west] (text14-7-4-2-6-9-6) at (18.192, 3.3129){$u_2$};

    \node[text=black,line width=0.0265cm,anchor=south west] (text14-7-4-2-6-7-1-3) at (107.4441, 3.2724){$u_4$};

    \node[text=black,line width=0.0265cm,anchor=south west] (text14-7-4-2-6-8) at (12.8653, 11.3782){$q_1$};

    \node[text=black,line width=0.0265cm,anchor=south west] (text14-7-4-2-6-9-7) at (29.5929, 14.1954){$w_1$};

    \node[text=black,line width=0.0265cm,anchor=south west] (text14-7-4-2-6-9-7-9) at (20.0402, 11.2569){$w_2$};

    \node[text=black,line width=0.0265cm,anchor=south west] (text14-7-4-2-6-9-7-9-2) at (48.3526, 9.2777){$w_4$};

    \node[text=black,line width=0.0265cm,anchor=south west] (text14-7-4-2-6-9-7-9-6) at (7.8776, 11.1878){$w_3$};

    \node[text=black,line width=0.0265cm,anchor=south west] (text14-7-4-2-6-9-7-9-9) at (72.4422, 9.2828){$w_5$};

    \node[text=black,line width=0.0265cm,shift={(19.0, -3.0)},anchor=south west] (text14-7-4-2-6-9-7-2) at (6.7231, 11.3466){$w_6$};

    \node[text=black,line width=0.0265cm,anchor=south west] (text14-7-4-2-6-9-7-9-9-8) at (87.9821, 11.3927){$q_2$};

    \node[text=black,line width=0.0265cm,anchor=south west] (text1) at (82.2749, 14.1954){$v_1$};

    \node[text=black,line width=0.0265cm,anchor=south west] (text2) at (86.1839, 8.3466){$v_2$};

  \end{scope}

\end{tikzpicture}}
    \caption{The edge gadget in metric details (colored core, black bark).}
    \label{fig:edge-gadget}
\end{figure}
\end{definition}

Now to a connected plane graph $G \hookrightarrow \Ball(0,1)$ we associate a loop $\gamma \colon \Sph^1 \looparrowright \R^2$, by positioning edge-gadgets meticulously and assembling them.
One may check that every step of the construction can be performed in polynomial time.

\begin{description}[align=left]
\item[Embed graph:]
According to \cite{Wagner1936,Fary1948,Stein_convex-maps_1951}, we may construct in polynomial time a plane embedding of $G\hookrightarrow \Ball(0,1)$ with straight edges. 

\item[\hypertarget{edge-configuration}{Edge configuration}:]
We may assume (after a generic perturbation of the vertices) that no edges of $G$ are parallel. Each edge $e$ of $G$ spans a line $L_e$ in $\R^2$, and we rescale the embedding so that all intersections $L_e\pitchfork L_{e'}$ lie inside $\Ball(0,1)$.

\item[Fix labels:] Fix a total order on the vertices $V$ and orient each edge $e\in E$ from its smallest vertex $e^- \in V$ to its largest vertex $e^+\in V$. Each oriented edge $e\in E$ generates an oriented affine line intersecting $\Ball(0,1)$ along an oriented chord $c_e=(c_e^-, c_e^+)$. By the \hyperlink{edge-configuration}{edge configuration}, we have $2\Card(E)$ distinct points $c_e^{\pm}\in  \partial\Ball(0,1)$ on the unit circle.

\item[\hypertarget{deformation-retracts}{Ensure deformation retracts}:] 
Choose $\epsilon>0$ such that $\Ball(G,\epsilon)\subset \Ball(0,1)$, and $\forall e_1,e_2\in E$ with $e_1\ne e_2$ we have $\Ball(c_{e_1},3\epsilon\tan(\frac{\epsilon}{2}))\cap \Ball(c_{e_2},3\epsilon\tan(\frac{\epsilon}{2})) \subset \Ball(0,1)$ (possible by the \hyperlink{edge-configuration}{edge configuration}), and for all $\{e_i\}_{i=1}^k\subset E$:
\begin{equation*} 
    \bigcap_{i=1}^k \overline{\Ball(e_i, \epsilon)} 
    \quad
    \text{retracts by deformation onto} 
    \quad
    \bigcap_{i=1}^k e_{i}.
\end{equation*}

\item[Embed gadgets:] 
For each $e\in E$, embed an edge gadget $g_e(\epsilon)$. 
By \hyperlink{edge-configuration}{edge configuration} and \hyperlink{deformation-retracts}{deformation retracts}, if $e\neq e'$ then $g_e\cap  g_{e'}$ consists of four points inside $\Ball(0,1)$.
For each $g_e$, remove the arc $\partial \Ball(0,1) \cap \Ball(c_e,3\epsilon\tan(\frac{\epsilon}{2}))$ between its endpoints. 
Fix an orientation of the resulting closed loop and call it $\gamma$. 

\item[Ensure generic:] We now argue that $\gamma$ may be assumed to be in a generic position.
\begin{itemize}[align=left]
        \item[Transverse intersections:] Since no two edges are parallel, and slopes of segments of $g_e\cap \Ball(0,1)$ are within angle $\epsilon/2$ of the slope of $e$, we may repeat the construction above after choosing $\epsilon$ small enough so that all intersection points between distinct gadgets are transverse double points. 
        
        \item[Double points:] What remains are possible isolated multiple points. We compute all intersection points between all gadgets in polynomial time, as well as the minimum distance $\epsilon'$ between them, and resolve multiple points locally (inside balls of radius $\epsilon'/2$, say) into sets of transverse double points so as to preserve the key intersection property below.
\end{itemize}
\end{description}
The polynomial construction of $\gamma$ is finished. See Figure \ref{fig:precise-gadget-assembly}.

Note that for all $e\in E$ we have the bigon $M_e$ satisfying $e\subset M_e\subset \overline{\Ball(e,\epsilon)}$ whence by \hyperlink{deformation-retracts}{deformation-retracts}, for all $\{e_i\}_{i=1}^k\subset E$:
\begin{equation}
\label{eqn:bigons-meet-vertices}
\tag{deformation-retract}
\bigcap_{i=1}^k M_{e_i}
\quad 
\text{retracts by deformation onto}
\quad 
\bigcap_{i=1}^k e_{i}.
\end{equation}
Both intersections are nonempty if and only if all the $e_i$ share a common vertex.
\end{proof}

\begin{figure}[h]
    \centering
    \scalebox{0.45}{\input{images/tikz/gadgets_around_graph_2.tex}}
    \caption{Assembling edge gadgets around the graph.}
    \label{fig:precise-gadget-assembly}
\end{figure}

\begin{proof}[Proof of the equivalence]
We now show that the instances $(G=(V,E), k)$ of \textsc{planar vertex cover} and $(\gamma,6\Card(E)+k)$ of \textsc{LooPinNum} yield equivalent problems.

\paragraph{From pinning set to vertex cover.}
Suppose that $\gamma$ has a pinning set of cardinal $6\Card(E)+k$. 
There must be at least $6\Card(E)$ pins ($6$ per gadget) pinning the embedded monorbigon regions outside $\Ball(0,1)$, and at most $k$ pins remain. 
There must also be a pin inside every $M_e$. 
We construct a map from the remaining $\le k$ pins to vertices of $G$ as follows: if the pin is in no $M_e$, then choose any vertex; if the pin is only in one $M_e$, then choose either vertex associated to that bigon; if it is in more than one $M_e$, these intersect in a nonempty region containing a unique vertex of $G$ by the \eqref{eqn:bigons-meet-vertices} property which is the chosen vertex. 
This yields a vertex cover of the graph using at most $k$ vertices.

\paragraph{From vertex cover to pinning set.}
Now $G$ has a vertex cover of cardinal $k$. 
Consider a set of $6\Card(E)+k$ pins $P$, and place $6\Card(E)$ pins in the embedded monogons and bigons outside $\Ball(0,1)$ (that is $6$ per gadget), so that $k$ pins remain. 
For each vertex of the cover, list its incident edges $e_1,\dots, e_m$ and place a pin in the region $\cap_{i=1}^m M_{e_i}$, which is homotopic to a disc by the \eqref{eqn:bigons-meet-vertices} property. 
We claim that $\gamma$ is pinned (note that we did not place a pin at infinity).

By contradiction, assume $\gamma$ is not pinned: it is not taut in the complement of the pins so by Theorem \ref{thm:Hass-Scott-2.7-4.2} it has a singular monogon (with marked point $x$, say) or bigon (with marked points $x,y$, say). 
Let $\alpha$ be the restriction of $\gamma$ to this monorbigon. 

We first argue that $\alpha$ are contained in the colored cores of the gadgets or the strands.
By Lemma \ref{lem:link(singular-monorbigons)=0}, for all pins $p, q$ we have $\lk(\alpha, [p]-[q])=0$. 
Thus $\alpha$ cannot pass through a portion of $\Sph^1$ between two distinct gadgets since otherwise it would link with the consecutive monogons. 
Moreover $\alpha$ cannot have a marked point among the double points $u_i$ of a gadget (as in Figure \ref{fig:edge-gadget}), since otherwise the neighboring pin $p_i$ together with at least one pin $p'_j$ in another monogon region of this or an adjacent gadget would yield $\lk(\alpha, [p_i]-[p'_j])\in\{\pm 1,\pm 2\}$ (this holds regardless of whether $\alpha$ is a monogon or a bigon). 

For the rest of this explanation, we now assign colors to portions of $\gamma$.
Each gadget gets assigned its own color which appears only on its core while its bark remains black, and the arcs of $\Sph^1$ between gadgets also remain black (as illustrated in Figure \ref{fig:edge-gadget}). 
If $\alpha$ is a singular monogon, label the colors of its outgoing edges as $c_1$ and $c_2$. If $\alpha$ is a bigon between marked points $x$ and $y$, label its colors $c_1^x,c_2^x,c_1^y,c_2^y$ where $c_1^x$ joins to $c_1^y$ and $c_2^x$ joins to $c_2^y$. See Figure \ref{fig:singular_possibilities}.

\begin{figure}[h]
	\centering
	\scalebox{0.35}{\input{images/tikz/reduction_supporting_figs_v2.tex}}
	
	\caption{Left: Arguing that $\alpha$ must avoid the portion of $\gamma$ running between gadgets. Middle: Singular monorbigon $\alpha$ with marked points and colors of initial segments. Right: Ruling out the possibility that $\alpha$ is a ``bigon between gadgets''.}
	\label{fig:singular_possibilities}
\end{figure}

Note that if $\alpha$ is a monogon with marked point $x$, then $c_1=c_2$ and are not black. Namely, $x$ is not a point of the form $w_1, w_6$ or $u_1, u_2$. 
This follows from the argument above, since $x$ does not border a gadget monogon, and $\alpha$ does not contain a portion of $\Sph^1$ connecting two distinct gadgets.

Similarly, if $\alpha$ is a bigon, then $c_1^x=c_1^y$, $c_2^x=c_2^y$, and neither $c_1^x$ nor $c_2^y$ is black. 
If $c_1^x\neq c_2^x$, then $\alpha$ is a ``bigon between gadgets'' and the marked points are among the vertices labeled $\{z_1,z_2,z_3,z_4\}$ in the right side of Figure \ref{fig:singular_possibilities}. However, since the colored portions of distinct gadgets intersects exactly $4$ times, the reader may verify that the few possible cases may be ruled out considering the pins $p_1$, $p_2$, and $p_3$ and applying the linking Lemma \ref{lem:link(singular-monorbigons)=0}.

We have reduced to the case where $\alpha$ is contained within a single gadget, monochromatic, and its marked point(s) of $\alpha$ are among $\{w_1,w_2,w_3,w_4,w_5,w_6\}$.
As the reader may check, Lemma \ref{lem:link(singular-monorbigons)=0} rules all out remaining cases and the argument is complete.

This completes the proof of the reduction \ref{algo:planar-vertex-cover-to-LooPinNum}.
\end{proof}

\begin{corollary}[\textsc{LooPinNum} is \textsf{NP}-hard]
\label{cor:LoopPin-NP-Hard}
The \textsc{LooPinNum} problem is \textsf{NP}-hard, even in restriction to loops in the plane.
\end{corollary}

\begin{proof}
The \textsc{planar vertex cover} problem is \textsf{NP}-complete \cite{Garey-Johnson-Stockmeyer_simplified-NP-complete_1974}.
The Algorithm \ref{algo:planar-vertex-cover-to-LooPinNum} returns a loop whose size is proportional to the size of its input planar graph.
Hence the \textsc{LooPinNum} problem is \textsf{NP}-hard.
\end{proof}

\begin{remark}[vertex degrees $\le 3$]
The \textsc{planar vertex cover} problem remains \textsf{NP}-complete even for graphs with $\max$-degree $3$ by \cite{Garey-Johnson_rectilinear-Steiner-tree_1977}[Lemma 1], but such an assumption would not simplify our construction.

In fact, our construction could be adapted without significant modification to work for any graph (planar or not), but the proof would be more involved.
\end{remark}

\section{The pinning ideal of a multiloop}

\label{sec:pinning-ideal}

In this last section, we first recall the definition of the pinning ideals and semi-lattices of a multiloop, and ask a few questions relating their combinatorial structure to the topology of multiloops.
Then we present counterexamples to certain naive conjectures about pinning ideals or semi-lattices discovered by systematic exploration.
Finally we propose a few heuristics for approximation algorithms.

\subsection{Structure of the pinning ideal}

\begin{definition}[pinning ideal]
    Consider a multiloop $\gamma\colon \sqcup_1^s \Sph^1 \looparrowright \Surf$ with regions $\Regions$. 

    The poset of regions $(\Parts(\Regions), \subset)$ forms a semi-lattice with respect to $\cup$. 
    The pinning sets form the \emph{pinning ideal} $\mathcal{PI} \subset \Parts(\Regions)$: a sub-poset which is absorbing under union.
    The unions of minimal pinning sets form the \emph{pinning semi-lattice} $\mathcal{PSL}\subset \mathcal{PI}$. 
\end{definition}

\begin{question}[structure of pinning ideals]
    A general question is: what can we say about the join posets which arise as pinning ideals of filling loops or multiloops?

    For instance, for which $s ,d, u, v \in \N$ can we find a multiloop with $s$ strands and $d$ double points, having $u$ optimal pinning sets and $v$ minimal pinning sets?
\end{question}

By Theorem \ref{thm:LooPin-to-mobidisc-CNF}, the pinning ideal of a loop is isomorphic to the solution ideal of its mobidisc formula.
The pinning ideal of multiloops remains mysterious.

\begin{question}[basis of the pinning ideal]
    Given a multiloop, can one compute a basis for its pinning ideal in polynomial time?
\end{question}

The following aims at qualitative differences between pinning ideals of loops and multiloops.

\begin{question}[counting strands]
    Can we compute or estimate lower bounds for the number of strands of a multiloop from its pinning ideal?
\end{question}

For a filling multiloop $\gamma \colon \sqcup_1^s \Sph^1 \looparrowright \Surf$, consider the increasing functionals on its lattice of regions $\Parts(\Regions)$, which to a set $P\in \Regions$ associate: 
\begin{itemize}[noitemsep]
    \item its cardinal $\Card(P) = \sum_{P_i\in P} 1$, its \emph{total-degree} $\operatorname{Deg}(P) = \sum_{P_i\in P} \deg(P_i)$;
    \item minus the self-intersection number $-\si_P(\gamma)$ of the multicurve $\gamma \colon \sqcup_1^s \Sph^1 \looparrowright \Surf_P$.
\end{itemize}
Note that $\Card(P)$ and $\deg(P)$ are valuations in the sense of lattice theory.

\begin{question}[multiloops with the same self-intersection functional]
    We have many examples of loops in the sphere with isomorphic pinning ideals.
    
    More precisely, we may say that two multiloops $\alpha, \beta$ to be $\si$-equivalent when their self-intersection functionals $\Parts(\Regions(\alpha)) \ni P\mapsto \si_P(\alpha) \in \N$ and $\Parts(\Regions(\beta)) \ni P\mapsto \si_P(\beta) \in \N$ are conjugate by a bijection $\varphi \colon \Regions(\alpha) \to \Regions(\beta)$.
    Can we describe the equivalence classes?
\end{question}

\begin{question}[topological moves] 

How do pinning ideals of multiloops behave under Reidemeister moves, flypes, and crossing resolutions? We will see in Figures \ref{fig:pinnum-R3} and \ref{fig:pin-flype} that the pinning number can change under $R3$-moves and flypes on indecomposible loops in the sphere.

How do pinning ideals of loops behave under the operations of spheric-sums and toric-sums introduced in \cite{CLS_ChorDiaGraFiloop_2023}?
\end{question}

\subsection{Variance under Reidemeister triangle moves and flypes}

One may hope to interpret certain link invariants in terms of certain numerical invariants of multiloops related to the pinning problem (such as the pinning number, the number of optimal or minimal pinning sets, etc.).

For this one must first choose a map from multiloops to links, and show that if two multiloops yield the same link, then they have equal pinning quantities.
The next paragraphs records the failure of two such approaches.

\subsubsection*{Lifting multiloops to legendrian links: pinning varies under $R3$-moves}

A multiloop $\gamma \colon \sqcup_1^s  \Sph^1\looparrowright \Surf$ lifts to an unoriented Legendrian link in the projective tangent bundle $\overline{\gamma}\colon \sqcup_1^s  \Sph^1 \looparrowright \overline{\Surf}$.

Two such links are isotopic if and only if their projections differ by sequences of certain combinatorial moves, including Reidemeister moves $R3$.

If a quantity related to the pinning problem for multiloops were invariant under the Reidemeister move $R3$, then one may hope to relate it to an isotopy invariant of that Legendrian link.
Alas, the pinning number is not, as one can see in Figure \ref{fig:pinnum-R3}. 

\begin{figure}[H]
    \centering
    \scalebox{0.8}{\def \globalscale {1.000000}
\begin{tikzpicture}[y=1cm, x=1cm, yscale=\globalscale,xscale=\globalscale, every node/.append style={scale=\globalscale}, inner sep=0pt, outer sep=0pt]
  \path[draw=black,line cap=butt,line join=miter,line width=0.0794cm] (2.1783, 5.9281) -- (0.0132, 4.6781) -- (0.0132, 2.1781) -- (0.8793, 1.6781) -- (1.3123, 1.9281) -- (1.3123, 3.4281) -- (2.1783, 3.9281) -- (3.0443, 3.4281) -- (3.0443, 2.4281) -- (2.1783, 1.9281) -- (0.8793, 2.6781) -- (0.4462, 2.4281) -- (0.4462, 1.4281) -- (2.6113, 0.1781) -- (4.7764, 1.4281) -- (4.7764, 2.4281) -- (4.3434, 2.6781) -- (3.0443, 1.9281) -- (2.1783, 2.4281) -- (2.1783, 3.4281) -- (3.0443, 3.9281) -- (3.9103, 3.4281) -- (3.9103, 1.9281) -- (4.3434, 1.6781) -- (5.2094, 2.1781) -- (5.2094, 4.6781) -- (3.0443, 5.9281) -- (2.1783, 5.4281) -- (2.1783, 4.9281) -- (3.4773, 4.1781) -- (3.4773, 3.1781) -- (2.6113, 2.6781) -- (1.7453, 3.1781) -- (1.7453, 4.1781) -- (3.0443, 4.9281) -- (3.0443, 5.4281) -- cycle;

  \path[draw=black,line cap=butt,line join=miter,line width=0.0794cm] (8.2405, 5.9281) -- (6.0754, 4.6781) -- (6.0754, 2.1781) -- (6.9414, 1.6781) -- (7.3744, 1.9281) -- (7.3744, 3.4281) -- (8.2405, 3.9281) -- (8.2405, 2.9281) -- (9.1065, 2.4281) -- (8.2405, 1.9281) -- (6.9414, 2.6781) -- (6.5084, 2.4281) -- (6.5084, 1.4281) -- (8.6735, 0.1781) -- (10.8385, 1.4281) -- (10.8385, 2.4281) -- (10.4055, 2.6781) -- (9.1065, 1.9281) -- (8.2405, 2.4281) -- (9.1065, 2.9281) -- (9.1065, 3.9281) -- (9.9725, 3.4281) -- (9.9725, 1.9281) -- (10.4055, 1.6781) -- (11.2716, 2.1781) -- (11.2716, 4.6781) -- (9.1065, 5.9281) -- (8.2405, 5.4281) -- (8.2405, 4.9281) -- (9.5395, 4.1781) -- (9.5395, 3.1781) -- (8.6735, 3.6781) -- (7.8075, 3.1781) -- (7.8075, 4.1781) -- (9.1065, 4.9281) -- (9.1065, 5.4281) -- cycle;

\end{tikzpicture}}
    \caption{The pinning number is not invariant under $R3$ moves: \cite[\href{https://github.com/ChristopherLloyd/LooPin/blob/main/counterexamples/R3_variance/tex/R3.pdf}{link}]{Simon_LooPin-Code_2023}.}
    \label{fig:pinnum-R3}
\end{figure}

\subsubsection*{Lifting multiloops to alternating links: pinning varies under flypes}
A multiloop $\gamma \colon \sqcup_1^s  \Sph^1\looparrowright \R^2$ yields two alternating diagrams for a link $\hat{\gamma} \colon \sqcup_1^s \Sph^1 \hookrightarrow \R^2 \times \R$ which are related by a mirror image reflection.

An oriented link may have several alternating diagrams, by the Tait flyping conjecture proved by Thistlethwaite--Menasco \cite{Menasco-Thistlethwaite_flyping-announcement_1991, Menasco-Thistlethwaite_classification-alternating-links_1993}, the indecomposible alternating diagrams of a prime link are related by sequences of flypes.

Hence if a pinning quantity of indecomposible multiloops in the sphere is invariant under flypes moves, then one may hope to relate it to an isotopy invariant of prime alternating links.
Alas, the pinning number is not, as one can see in Figure \ref{fig:pin-flype}.

\begin{figure}[H]
    \centering
    \scalebox{1.2}{\def \globalscale {1.000000}
\begin{tikzpicture}[y=1cm, x=1cm, yscale=\globalscale,xscale=\globalscale, every node/.append style={scale=\globalscale}, inner sep=0pt, outer sep=0pt]
  \path[draw=black,line cap=butt,line join=miter,line width=0.0529cm] (1.3536, 0.7471).. controls (1.0202, 0.7471) and (0.5893, 0.5114) .. (0.3536, 0.7471).. controls (-0.1179, 1.2185) and (-0.1179, 2.2757) .. (0.3536, 2.7471).. controls (0.5893, 2.9828) and (1.0456, 2.8747) .. (1.3536, 2.7471).. controls (1.5713, 2.6569) and (1.6358, 2.3373) .. (1.8536, 2.2471).. controls (2.3155, 2.0558) and (3.0, 2.6007) .. (3.3536, 2.2471).. controls (3.5893, 2.0114) and (3.5893, 1.4828) .. (3.3536, 1.2471).. controls (3.1179, 1.0114) and (2.6615, 1.3747) .. (2.3536, 1.2471).. controls (2.1358, 1.1569) and (1.9438, 0.9649) .. (1.8536, 0.7471).. controls (1.7898, 0.5932) and (1.7357, 0.365) .. (1.8536, 0.2471).. controls (1.9714, 0.1293) and (2.2357, 0.1293) .. (2.3536, 0.2471).. controls (2.4714, 0.365) and (2.4173, 0.5932) .. (2.3536, 0.7471).. controls (2.2634, 0.9649) and (2.0713, 1.1569) .. (1.8536, 1.2471).. controls (1.5456, 1.3747) and (1.0893, 1.0114) .. (0.8536, 1.2471).. controls (0.6179, 1.4828) and (0.6179, 2.0114) .. (0.8536, 2.2471).. controls (0.9714, 2.365) and (1.1996, 2.1834) .. (1.3536, 2.2471).. controls (1.5713, 2.3373) and (1.6358, 2.6569) .. (1.8536, 2.7471).. controls (2.0075, 2.8109) and (2.1996, 2.8109) .. (2.3536, 2.7471).. controls (2.5713, 2.6569) and (2.7634, 2.4649) .. (2.8536, 2.2471).. controls (2.9173, 2.0932) and (2.9714, 1.865) .. (2.8536, 1.7471).. controls (2.7357, 1.6293) and (2.4714, 1.6293) .. (2.3536, 1.7471).. controls (2.2357, 1.865) and (2.2898, 2.0932) .. (2.3536, 2.2471).. controls (2.4438, 2.4649) and (2.6358, 2.6569) .. (2.8536, 2.7471).. controls (3.1615, 2.8747) and (3.6179, 2.9828) .. (3.8536, 2.7471).. controls (4.325, 2.2757) and (4.325, 1.2185) .. (3.8536, 0.7471).. controls (3.2643, 0.1579) and (2.1869, 0.7471) .. (1.3536, 0.7471) -- cycle;

  \path[draw=black,line cap=butt,line join=miter,line width=0.0529cm] (5.8536, 0.7471).. controls (5.5202, 0.7471) and (5.0893, 0.5114) .. (4.8536, 0.7471).. controls (4.3821, 1.2185) and (4.3821, 2.2757) .. (4.8536, 2.7471).. controls (5.325, 3.2185) and (6.2376, 3.0023) .. (6.8536, 2.7471).. controls (7.0713, 2.6569) and (7.1358, 2.3373) .. (7.3536, 2.2471).. controls (7.5075, 2.1834) and (7.7357, 2.365) .. (7.8536, 2.2471).. controls (8.0893, 2.0114) and (8.0893, 1.4828) .. (7.8536, 1.2471).. controls (7.6179, 1.0114) and (7.1615, 1.3747) .. (6.8536, 1.2471).. controls (6.6358, 1.1569) and (6.4438, 0.9649) .. (6.3536, 0.7471).. controls (6.2898, 0.5932) and (6.2357, 0.365) .. (6.3536, 0.2471).. controls (6.4714, 0.1293) and (6.7357, 0.1293) .. (6.8536, 0.2471).. controls (6.9714, 0.365) and (6.9173, 0.5932) .. (6.8536, 0.7471).. controls (6.7634, 0.9649) and (6.5713, 1.1569) .. (6.3536, 1.2471).. controls (6.0456, 1.3747) and (5.5893, 1.0114) .. (5.3536, 1.2471).. controls (5.1179, 1.4828) and (5.1179, 2.0114) .. (5.3536, 2.2471).. controls (5.4714, 2.365) and (5.6996, 2.1834) .. (5.8536, 2.2471).. controls (6.0713, 2.3373) and (6.2634, 2.5294) .. (6.3536, 2.7471).. controls (6.4173, 2.9011) and (6.4714, 3.1293) .. (6.3536, 3.2471).. controls (6.2357, 3.365) and (5.9714, 3.365) .. (5.8536, 3.2471).. controls (5.7357, 3.1293) and (5.7898, 2.9011) .. (5.8536, 2.7471).. controls (5.9438, 2.5294) and (6.1358, 2.3373) .. (6.3536, 2.2471).. controls (6.5075, 2.1834) and (6.6996, 2.1834) .. (6.8536, 2.2471).. controls (7.0713, 2.3373) and (7.1358, 2.6569) .. (7.3536, 2.7471).. controls (7.6615, 2.8747) and (8.1179, 2.9828) .. (8.3536, 2.7471).. controls (8.825, 2.2757) and (8.825, 1.2185) .. (8.3536, 0.7471).. controls (7.7643, 0.1579) and (6.6869, 0.7471) .. (5.8536, 0.7471) -- cycle;

\end{tikzpicture}}
    \caption{The pinning number may change under flypes \& mutations: \cite[\href{https://github.com/ChristopherLloyd/LooPin/blob/main/counterexamples/flype_and_mutation_variance/tex/flype_and_mutation_variance.pdf}{link}]{Simon_LooPin-Code_2023}.}
    \label{fig:pin-flype}
\end{figure}

\subsection{Databases of pinning ideals and related quantities}

Using Algorithm \ref{algo:Intro:self-intersection} and some functionalities from \texttt{plantri} \cite{plantri}, we computed the pinning ideals and semi-lattices of all irreducible indecomposible spherical multiloops with at most $12$ regions (unoriented and up to reflection), and the statistics of certain numerical parameters.
The number of such loops follows \cite[A264759]{oeis} and the number of multiloops follows \cite[A113201]{oeis}.
The results are available in the \href{https://github.com/ChristopherLloyd/LooPindex}{LooPindex} \cite{Simon-Stucky_LooPindex_2024}.
Let us mention a few observations. 

\begin{remark}[region degrees]
    On average, there is a strong correlation between the region's degrees and probability of being pinned.
    It is false however that the average degree increases from optimal pinning sets to minimal pinning sets. See Figure \ref{fig:counterexample_naive_gonality_conjecture}.
\end{remark}

\begin{figure}[H]
\centering
{\scalefont{1}\scalebox{0.35}{\definecolor{cd82626}{RGB}{216,38,38}

\def \globalscale {1.000000}
\begin{tikzpicture}[font=\bf,y=1cm, x=1cm, yscale=\globalscale,xscale=\globalscale, every node/.append style={scale=\globalscale}, inner sep=0pt, outer sep=0pt]
  \path[draw=cd82626,fill=cd82626] (10.4951, 4.4254) circle (0.0529cm);

  \path[draw=cd82626,fill=cd82626] (10.4951, 2.4834) circle (0.0529cm);

  \path[draw=cd82626,fill=cd82626] (6.6146, 2.4834) circle (0.0529cm);

  \path[draw=cd82626,fill=cd82626] (6.6146, 10.2515) circle (0.0529cm);

  \path[draw=cd82626,fill=cd82626] (2.734, 10.2515) circle (0.0529cm);

  \path[draw=cd82626,fill=cd82626] (2.734, 2.4834) circle (0.0529cm);

  \path[draw=cd82626,fill=cd82626] (4.6743, 2.4834) circle (0.0529cm);

  \path[draw=cd82626,fill=cd82626] (4.6743, 6.3675) circle (0.0529cm);

  \path[draw=cd82626,fill=cd82626] (12.4354, 6.3675) circle (0.0529cm);

  \path[draw=cd82626,fill=cd82626] (12.4354, 0.5413) circle (0.0529cm);

  \path[draw=cd82626,fill=cd82626] (8.5549, 0.5413) circle (0.0529cm);

  \path[draw=cd82626,fill=cd82626] (8.5549, 8.3095) circle (0.0529cm);

  \path[draw=cd82626,fill=cd82626] (0.7938, 8.3095) circle (0.0529cm);

  \path[draw=cd82626,fill=cd82626] (0.7938, 4.4254) circle (0.0529cm);

  \path[draw=cd82626,line width=0.1588cm] (10.4951, 4.4254) -- (10.4951, 2.4834);

  \path[draw=cd82626,line width=0.1588cm] (10.4951, 2.4834) -- (6.6146, 2.4834);

  \path[draw=cd82626,line width=0.1588cm] (6.6146, 2.4834) -- (6.6146, 10.2515);

  \path[draw=cd82626,line width=0.1588cm] (6.6146, 10.2515) -- (2.734, 10.2515);

  \path[draw=cd82626,line width=0.1588cm] (2.734, 10.2515) -- (2.734, 2.4834);

  \path[draw=cd82626,line width=0.1588cm] (2.734, 2.4834) -- (4.6743, 2.4834);

  \path[draw=cd82626,line width=0.1588cm] (4.6743, 2.4834) -- (4.6743, 6.3675);

  \path[draw=cd82626,line width=0.1588cm] (4.6743, 6.3675) -- (12.4354, 6.3675);

  \path[draw=cd82626,line width=0.1588cm] (12.4354, 6.3675) -- (12.4354, 0.5413);

  \path[draw=cd82626,line width=0.1588cm] (12.4354, 0.5413) -- (8.5549, 0.5413);

  \path[draw=cd82626,line width=0.1588cm] (8.5549, 0.5413) -- (8.5549, 8.3095);

  \path[draw=cd82626,line width=0.1588cm] (8.5549, 8.3095) -- (0.7938, 8.3095);

  \path[draw=cd82626,line width=0.1588cm] (0.7938, 8.3095) -- (0.7938, 4.4254);

  \path[draw=cd82626,line width=0.1588cm] (0.7938, 4.4254) -- (10.4951, 4.4254);

\end{tikzpicture}}}
{\scalefont{1}\scalebox{0.35}{\definecolor{cd82626}{RGB}{216,38,38}

\def \globalscale {1.000000}
\begin{tikzpicture}[font=\bf,y=1cm, x=1cm, yscale=\globalscale,xscale=\globalscale, every node/.append style={scale=\globalscale}, inner sep=0pt, outer sep=0pt]
  \path[draw=cd82626,fill=cd82626] (8.5549, 10.257) circle (0.0529cm);

  \path[draw=cd82626,fill=cd82626] (4.6743, 10.257) circle (0.0529cm);

  \path[draw=cd82626,fill=cd82626] (4.6743, 6.3765) circle (0.0529cm);

  \path[draw=cd82626,fill=cd82626] (12.4354, 6.3765) circle (0.0529cm);

  \path[draw=cd82626,fill=cd82626] (12.4354, 0.5556) circle (0.0529cm);

  \path[draw=cd82626,fill=cd82626] (2.734, 0.5556) circle (0.0529cm);

  \path[draw=cd82626,fill=cd82626] (2.734, 12.1973) circle (0.0529cm);

  \path[draw=cd82626,fill=cd82626] (6.6146, 12.1973) circle (0.0529cm);

  \path[draw=cd82626,fill=cd82626] (6.6146, 2.4959) circle (0.0529cm);

  \path[draw=cd82626,fill=cd82626] (10.4951, 2.4959) circle (0.0529cm);

  \path[draw=cd82626,fill=cd82626] (10.4951, 8.3167) circle (0.0529cm);

  \path[draw=cd82626,fill=cd82626] (0.7938, 8.3167) circle (0.0529cm);

  \path[draw=cd82626,fill=cd82626] (0.7938, 4.4362) circle (0.0529cm);

  \path[draw=cd82626,fill=cd82626] (8.5549, 4.4362) circle (0.0529cm);

  \path[draw=cd82626,line width=0.1588cm] (8.5549, 10.257) -- (4.6743, 10.257);

  \path[draw=cd82626,line width=0.1588cm] (4.6743, 10.257) -- (4.6743, 6.3765);

  \path[draw=cd82626,line width=0.1588cm] (4.6743, 6.3765) -- (12.4354, 6.3765);

  \path[draw=cd82626,line width=0.1588cm] (12.4354, 6.3765) -- (12.4354, 0.5556);

  \path[draw=cd82626,line width=0.1588cm] (12.4354, 0.5556) -- (2.734, 0.5556);

  \path[draw=cd82626,line width=0.1588cm] (2.734, 0.5556) -- (2.734, 12.1973);

  \path[draw=cd82626,line width=0.1588cm] (2.734, 12.1973) -- (6.6146, 12.1973);

  \path[draw=cd82626,line width=0.1588cm] (6.6146, 12.1973) -- (6.6146, 2.4959);

  \path[draw=cd82626,line width=0.1588cm] (6.6146, 2.4959) -- (10.4951, 2.4959);

  \path[draw=cd82626,line width=0.1588cm] (10.4951, 2.4959) -- (10.4951, 8.3167);

  \path[draw=cd82626,line width=0.1588cm] (10.4951, 8.3167) -- (0.7938, 8.3167);

  \path[draw=cd82626,line width=0.1588cm] (0.7938, 8.3167) -- (0.7938, 4.4362);

  \path[draw=cd82626,line width=0.1588cm] (0.7938, 4.4362) -- (8.5549, 4.4362);

  \path[draw=cd82626,line width=0.1588cm] (8.5549, 4.4362) -- (8.5549, 10.257);

\end{tikzpicture}}}
{\scalefont{1}\scalebox{0.35}{\definecolor{cd82626}{RGB}{216,38,38}

\def \globalscale {1.000000}
\begin{tikzpicture}[font=\bf,y=1cm, x=1cm, yscale=\globalscale,xscale=\globalscale, every node/.append style={scale=\globalscale}, inner sep=0pt, outer sep=0pt]
  \path[draw=cd82626,fill=cd82626] (5.1594, 4.9213) circle (0.0529cm);

  \path[draw=cd82626,fill=cd82626] (10.9802, 4.9213) circle (0.0529cm);

  \path[draw=cd82626,fill=cd82626] (10.9802, 2.0108) circle (0.0529cm);

  \path[draw=cd82626,fill=cd82626] (9.525, 2.0108) circle (0.0529cm);

  \path[draw=cd82626,fill=cd82626] (9.525, 6.3765) circle (0.0529cm);

  \path[draw=cd82626,fill=cd82626] (3.7042, 6.3765) circle (0.0529cm);

  \path[draw=cd82626,fill=cd82626] (3.7042, 10.7421) circle (0.0529cm);

  \path[draw=cd82626,fill=cd82626] (2.249, 10.7421) circle (0.0529cm);

  \path[draw=cd82626,fill=cd82626] (2.249, 7.8317) circle (0.0529cm);

  \path[draw=cd82626,fill=cd82626] (8.0698, 7.8317) circle (0.0529cm);

  \path[draw=cd82626,fill=cd82626] (8.0698, 0.5556) circle (0.0529cm);

  \path[draw=cd82626,fill=cd82626] (12.4354, 0.5556) circle (0.0529cm);

  \path[draw=cd82626,fill=cd82626] (12.4354, 3.466) circle (0.0529cm);

  \path[draw=cd82626,fill=cd82626] (6.6146, 3.466) circle (0.0529cm);

  \path[draw=cd82626,fill=cd82626] (6.6146, 9.2869) circle (0.0529cm);

  \path[draw=cd82626,fill=cd82626] (0.7938, 9.2869) circle (0.0529cm);

  \path[draw=cd82626,fill=cd82626] (0.7938, 12.1973) circle (0.0529cm);

  \path[draw=cd82626,fill=cd82626] (5.1594, 12.1973) circle (0.0529cm);

  \path[draw=cd82626,line width=0.1588cm] (5.1594, 4.9213) -- (10.9802, 4.9213);

  \path[draw=cd82626,line width=0.1588cm] (10.9802, 4.9213) -- (10.9802, 2.0108);

  \path[draw=cd82626,line width=0.1588cm] (10.9802, 2.0108) -- (9.525, 2.0108);

  \path[draw=cd82626,line width=0.1588cm] (9.525, 2.0108) -- (9.525, 6.3765);

  \path[draw=cd82626,line width=0.1588cm] (9.525, 6.3765) -- (3.7042, 6.3765);

  \path[draw=cd82626,line width=0.1588cm] (3.7042, 6.3765) -- (3.7042, 10.7421);

  \path[draw=cd82626,line width=0.1588cm] (3.7042, 10.7421) -- (2.249, 10.7421);

  \path[draw=cd82626,line width=0.1588cm] (2.249, 10.7421) -- (2.249, 7.8317);

  \path[draw=cd82626,line width=0.1588cm] (2.249, 7.8317) -- (8.0698, 7.8317);

  \path[draw=cd82626,line width=0.1588cm] (8.0698, 7.8317) -- (8.0698, 0.5556);

  \path[draw=cd82626,line width=0.1588cm] (8.0698, 0.5556) -- (12.4354, 0.5556);

  \path[draw=cd82626,line width=0.1588cm] (12.4354, 0.5556) -- (12.4354, 3.466);

  \path[draw=cd82626,line width=0.1588cm] (12.4354, 3.466) -- (6.6146, 3.466);

  \path[draw=cd82626,line width=0.1588cm] (6.6146, 3.466) -- (6.6146, 9.2869);

  \path[draw=cd82626,line width=0.1588cm] (6.6146, 9.2869) -- (0.7938, 9.2869);

  \path[draw=cd82626,line width=0.1588cm] (0.7938, 9.2869) -- (0.7938, 12.1973);

  \path[draw=cd82626,line width=0.1588cm] (0.7938, 12.1973) -- (5.1594, 12.1973);

  \path[draw=cd82626,line width=0.1588cm] (5.1594, 12.1973) -- (5.1594, 4.9213);

\end{tikzpicture}}}
\caption{Left and center: Loops with the the property that the average degree of every optimal pinning set is larger than the average degree of every minimal, suboptimal pinning set. \\
Right: A loop with the property that the sum of degrees of some optimal pinning set is greater than the sum of degrees of some minimal suboptimal pinning set. \\ See \cite[\href{https://github.com/ChristopherLloyd/LooPin/blob/main/counterexamples/falsehood_of_naive_degree_conjectures/tex/falsehood_of_naive_degree_conjectures.pdf}{link}]{Simon_LooPin-Code_2023}.}
\label{fig:counterexample_naive_gonality_conjecture}
\end{figure}

One may turn the previous remark into a heuristic leading to efficient algorithms computing almost optimal pinning sets of most multiloops.

\begin{remark}[heuristics]
    The strategy consisting in solving the boolean formula whose clauses correspond to the regions bounded by embedded monorbigons often yields an almost pinning set.
    
    One may construct arbitrarily complicated loops with fixed pinning number $\varpi \ge 3$ by taking long words in $\Free_{\varpi-1}$ and corresponding geodesics in a $\varpi$-punctured sphere.
    However for long geodesics, the regions to be pinned appear obvious.
\end{remark}

\begin{example}[smallest multiloops with $\deg \ge 3$] 

Using lemma \ref{lem:average-degree-regions} and an exhaustive computation we found all multiloops in the sphere with at most $12$ regions, all of which have degree $\ge 3$ (see Figure \ref{fig:smallest_multiloops} and the front page of the \href{https://christopherlloyd.github.io/LooPindex/}{LooPindex}). Their pinning data is illustrated in the figures that follow. Optimal pinning sets are labeled with capital letters and shades of red, and the other minimal pinning sets are labeled with lowercase letters and shades of green.  For better visibility, we do not plot the entire pinning ideal but the pinning semi-lattice, together with the set of all regions.  The heights of vertices in the poset (and the labels therein) correspond to their cardinals. A lighter edge emphasizes a greater difference between its endpoint's cardinals.
\end{example}

\newpage

\begin{multicols}{2}
\begin{figure}[H]
{\scalefont{1}\scalebox{0.5}{\definecolor{cd82626}{RGB}{216,38,38}
\definecolor{c2626d8}{RGB}{38,38,216}
\definecolor{c26d826}{RGB}{38,216,38}
\definecolor{c990000}{RGB}{153,0,0}
\definecolor{c00cc00}{RGB}{0,204,0}
\definecolor{c009900}{RGB}{0,153,0}
\definecolor{c006600}{RGB}{0,102,0}
\definecolor{lime}{RGB}{0,255,0}

\def \globalscale {1.000000}
\begin{tikzpicture}[font=\bf,y=1cm, x=1cm, yscale=\globalscale,xscale=\globalscale, every node/.append style={scale=\globalscale}, inner sep=0pt, outer sep=0pt]
  \path[draw=cd82626,fill=cd82626] (12.4301, 7.9428) circle (0.0529cm);

  \path[draw=cd82626,fill=cd82626] (12.4301, 0.9578) circle (0.0529cm);

  \path[draw=cd82626,fill=cd82626] (0.7885, 0.9578) circle (0.0529cm);

  \path[draw=cd82626,fill=cd82626] (0.7885, 7.9428) circle (0.0529cm);

  \path[draw=c2626d8,fill=c2626d8] (5.4451, 12.5995) circle (0.0529cm);

  \path[draw=c2626d8,fill=c2626d8] (10.1018, 12.5995) circle (0.0529cm);

  \path[draw=c2626d8,fill=c2626d8] (10.1018, 3.2861) circle (0.0529cm);

  \path[draw=c2626d8,fill=c2626d8] (5.4451, 3.2861) circle (0.0529cm);

  \path[draw=c26d826,fill=c26d826] (7.7735, 5.6145) circle (0.0529cm);

  \path[draw=c26d826,fill=c26d826] (3.1168, 5.6145) circle (0.0529cm);

  \path[draw=c26d826,fill=c26d826] (3.1168, 10.2711) circle (0.0529cm);

  \path[draw=c26d826,fill=c26d826] (7.7735, 10.2711) circle (0.0529cm);

  \path[draw=cd82626,line width=0.1588cm] (12.4301, 7.9428) -- (12.4301, 0.9578);

  \path[draw=cd82626,line width=0.1588cm] (12.4301, 0.9578) -- (0.7885, 0.9578);

  \path[draw=cd82626,line width=0.1588cm] (0.7885, 0.9578) -- (0.7885, 7.9428);

  \path[draw=cd82626,line width=0.1588cm] (0.7885, 7.9428) -- (12.4301, 7.9428);

  \path[draw=c2626d8,line width=0.1588cm] (5.4451, 12.5995) -- (10.1018, 12.5995);

  \path[draw=c2626d8,line width=0.1588cm] (10.1018, 12.5995) -- (10.1018, 3.2861);

  \path[draw=c2626d8,line width=0.1588cm] (10.1018, 3.2861) -- (5.4451, 3.2861);

  \path[draw=c2626d8,line width=0.1588cm] (5.4451, 3.2861) -- (5.4451, 12.5995);

  \path[draw=c26d826,line width=0.1588cm] (7.7735, 5.6145) -- (3.1168, 5.6145);

  \path[draw=c26d826,line width=0.1588cm] (3.1168, 5.6145) -- (3.1168, 10.2711);

  \path[draw=c26d826,line width=0.1588cm] (3.1168, 10.2711) -- (7.7735, 10.2711);

  \path[draw=c26d826,line width=0.1588cm] (7.7735, 10.2711) -- (7.7735, 5.6145);

  \path[draw=c990000,fill=c990000,line width=0.0cm] (1.11, 0.6363) ellipse (0.1663cm and 0.1663cm);

  \node[text=white,anchor=south,line width=0.0cm] (text6950) at (1.11, 0.5027){B};

  \path[draw=c00cc00,fill=c00cc00,line width=0.0cm] (1.4759, 0.6363) ellipse (0.1663cm and 0.1663cm);

  \node[text=white,anchor=south,line width=0.0cm] (text1010) at (1.4759, 0.5027){b};

  \path[draw=c009900,fill=c009900,line width=0.0cm] (1.8418, 0.6363) ellipse (0.1663cm and 0.1663cm);

  \node[text=white,anchor=south,line width=0.0cm] (text8923) at (1.8418, 0.5027){c};

  \path[draw=c006600,fill=c006600,line width=0.0cm] (2.2076, 0.6363) ellipse (0.1663cm and 0.1663cm);

  \node[text=white,anchor=south,line width=0.0cm] (text5139) at (2.2076, 0.5027){d};

  \path[draw=red,fill=red,line width=0.0cm] (1.11, 7.6213) ellipse (0.1663cm and 0.1663cm);

  \node[text=white,anchor=south,line width=0.0cm] (text4408) at (1.11, 7.4877){A};

  \path[draw=lime,fill=lime,line width=0.0cm] (1.4759, 7.6213) ellipse (0.1663cm and 0.1663cm);

  \node[text=white,anchor=south,line width=0.0cm] (text5128) at (1.4759, 7.4877){a};

  \path[draw=c00cc00,fill=c00cc00,line width=0.0cm] (1.8418, 7.6213) ellipse (0.1663cm and 0.1663cm);

  \node[text=white,anchor=south,line width=0.0cm] (text2787) at (1.8418, 7.4877){b};

  \path[draw=c006600,fill=c006600,line width=0.0cm] (2.2076, 7.6213) ellipse (0.1663cm and 0.1663cm);

  \node[text=white,anchor=south,line width=0.0cm] (text1287) at (2.2076, 7.4877){d};

  \path[draw=red,fill=red,line width=0.0cm] (3.4383, 9.9496) circle (0.1663cm);

  \node[text=white,anchor=south,line width=0.0cm] (text9510) at (3.4383, 9.816){A};

  \path[draw=lime,fill=lime,line width=0.0cm] (3.8042, 9.9496) ellipse (0.1663cm and 0.1663cm);

  \node[text=white,anchor=south,line width=0.0cm] (text1895) at (3.8042, 9.816){a};

  \path[draw=c00cc00,fill=c00cc00,line width=0.0cm] (4.1701, 9.9496) ellipse (0.1663cm and 0.1663cm);

  \node[text=white,anchor=south,line width=0.0cm] (text9462) at (4.1701, 9.816){b};

  \path[draw=c009900,fill=c009900,line width=0.0cm] (4.536, 9.9496) ellipse (0.1663cm and 0.1663cm);

  \node[text=white,anchor=south,line width=0.0cm] (text7499) at (4.536, 9.816){c};

  \path[draw=c990000,fill=c990000,line width=0.0cm] (3.4383, 7.6213) ellipse (0.1663cm and 0.1663cm);

  \node[text=white,anchor=south,line width=0.0cm] (text1955) at (3.4383, 7.4877){B};

  \path[draw=lime,fill=lime,line width=0.0cm] (3.8042, 7.6213) circle (0.1663cm);

  \node[text=white,anchor=south,line width=0.0cm] (text3414) at (3.8042, 7.4877){a};

  \path[draw=c009900,fill=c009900,line width=0.0cm] (4.1701, 7.6213) circle (0.1663cm);

  \node[text=white,anchor=south,line width=0.0cm] (text7154) at (4.1701, 7.4877){c};

  \path[draw=c006600,fill=c006600,line width=0.0cm] (4.536, 7.6213) circle (0.1663cm);

  \node[text=white,anchor=south,line width=0.0cm] (text9651) at (4.536, 7.4877){d};

  \path[draw=c990000,fill=c990000,line width=0.0cm] (5.7667, 9.9496) ellipse (0.1663cm and 0.1663cm);

  \node[text=white,anchor=south,line width=0.0cm] (text9885) at (5.7667, 9.816){B};

  \path[draw=lime,fill=lime,line width=0.0cm] (6.1325, 9.9496) ellipse (0.1663cm and 0.1663cm);

  \node[text=white,anchor=south,line width=0.0cm] (text3211) at (6.1325, 9.816){a};

  \path[draw=c00cc00,fill=c00cc00,line width=0.0cm] (6.4984, 9.9496) ellipse (0.1663cm and 0.1663cm);

  \node[text=white,anchor=south,line width=0.0cm] (text9723) at (6.4984, 9.816){b};

  \path[draw=c006600,fill=c006600,line width=0.0cm] (6.8643, 9.9496) ellipse (0.1663cm and 0.1663cm);

  \node[text=white,anchor=south,line width=0.0cm] (text5266) at (6.8643, 9.816){d};

  \path[draw=red,fill=red,line width=0.0cm] (5.7667, 7.6213) circle (0.1663cm);

  \node[text=white,anchor=south,line width=0.0cm] (text1917) at (5.7667, 7.4877){A};

  \path[draw=c00cc00,fill=c00cc00,line width=0.0cm] (6.1325, 7.6213) circle (0.1663cm);

  \node[text=white,anchor=south,line width=0.0cm] (text6364) at (6.1325, 7.4877){b};

  \path[draw=c009900,fill=c009900,line width=0.0cm] (6.4984, 7.6213) circle (0.1663cm);

  \node[text=white,anchor=south,line width=0.0cm] (text1793) at (6.4984, 7.4877){c};

  \path[draw=c006600,fill=c006600,line width=0.0cm] (6.8643, 7.6213) circle (0.1663cm);

  \node[text=white,anchor=south,line width=0.0cm] (text3265) at (6.8643, 7.4877){d};

  \path[draw=red,fill=red,line width=0.0cm] (5.7667, 12.2779) ellipse (0.1663cm and 0.1663cm);

  \node[text=white,anchor=south,line width=0.0cm] (text9265) at (5.7667, 12.1444){A};

  \path[draw=lime,fill=lime,line width=0.0cm] (6.1325, 12.2779) ellipse (0.1663cm and 0.1663cm);

  \node[text=white,anchor=south,line width=0.0cm] (text9620) at (6.1325, 12.1444){a};

  \path[draw=c009900,fill=c009900,line width=0.0cm] (6.4984, 12.2779) ellipse (0.1663cm and 0.1663cm);

  \node[text=white,anchor=south,line width=0.0cm] (text5940) at (6.4984, 12.1444){c};

  \path[draw=c006600,fill=c006600,line width=0.0cm] (6.8643, 12.2779) ellipse (0.1663cm and 0.1663cm);

  \node[text=white,anchor=south,line width=0.0cm] (text1972) at (6.8643, 12.1444){d};

  \path[draw=c990000,fill=c990000,line width=0.0cm] (5.7667, 5.2929) circle (0.1663cm);

  \node[text=white,anchor=south,line width=0.0cm] (text4275) at (5.7667, 5.1594){B};

  \path[draw=lime,fill=lime,line width=0.0cm] (6.1325, 5.2929) circle (0.1663cm);

  \node[text=white,anchor=south,line width=0.0cm] (text6735) at (6.1325, 5.1594){a};

  \path[draw=c00cc00,fill=c00cc00,line width=0.0cm] (6.4984, 5.2929) circle (0.1663cm);

  \node[text=white,anchor=south,line width=0.0cm] (text1497) at (6.4984, 5.1594){b};

  \path[draw=c009900,fill=c009900,line width=0.0cm] (6.8643, 5.2929) circle (0.1663cm);

  \node[text=white,anchor=south,line width=0.0cm] (text354) at (6.8643, 5.1594){c};

\end{tikzpicture}}}
\end{figure}
\columnbreak

\begin{figure}[H]
\hspace{-1cm}
{\scalefont{1}\scalebox{0.75}{\input{images/tikz/8_3_2lattice.tex}}}
\end{figure}
\end{multicols}

\vspace{-1cm}

\begin{multicols}{2}
\begin{figure}[H]
{\scalefont{1}\scalebox{0.5}{\input{images/tikz/10_1_18pindata_clean.tex}}}
\end{figure}
\columnbreak

\begin{figure}[H]
\hspace{-2cm}
{\scalefont{1}\scalebox{0.75}{\input{images/tikz/10_1_18lattice.tex}}}
\end{figure}
\end{multicols}

\vspace{-1cm}

\begin{multicols}{2}
\begin{figure}[H]
\hspace{-1cm}
{\scalefont{0.7}\scalebox{0.5}{\input{images/tikz/11_1_97pindata_clean.tex}}}
\end{figure}
\columnbreak

\begin{figure}[H]
\hspace{-3cm}
{\scalefont{1}\scalebox{0.75}{\input{images/tikz/11_1_97lattice.tex}}}
\end{figure}
\end{multicols}

\begin{multicols}{2}
\begin{figure}[H]
\hspace{-1cm}
{\scalefont{0.7}\scalebox{0.5}{\input{images/tikz/12_1_262pindata_clean.tex}}}
\end{figure}
\columnbreak

\begin{figure}[H]
\hspace{-3.5cm}
{\scalefont{1}\scalebox{0.75}{\input{images/tikz/12_1_262lattice.tex}}}
\end{figure}
\end{multicols}

\vspace{-1cm}

\begin{multicols}{2}
\begin{figure}[H]
\hspace{-1cm}
{\scalefont{0.7}\scalebox{0.5}{\input{images/tikz/12_2_301pindata_clean.tex}}}
\end{figure}
\columnbreak

\begin{figure}[H]
\hspace{-3.5cm}
{\scalefont{1}\scalebox{0.75}{\input{images/tikz/12_2_301lattice.tex}}}
\end{figure}
\end{multicols}

\vspace{-1cm}

\begin{multicols}{2}
\begin{figure}[H]
\hspace{-1cm}
{\scalefont{0.7}\scalebox{0.5}{\input{images/tikz/12_4_11pindata_clean.tex}}}
\end{figure}
\columnbreak

\begin{figure}[H]
\hspace{-3cm}
{\scalefont{1}\scalebox{0.75}{\input{images/tikz/12_4_11lattice.tex}}}
\end{figure}
\end{multicols}


\newpage
\section*{Acknowledgments}

We wish to thank Yago Antol\'{i}n and Alan Reid for the invitation and financial support to attend the workshop on orderings and groups at ICMAT in Madrid during summer 2023, where our friendship and collaboration began.

We also recognize the help of Nathan Dunfield who gave us tips for adapting \texttt{SnapPy}  \cite{SnapPy} to our purposes (used here for many of the figures), as well as Gunnar Brinkmann and Brendan McKay who gave us tips and custom plugins for \texttt{plantri}. We thank Matt Clay for help with Figure \ref{fig:lifts_intersecting_domain}.

The second author gratefully acknowledges the unwavering support of his partner Madeleine. 
He is also grateful to the faculty and staff of the Department of Mathematical Sciences at Northern Illinois University for welcoming him as a visiting scholar during the Spring of 2024, where a large portion of this work was completed.

Last but not least, we are pleased to thank the referee for their thorough reading, pointing out several mistakes in our proofs and providing relevant critiques.

\section*{Declarations}

\paragraph{Ethical approval:}
This declaration is not applicable.
\paragraph{Funding:}
This declaration is not applicable.

\paragraph{Availability of data and materials:}

The experimental results \cite{Simon_LooPin-Code_2023} cited throughout this work (including implementations of our algorithms, catalogs of multiloops and their pinning ideals, and counterexamples to various naive conjectures) are available at our \href{https://github.com/ChristopherLloyd/LooPin}{GitHub repository}.

\bibliographystyle{alpha}
\bibliography{biblio.bib}

\end{document}